\documentclass[11pt]{amsart} 

\usepackage{amssymb,amsmath,amscd,amsxtra,diagrams}
\usepackage[margin=4cm]{geometry}
\usepackage[all]{xy}
\usepackage[pdftex]{graphicx,color} 

\newtheorem{theorem}{Theorem}[section]
\newtheorem{lemma}[theorem]{Lemma}

\newtheorem{proposition}[theorem]{Proposition}  
\newtheorem{corollary}[theorem]{Corollary} 

\theoremstyle{definition}
\newtheorem{remark}[theorem]{Remark}

\newtheorem{definition}[theorem]{Definition}

\newcommand{\beq}{\begin{equation}}
\newcommand{\eeq}{\end{equation}}
\newcommand{\beqa}{\begin{eqnarray}}
\newcommand{\eeqa}{\end{eqnarray}}
\newcommand{\beaa}{\begin{eqnarray*}}
\newcommand{\ben}{\begin{eqnarray*}}
\newcommand{\eaa}{\end{eqnarray*}}
\newcommand{\een}{\end{eqnarray*}}

\newcommand{\CC}{\mathbb{C}}

\newcommand{\F}{\mathcal{F}}
\newcommand{\PP}{\mathbb{P}}
\newcommand{\M}{\mathcal{M}}
\newcommand{\W}{\mathcal{W}}

\newcommand{\D}{\mathcal{D}}

\renewcommand{\H}{\mathcal{H}}
\renewcommand{\L}{\mathcal{L}}
\renewcommand{\O}{\mathcal{O}}
\newcommand{\RR}{\mathbb{R}}
\newcommand{\T}{\mathcal{T}}
\newcommand{\U}{\mathcal{U}}
\newcommand{\ZZ}{\mathbb{Z}}
\newcommand{\HH}{\mathbb{H}}

\newcommand{\ff}{\mathbf{f}}
\newcommand{\pp}{\mathbf{p}}
\newcommand{\qq}{\mathbf{q}}
\renewcommand{\tt}{\mathbf{t}}

\def\diag{\mathop{\rm diag} \nolimits}

\begin{document} 

\title[Quantum cohomology of $\mathbb{P}^2$]
{\textbf{ The Period map for quantum cohomology of $\mathbb{P}^2$}}
\author{Todor Milanov}
\address{Kavli IPMU (WPI), UTIAS, The University of Tokyo, Kashiwa, Chiba 277-8583, Japan}
\email{todor.milanov@ipmu.jp} 

\begin{abstract}
We invert the period map defined by the second structure connection of
quantum cohomology of $\mathbb{P}^2$. For small quantum cohomology the
inverse is given explicitly in terms of the Eisenstein series $E_4$
and $E_6$, while for big quantum cohomology the inverse is determined
perturbatively as a Taylor series expansion whose coefficients are
quasi-modular forms. 
\end{abstract}
\maketitle

\section{Introduction}
The results of this paper are in the settings of quantum cohomology of
$\PP^2$. Nevertheless, the problems that we solve can be given in much
more general settings. Let us start by giving the general picture and
providing some background and motivation for our results. 
 
\subsection{The second structure connection}\label{sec:2nd-str-conn}
We assume that the reader is familiar with the definition of a semi-simple
Frobenius manifold (see \cite{Du} for some background). 
Let $M$ be a complex semi-simple Frobenius manifold and let $\T_M$ be
the sheaf of holomorphic vector fields on $M$. By definition the data
of Frobenius structure is given by the following list of objects
\begin{itemize}
\item[(1)]
A non-degenerate symmetric bi-linear pairing  
$(\ ,\ )$ on $\T_M$. 
\item[(2)]
A commutative associative multiplication $\bullet :\T_M\otimes \T_M\to
\T_M$ .
\item[(3)]
A flat vector field ${\mathbf 1}\in \Gamma(M,\T_M)$ that is a unity,
i.e., $\mathbf{1}\bullet v=v$ for all $v\in \T_M.$ 
\item[(4)]
An Euler vector field $E\in \Gamma(M,\T_M)$.
\end{itemize}
We are going to work only with Frobenius manifolds satisfying the
following 4 additional conditions:
\begin{enumerate}
\item[(i)] The tangent bundle $TM$ is trivial and it admits a
  trivialization given by a frame of global flat vector fields.
\item[(ii)]
Recall that the operator 
$$
\operatorname{ad}_E:\T_M\to \T_M,\quad v\mapsto [E,v]
$$
preserves the space of flat vector fields. We require that the
restriction of $\operatorname{ad}_E$ to the space of flat vector
fields is a diagonalizable operator with rational eigenvalues.
\item[(iii)]
  The Frobenius manifold has a {\em calibration} (see Section \ref{sec:pv}).
\item[(iv)]
  The Frobenius manifold has a direct product decomposition
  $M=\CC\times B$ such that if we denote by $t_1:M\to
  \CC$ the projection along $B$, then $dt_1$ is a flat 1-form and
  $\langle dt_1,\mathbf{1}\rangle=1.$ 
\end{enumerate}
Conditions (i)--(iv) are satisfied for all Frobenius manifolds
constructed by quantum cohomology or by the primitive forms in
singularity theory. 

Let us fix a base point $t^\circ\in M$ and a basis
$\{\phi_i\}_{i=1}^N$ of the reference tangent space
$H:=T_{t^\circ}M$. Furthermore, let $(t_1,\dots,t_N)$ be a local flat
coordinate system on an open neighborhood of $t^\circ$ such that $\partial/\partial t_i=\phi_i$ in
$H$. The flat vector fields $\partial/\partial t_i$
($1\leq i\leq N$) 
extend to global flat vector fields on $M$ and provide a
trivialization of the tangent bundle $TM\cong M\times H$. This allows us to
identify the Frobenius multiplication $\bullet$ with a family of
associative  commutative multiplications $\bullet_t :H\otimes H\to H$
depending analytically on $t\in M$. Modifying our choice of 
$\{\phi_i\}_{i=1}^N$ and $\{t_i\}_{i=1}^N$ if necessary we may arrange
that 
\ben
E=\sum_{i=1}^N ((1-d_i) t_i + r_i)\partial/\partial t_i,
\een
where $\partial/\partial t_1$ coincides with the unit vector field
$\mathbf{1}$ and the numbers 
\ben
0=d_1\leq d_2\leq \cdots \leq d_N=:D
\een
are symmetric with respect to the middle of the interval $[0,D]$. The
number $D$ is known as the {\em conformal dimension} of $M$. The
operator 
\ben
\theta:\T_M\to \T_M,\quad v\mapsto [E,v]-\frac{1}{2}(2-D) v
\een
preserves the subspace of flat vector fields. It induces a linear
operator on $H$ which is known to be skew symmetric with respect to
the Frobenius pairing $(\ ,\ )$. Following Givental, we refer to
$\theta$ as the {\em Hodge grading operator}.  

There are two flat connections that one can associate with the
Frobenius structure. The first one is usually called {\em Dubrovin's
  connection}. It is a connection on the $H$-trivial bundle on $M\times
\CC^*$ defined by 
\ben
\nabla_{\partial/\partial t_i} & = &  \frac{\partial}{\partial t_i} -
z^{-1} \phi_i\bullet \\
\nabla_{\partial/\partial z} & = & \frac{\partial}{\partial z} +z^{-1}
\theta -z^{-2} E\bullet
\een
where $z$ is the standard coordinate on $\CC^*=\CC-\{0\}$ and
for $v\in \Gamma(M,\T_M)$ we denote by $v\bullet :H\to H$ the linear
operator of Frobenius multiplication 
by $v$.  

Our main interest is in the 2nd structure connection 
\ben
\nabla^{(n)}_{\partial/\partial t_i} & = &
\frac{\partial}{\partial t_i} + (\lambda-E\bullet_t)^{-1} (\phi_i \bullet_t) (\theta-n-1/2) \\
\nabla^{(n)}_{\partial/\partial\lambda} & = & 
\frac{\partial}{\partial \lambda}-(\lambda-E\bullet_t)^{-1} (\theta-n-1/2),
\een
where $n\in \CC$ is a complex parameter. 
This is a connection on the trivial bundle
\ben
(M\times \CC)'\times H \to (M\times \CC)',
\een
where 
\ben
(M\times \CC)'=\{ (t,\lambda)\ |\ \det (\lambda-E\bullet_t)\neq 0\}.
\een
The hypersurface $\det (\lambda-E\bullet_t)=0$ in $M\times \CC$ is
called the {\em discriminant}.

\subsection{Period vectors}\label{sec:pv}
The definition of the period map depends on the choice of a {\em
  calibration} $S(t,z)$ of $M$. By definition (see \cite{G2}), the
calibration is an operator series $S=1+\sum_{k=1}^\infty
S_k(t)z^{-k}$, $S_k\in \operatorname{End}(H)$, such that
the Dubrovin's connection has a fundamental solution near $z=\infty$
of the form 
\ben
S(t,z) z^{\theta} z^{-\rho},
\een
where $\rho\in \operatorname{End}(H)$ is a nilpotent operator, $[\theta,\rho]=-\rho$, and 
the following symplectic condition holds
\ben
S(t,z)S(t,-z)^T=1,
\een
where ${}^T$ denotes transposition with respect to the Frobenius pairing.

Let us fix a reference point $(t^\circ,\lambda^\circ)\in (M\times
\CC)'$ such that $\lambda^\circ$ is a sufficiently
large real number. It is
easy to check that the following functions provide a fundamental
solution to the 2nd structure connection
\ben
I^{(n)}(t,\lambda) = \sum_{k=0}^\infty (-1)^k S_k(t) \widetilde{I}^{(n+k)}(\lambda),
\een 
where
\ben
\widetilde{I}^{(m)}(\lambda) = e^{-\rho \partial_\lambda \partial_m}
\Big(
\frac{\lambda^{\theta-m-\frac{1}{2}} }{ \Gamma(\theta-m+\frac{1}{2}) }
\Big).
\een
The 2nd structure connection has a
Fuchsian singularity at infinity, therefore the series $I^{(n)}(t,\lambda)$ is
convergent for all $(t,\lambda)$ sufficiently close to
$(t^\circ,\lambda^\circ)$. Using the differential equations we extend
$I^{(n)}$ to a multi-valued analytic function on $(M\times \CC)'$. We
define the following multi-valued functions taking values in $H$:
\ben
I^{(n)}_a(t,\lambda):=I^{(n)}(t,\lambda)\, a, 
\quad a\in H,\quad n\in \ZZ.
\een 
These functions will be called {\em period vectors}. Using analytic continuation
we get a representation 
\beq\label{mon-repr}
\pi_1((M\times\CC)',(t^\circ,\lambda^\circ) )\to \operatorname{GL}(H)
\eeq
called the {\em monodromy representation} of the Frobenius
manifold. The image $W$ of the monodromy representation is called the
{\em monodromy group}. 

Under the semi-simplicity
assumption, we may choose a generic reference point $t^\circ$ on $M$, such that 
the Frobenius multiplication $\bullet_{t^\circ}$ is semi-simple and the
operator $E\bullet_{t^\circ}$ has $N$ pairwise different eigenvalues
$u_i^\circ$ ($1\leq i\leq N$). 
The fundamental group $\pi_1((M\times \CC)', (t^\circ,\lambda^\circ))$ fits into the following
exact sequence
\beq\label{fund-es}
\pi_1(F^\circ,\lambda^\circ) \rTo^{i_*}
\pi_1((M\times \CC)', (t^\circ,\lambda^\circ)  )\rTo^{p_*} 
\pi_1(M,t^\circ)\rTo 1,
\eeq
where $p: (M\times \CC)'\to M$ is the projection on $M$, 
$F^\circ=p^{-1}(t^\circ)=\CC\setminus{\{u_1^\circ,\dots,u_N^\circ\}}$ is the fiber over $t^\circ$, and
$i:F^\circ\to (M\times \CC)'$ is the natural inclusion.  For a proof we refer to
\cite{Shi}, Proposition 5.6.4 or \cite{No}, Lemma 1.5 C.  
Using the exact sequence \eqref{fund-es} we get that the monodromy
group $W$ is generated by the monodromy transformations 
representing the lifts of the generators of $\pi_1(M,t^\circ)$ in
$\pi_1((M\times \CC)', (t^\circ, \lambda^\circ))$ and the
generators of $\pi_1(F^\circ,\lambda^\circ)$. 

The image of $\pi_1(F^\circ,\lambda^\circ)$ under the monodromy
representation is a reflection group that can be described as
follows. Using the differential equations of the 2nd structure connection it is
easy to prove that the pairing
\beq\label{inters-pairing}
(a|b):=(I^{(0)}_a(t,\lambda),(\lambda-E\bullet)I^{(0)}_b(t,\lambda))
\eeq
is independent of $t$ and $\lambda$. This pairing is known as the {\em
  intersection pairing}.
Suppose now that $\gamma$ is a simple loop in $F^\circ$, i.e., a loop that starts at $\lambda^\circ$,
approaches one of the punctures $u_i^\circ$ along a path $\gamma'$
that ends at a point sufficiently close to $u_i^\circ$, goes around
$u_i^\circ$, and finally returns back to $\lambda^\circ$ along
$\gamma'$. By analyzing the second structure connection near
$\lambda=u_i$ it is easy to see that up to a sign there exists a unique $a\in H$
such that $(a|a)=2$ and the monodromy transformation of $a$ along
$\gamma$ is $-a$. The monodromy transformation representing $\gamma\in
\pi_1(F^\circ,\lambda^\circ)$ is the reflection defined by the
following formula:
\ben
w_a(x)=x-(a|x) a.
\een
Let us denote by $R$ the set of all $a\in H$ as above determined by
all possible choices of simple loops in $F^\circ$. We refer to the
elements of $R$ as reflection vectors.

\subsection{The ring of modular functions}

Our main interest is in the period map
\ben
Z:((M\times \CC)')^\sim \to H^*,
\quad 
(t,\lambda)\mapsto Z(t,\lambda) 
\een
where $((M\times \CC)')^\sim$ is the universal cover of $(M\times
\CC)'$ and $Z(t,\lambda)\in H^*$ is defined by 
\ben
\langle Z(t,\lambda),\alpha\rangle :=Z_\alpha(t,\lambda) = (I^{(-1)}_\alpha(t,\lambda),1).
\een
Recall that we require that the Frobenius manifold $M$ satisfies
condition (iv) from Section \ref{sec:2nd-str-conn}.  Under this
condition 
the flow of the unit vector field $\mathbf{1}$ defines a free action
of $\CC$ on $M=\CC\times B$ 
\ben
\CC\times M\to M,\quad  (x,t)\mapsto t+x\mathbf{1},
\een
where for $t=(t_1,{}'t)\in \CC\times B$ we define
$t+x\mathbf{1}:=(t_1+x,{}'t)$. 
The period map has the following translation symmetry
\beq\label{tr-sym}
Z(t,\lambda) = Z(t-\lambda\mathbf{1},0).
\eeq
Therefore, we will restrict our analysis to the case $t_1=0$, i.e., we
will assume that $t\in B$ and that the period map is defined on the
universal cover of 
\ben
X:=(B\times \CC)' = \{ (t,\lambda)\in B\times \CC \ |\
\operatorname{det}(\lambda-E\bullet)\neq 0\}.
\een

Let us denote by $\Omega\subset H^*$ the image of the period map
$Z$. This is a $W$-invariant subset which will be called the
{\em period domain}. In general very little is known about such
domains. For example it would be interesting to classify
semi-simple Frobenius manifolds such that the action of $W$ on $\Omega$ is properly
discontinuous and the quotient $[\Omega/W]$ is an orbifold whose
coarse moduli space is isomorphic to the Frobenius manifold
$M$. Furthermore, we would like to introduce the ring of
{\em modular functions}
\ben
\mathcal{M}(\Omega,W) := \{f\in \Gamma(\Omega,\O_{H^*})^W\ |\ f\circ
Z\in \O(B\times \CC)\},
\een
where $\Gamma(\Omega,\O_{H^*})^W$ is the ring of $W$-invariant
holomorphic functions in $\Omega$. Note that in general if $f\in
\Gamma(\Omega,\O_{H^*})^W$ is an arbitrary function, then the
composition $f\circ Z$ defines a holomorphic function on  $(B\times
\CC)'$. The condition in the above definition requires that $f\circ Z$
extends analytically across the discriminant.

\subsection{Example}
In general, one might try to investigate a more general period map
defined by $(I^{(n)}(t,\lambda),\phi_i)$ for any $n\in \CC$ and $1\leq
i\leq N$. Since $ (I^{(n)}(t,\lambda),\phi_i) =
-\partial_{t_i}(I^{n-1}(t,\lambda),\mathbf{1})$ the choice $i=1$,
 which yields $\phi_1=\mathbf{1}$, is quite natural. 
Let us discuss the possibility of choosing different values of $n\in
\CC$ in the case of $A_1$-singularity. Since the periods defined by
the second structure connection locally near a generic point on the
discriminat have the same leading order terms as the periods of
$A_1$-singularity, one can get a good intuition of what values
of $n$ could be interesting to investigate.

The period map takes the form
\ben
Z:(\CC\setminus{\{0\}})^\sim \rTo H^*\cong \CC,\quad
Z(\lambda)=\frac{\lambda^{-n-\tfrac{1}{2}}}{\Gamma(-n+\tfrac{1}{2})}.
\een
The monodromy group is a cyclic group $W=(w) $ and the action of $w$ on $\CC$ is multiplication by
$e^{-\pi\sqrt{-1} (2n+1)}$. If $n$ is a complex non-rational number
then the quotient $\CC/W$ might even fail
to be a Hausdorf space.   Let us assume that $n\in \mathbb{Q}$, so that the
quotient  is an orbifold. The quotient space $\CC/W$ has the structure
of a smooth complex manifold isomorphic to $\CC$. The isomorphism is
induced from
\ben
\pi:\CC/W\to \CC,\quad \lambda\mapsto \lambda^q,
\een
where we write $2n+1=\tfrac{2p}{q}$ with $q>0$ and $p$ relatively prime
integers. The period map induces a holomorphic map
\ben
\CC\setminus{\{0\}}\rTo \CC/W=\CC,\quad
\lambda\mapsto = 
\frac{\lambda^{-p} }{\Gamma(1-\tfrac{p}{q} )}.
\een
The above map extends holomorphically across the discriminant
$\lambda=0$ if and only if $p\leq 0$. Moreover it is an isomorphism if
and only if $p=-1$, i.e., $n=-\tfrac{1}{2}-\tfrac{1}{q}$.  

The conclusion is that in general, there might be other period maps
which would allow us to identify the Frobenius manifold with the orbit
space of the corresponding monodromy group. Our choice $n=-1$ is
motivated by the applications of semi-simple Frobenius manifold to
integrable hierarchies and representations of lattice vertex algebras
(see the Appendix for more details).

\subsection{Riemann--Hilbert problem for Gromov--Witten invariants}
The notion of a Frobenius manifold was invented by Dubrovin \cite{Du}
in order to give a geometric interpretation of the properties of
quantum cohomology of a smooth projective variety $X$. It was
conjectured by Givental in \cite{G1} and proved by Teleman in
\cite{Te} that if the quantum cohomology is
semi-simple as a Frobenius manifold then the higher genus Gromov--Witten invariants are
uniquely determined by genus-0, i.e., by the underlying semi-simple Frobenius
structure. On the other hand, the entire semi-simple Frobenius
structure is uniquely determined by the second structure
connection. The latter is a Fuchsian connection and hence can be
recovered uniquely as a solution to a classical Riemann--Hilbert
problem (see \cite{Du2,Mi}). The problem that we are interested in is how to express the
higher genus Gromov--Witten invariants of a manifold with semi-simple
quantum cohomology in terms of the monodromy data of the second
structure connection. We refer to such a problem as the
{\em Riemann--Hilbert problem for Gromov--Witten invariants}. By
definition such a problem has a solution. Namely, by solving a
classical Riemann Hilbert problem we can recover the second structure
connection from its monodromy data, then we can recover the
semi-simple Frobenius structure, and finally it remains to recall
Givental's higher genus reconstruction. The solutions to the
Riemann--Hilbert problems are usually highly
transcendental. Therefore, it seems that the dependence of the
Gromov--Witten invariants on the monodromy data should also be quite
complicated. However, in a series of examples (see \cite{BM,GM,MT,MST}) the monodromy data
leads to a highest weight representation of a vertex algebra or to
the Hirota bilinear equations of an integrable hierarchy which allow
us to uniquely determine the invariants. In other words, we are
looking for a representation of a Lie algebra or more generally a
vertex algebra which will allow us to express all invariants via the
monodromy data in a simple combinatorial way. The main point is not
that we will find a way to compute Gromov--Witten invariants, but
rather that we can understand a conceptual question. Namely is there a
strong relation between semi-simple Frobenius manifolds and Lie
algebras. Note that the set of reflection vectors can be used to
generalize the notion of root systems and to propose various
constructions of Lie algebras. This was actually done in the settings
of singularity theory by several
authors (e.g. see \cite{SY,Sl}). The problem is whether such Lie algebras
have interesting applications. The Riemann--Hilbert problem for
Gromov--Witten invariants can be viewed as a motivation to develop Lie
algebra theory for semi-simple Frobenius manifolds. 

All examples in which some interesting relation to Lie theory was
established have conformal dimension
$D\leq 1$ ($D$ is the complex dimension of the manifold $X$ in the case of quantum
cohomology). The case of quantum cohomology of $\PP^2$ is a very good
candidate to make progress in conformal dimension $>1$ because the geometry
of $\PP^2$ relevant for studying quantum cohomology and mirror
symmetry is very well understood.  The problem in the current paper
comes from our attempt to generalize the work in \cite{BM}. Namely, we
would like to find differential operator constraints for the total
descendant potential of $\PP^2$. This is still a very difficult
problem. We will argue in the appendix that the genus-0 reduction of a differential operator
constraints yields a Hamilton--Jacobi equation given by a holomorphic
function in $\mathcal{M}(\Omega,W)$. This motivates to some extent our
interest in the ring of modular functions.

\begin{definition}\label{def:inverse}
The period map $Z$ is said to be {\em invertible} if there exists a set of
modular functions $f_i\in \mathcal{M}(\Omega,W)$ $(1\leq i\leq N)$
such that the set of holomorphic functions $f_i\circ Z$ $(1\leq i\leq
N)$ is a coordinate system on $B\times \CC$. A set of such modular
functions $\{f_i\}_{i=1}^N$ is called the {\em inverse of the period
  map}. 
\end{definition}
There are two reasons why we are interested in finding the inverse of
the period map. The first one is related to the discussion above. We
expect that if the period map is invertible then the corresponding modular
functions $f_i$ will give a complete set of recursion relations, which
would allow us to determine the genus-0 total descendant potential in
terms of the monodromy data of the Frobenius manifold via an explicit
recursion (see Appendix for more details).

The second reason is related to the problem of uniformizing a
semi-simple Frobenius manifold.
We expect that  semi-simple Frobenius manifolds relevant in the study
of mirror symmetry are quotients of a simply connected domain by a
discrete group. At this point we can only speculate, but we believe
that the problem of uniformizing the Frobenius manifold corresponding
to the quantum cohomology of some smooth projective variety $X$ is
related to the problem of constructing the manifold of stability
conditions of the bounded derived category $D^b(\operatorname{Coh}\, X)$.

{\bf Acknowledgements.}
I would like to thank M. Kapranov for pointing out to me that if the
solutions to some differential equations satisfy a quadratic relation,
then this is an indication that the differential equation itself  might
be a symmetric square. I would like to thank the anonymous referee for
the useful comments and remarks that helped me to improve the
exposition. This work is partially supported by JSPS
Grant-In-Aid (Kiban C) 17K05193 
and by the World Premier International Research Center Initiative (WPI
Initiative),  MEXT, Japan.

\section{Quantum cohomology of $\mathbb{P}^2$}

From now on we will work only in the settings of quantum cohomology of
$\PP^2$. The goal of this section is to introduce the necessary
notation and to state our results.

\subsection{Frobenius manifold structure}

Let $H=H^*(\mathbb{P}^2;\CC)$ and $t_i$ $(1\leq i\leq 3)$ be the
linear coordinates on $H$ corresponding to the basis
\ben
\phi_i:=p^{i-1},\quad i=1,2,3,
\een
where $p=c_1(\O(1))$ is the hyperplane class. Quantum cohomology
defines a Frobenius manifold structure on the space
\ben
M=\{(t_1,Q,t_3)\in \CC\times \CC^* \times \CC\ |\ |t_3Q^{1/3}|<\epsilon\},
\een
where $\epsilon$ is a sufficiently small positive  real number and the
coordinate $Q:=e^{t_2}$ is identified with the Novikov
variable. By definition, the linear coordinates $t_i$ are flat, the
Frobenius pairing is given by the Poincare pairing
\ben
(\partial_i,\partial_j) = \delta_{i+j,4},\quad 1\leq i,j\leq 3,
\een
where $\partial_i=\partial/\partial t_i$, while the multiplication is
given by the quantum cup product. The latter is defined by 
\ben
(\partial_i\bullet \partial_j,\partial_k) := \frac{\partial^3
  F(t)}{\partial t_i\partial t_j\partial t_k},
\een
where $F(t)$ is the genus-0 potential 
\ben
F(t)=\sum_{l,d=0}^\infty
\frac{1}{l!}\langle t,\dots,t\rangle_{0,l,d} ,
\een
where $t:=t_1+t_2p+t_3 p^2$. Following Kontsevich and Ruan--Tian we can
derive an explicit recursive formula for $F$ as follows. Using the
string equation, the divisor equation, and the dimension formula of the virtual fundamental
cycle we get that $F$ has the form
\beq\label{potential}
F(t) = \frac{1}{2} (t_1^2 t_3 +t_1t_2^2) +\sum_{d=1}^\infty \frac{N_d}{(3d-1)!}\, Q^d t_3^{3d-1},
\eeq
where the coefficient $N_d$ can be interpreted as the number of rational curves in
$\PP^2$ of degree $d$ passing through $3d-1$ points in general
position. The system of WDVV equations contains a  single non-trivial
equation
\ben
F_{333} = F_{233}^2-F_{222}F_{233},
\een 
where the index $i$, $i=2,3$, denotes partial derivative with respect
to $t_i$. Comparing the coefficients in front of $Q^d$ yields 
\ben
N_d= \sum_{m=1}^{d-1} 
\Big( 
{3d-4\choose 3m-2} \,m^2(d-m)^2 - 
{3d-4\choose 3m-3} \,m (d-m)^3
\Big) N_{m}N_{d-m},  
\een
which together with $N_1=1$ determines $N_d$ for all $d>1$. The first
few values are
\ben
N_1=N_2=1,\ N_3=12,\ N_4= 620,\ N_5=87304,\ N_6=26312976,\ \dots.
\een
Let us point out that the number $\epsilon$ in the definition of the
domain $M$ is chosen in such a way that the radius of convergence of
the series \eqref{potential} is $\epsilon |Q|^{-1/3}$.

Furthermore, the Euler vector field has the
form
\ben
E=t_1\partial_1 -t_3 \partial_3 + 3 \partial_2
\een
and the Hodge grading operator is
\ben
\theta: H\to H,\quad \theta  = \diag( 1,0,-1). 
\een

\subsection{The $\Gamma$-integral structure of Iritani}

Let us recall the notation of Section \ref{sec:pv}. Following
Givental (see \cite{G2}), we equip the quantum cohomology with
calibration 
\ben
S(t,z) = 1+S_1(t)z^{-1}+S_2(t)z^{-2}+\cdots,
\quad
S_k(t)\in \operatorname{End}(H)
\een
defined by 
\ben
(S(t,z)\phi_i,\phi_j) = (\phi_i,\phi_j)+\sum_{k=0}^\infty \langle
\phi_i \psi^k,\phi_j\rangle_{0,2}(t)\, z^{-k-1}.
\een
The fundamental solution corresponding to
such calibration is $S(t,z)z^\theta z^{-\rho}$, where the nilpotent
operator $\rho$ is given by classical cup product
multiplication by $c_1(T\PP^2)=3p$. 

There is a very elegant way to describe the reflection lattice, i.e.,
the $\ZZ$-submodule of $H$ spanned by all reflection vectors $a\in R.$
Namely, using Iritani's $\Gamma$-class modification of the Chern character
map, we will obtain an explicit description of all reflection
vectors in terms of $K^0(\PP^2)$ -- the $K$-ring of topological vector
bundles on $\PP^2$. Following Iritani's construction in \cite{Ir}, let
us introduce the map  
\ben
\Psi: K^0(\PP^2)\to H,
\een
defined by 
\ben
\Psi(E) = \frac{1}{\sqrt{2\pi}}\, \Gamma(1+p)^3 \, (2\pi\sqrt{-1})^{\rm
  deg} ch(E),
\een
where ${\rm deg}:H\to H$ is the degree operator ${\rm
  deg}(p^i)=ip^i$ and the $\Gamma$ function should be expanded as a
Taylor series at $p=0$, i.e., 
\ben
\Gamma(1+p) = 1+\Gamma'(1) p +\frac{1}{2} \Gamma''(1) p^2.
\een
Recall that $K^0(\PP^2)=\ZZ[L]/(L-1)^3$, where $L=\O(1)$. The above
formula gives 
\ben
\Psi(L^m) = 
\frac{1}{\sqrt{2\pi}}\, \Gamma(1+p)^3 \, e^{2\pi\sqrt{-1} m p}.
\een 
Slightly abusing the notation we identify $L^m$ with its image
$\Psi(L^m)$.

Let
us choose a reference point $(t^\circ,\lambda^\circ)\in (M\times
\CC)'$ such that $t_1^\circ=t_2^\circ=t_3^\circ=0$ and $\lambda^\circ$
is a sufficiently large real number. Recall that for $t_3=0$ the
quantum cup product 
$\bullet_t$ turns $H$ into the following algebra 
\ben
(H,\bullet_t) = \CC[p]/(p^3-Q).
\een
We get that the eigenvalues of $E\bullet_{t^\circ}$ are
$u_i^\circ=3\zeta^{-i+1}$ $(1\leq i\leq 3)$, where
$\zeta=e^{2\pi\sqrt{-1}/3}.$  Let us denote by $[a,b]$
$(a,b\in \CC)$ the straight segment in $\CC$ from $a$ to
$b$. 
\begin{figure}[htbp]
\begin{center}
\scalebox{0.55}{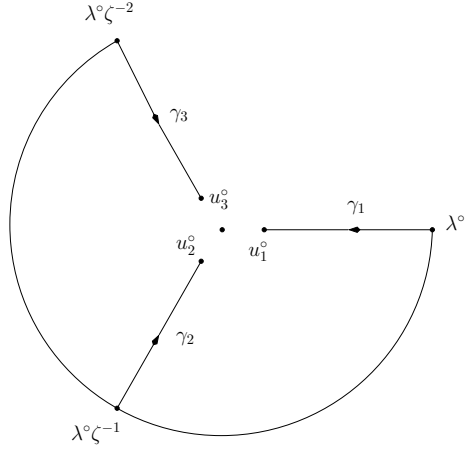}
\caption{Reflection vectors}
\label{fig:rv}
\end{center}
\end{figure}
Let $\gamma_i$ $(1\leq i\leq 3)$ be the composition of the arc 
\ben
\lambda(s)=\lambda^\circ
e^{-2\pi \sqrt{-1} s/3},\quad s\in [0, i-1]
\een 
and the straight segment
$[\zeta^{-i+1} \lambda^\circ,u_i^\circ]$ (see Figure \ref{fig:rv}).

\begin{proposition}\label{prop:refl-vec}
  Let $\gamma_i$ $(1\leq i\leq 3)$ be the paths constructed
  above. Then $L^{i-1}$ is the reflection vector
corresponding to the path $\gamma_i$.
\end{proposition}
\begin{corollary}\label{cor:refl}
The set of all reflection vectors is given by 
\ben
R=W\cdot 1\cup W\cdot L\cup W\cdot L^2.
\een
\end{corollary}
The proof of Proposition \ref{prop:refl-vec} and Corollary
\ref{cor:refl} will be given in Section \ref{sec:refl}.

Since the intersection pairing \eqref{inters-pairing} is independent
of $t$, by setting $t_1=t_3=0$ and passing to the limit
$\operatorname{Re}(t_2)\to -\infty$ (i.e. $Q\to 0$) we get the
following formula for the intersection pairing
\ben
(L^m|L^n) =
(\widetilde{I}^{(0)}_{L^m}(\lambda),(\lambda-\rho)
\widetilde{I}^{(0)}_{L^n}(\lambda))=
2+(m-n)^2. 
\een
Recalling Corollary \ref{cor:refl} we get that the lattice 
\ben
\operatorname{Im}(\Psi):= \ZZ+ \ZZ\, L +\ZZ\, L^2\subset H
\een
coincides with the reflection lattice. 
Furthermore, recalling Proposition \ref{prop:refl-vec} we get  Dubrovin's
result (see \cite{Du2}) that the monodromy group $W\cong 
\operatorname{PSL}_2(\ZZ)\times \{\pm 1\}$. The construction of this  
isomorphism amounts to choosing an appropriate $\mathbb{Q}$-basis
$(E_1,E_2,E_3)$ of the reflection lattice. In our notation this basis is given by 
\ben
E_1=1+2L-L^2,\quad
E_2=-3+4L-L^2,\quad
E_3=1-2L+L^2.
\een

\subsection{The period map for quantum cohomology of $\PP^2$}
Now we are in position to state our results. We identify $H^*\cong
\CC^3$ via the linear functions on $H^*$ corresponding to $E_i$,
$i=1,2,3$. The period map takes the form 
\ben
Z(t,\lambda)=(Z_1(t,\lambda),Z_2(t,\lambda),Z_3(t,\lambda)),
\een
where $Z_i(t,\lambda):=Z_{E_i}(t,\lambda)$. As we have explained in
the introduction, we may restrict our analysis to parameters 
$
t\in B:=\{t\in M\ |\ t_1=0\}.
$
Put $B_{\rm small}=\{t\in B\ |\ t_3=0\}=\CC^*$ and 
\ben
X_{\rm small}:=\left. (M\times \CC)'\right|_{t_1=t_3=0} =\{ (Q,\lambda)\in
\CC^*\times \CC\  |\ \lambda^3-27Q \neq 0\}.
\een
Let us introduce the domain 
\ben
\Omega_{\rm small} := \{ z\in (\CC^*)^3\ |\ z_2^2=4z_1z_3,\ \operatorname{Im}(-z_2/(2z_3))>0\}.
\een
Let us point out that there is a natural isomorphism 
\ben
\Phi_{\rm small}: \HH\times \CC^*\to \Omega_{\rm small},\quad
(\tau,y)\mapsto (\tau^2 y, -2\tau y, y)
\een
under which the action of the monodromy group takes a very simple form
(see Lemma \ref{le:mod-action}). 
We will prove later on (see Lemma \ref{X-retract-Y}) that $X_{\rm small}$ is a
deformation retract of $X$. Therefore, the universal cover
$\widetilde{X}_{\rm small}$ is an analytic submanifold of
$\widetilde{X}$ and we can introduce the restriction of  the period map $Z_{\rm
  small}:=Z|_{\widetilde{X}_{\rm small}}$. 
Recall the Eisenstein series 
\ben
E_2(\tau) & = & 1-24 \sum_{m=1}^\infty \frac{mq^m}{1-q^m},\\
E_4(\tau) & = & 1+240\,\sum_{m=1}^\infty \frac{m^3 q^m}{1-q^m},\\
E_6(\tau) & = & 1-504\,\sum_{m=1}^\infty \frac{m^5 q^m}{1-q^m},
\een
where $q=e^{2\pi\sqrt{-1} \tau}$. Our first result can be stated as follows.
\begin{theorem}\label{t1}
a) The image of $Z_{\rm small}$ is $\Omega_{\rm small}\setminus{\{E_6(-z_2/(2z_3))=0\}}$.

b) Let $\widetilde{\pi}_{\rm small}:\widetilde{X}_{\rm small}\to X_{\rm small}$ be the universal cover and
\ben
\pi_{\rm small} :\Omega_{\rm small} \to \CC^*\times \CC,\quad (\tau,y)\mapsto(Q(\tau,y),\lambda(\tau,y)),
\een 
be the map defined by 
\ben
Q(\tau,y)&:= & \frac{8}{27}\, (2\pi/y)^{6} (E_4^3(\tau)-E_6^2(\tau)) \\
\lambda(\tau,y) & := &  2(2\pi/y)^2  E_4(\tau),
\een
where $(\tau,y)$  is the coordinate system on $\Omega_{\rm small}$
introduced above. Then $\pi_{\rm small}\circ Z_{\rm small}=\widetilde{\pi}_{\rm small}$. 

c)
The fibers of the map $\pi_{\rm small}$ are the $W$-orbits in
$\Omega_{\rm small}$, i.e., 
the small quantum cohomology $B_{\rm small}\times \CC$ is the coarse moduli
space for the orbifold $[\Omega_{\rm small}/W]$.
\end{theorem}

Generalizing the results of Theorem \ref{t1} to big quantum
cohomology is a very challenging problem. We expect that $M_{\rm
  small}:=\CC\times B_{\rm  small}$ is an analytic subvariety in a
larger Frobenius manifold $N$ and that $M$ is just a tubular
neighborhood of $M_{\rm small}$ in $N$. We were able to prove an interesting
result about the holomorphic thickening of $Z_{\rm small}$, which
might be viewed as the first step towards constructing a global
Frobenius manifold.

The period map $Z$ maps a small open neighborhood of $\widetilde{X}_{\rm
  small}$ in $\widetilde{X}$ into a small open neighborhood of
$\Omega_{\rm small}$ in $\CC^3$. Therefore we have an induced map of
ringed spaces 
\ben
(\widetilde{X}_{\rm small},\widetilde{\iota}^{\,-1} \O_{\widetilde{X}}
)\to \Omega:=(\Omega_{\rm small},\iota^{-1} \O_{\CC^3}),
\een
where $\widetilde{\iota}:\widetilde{X}_{\rm small}\to \widetilde{X}$
and $\iota:\Omega_{\rm small}\to \CC^3$ are the natural inclusion
maps. The ring of regular functions on $\Omega$ is by definition
$\Gamma(\Omega_{\rm small},\O_{\CC^3})$, i.e., functions defined and
holomorphic in an open neighborhood of $\Omega_{\rm small}$ in
$\CC^3$. This ring is equipped with the action of the monodromy group
$W$. The pullback via the period map $Z$ defines a ring homomorphism
\ben
\Gamma(\Omega_{\rm small},\O_{\CC^3})^W\to \Gamma(X_{\rm
  small},\O_{X}),\quad f\mapsto f\circ Z,
\een
where $\Gamma(\Omega_{\rm small},\O_{\CC^3})^W $ is the subring of
$W$-invariant functions. The coordinate functions $Q=e^{t_2}$,
$t=t_3$, and $\lambda$ of $B\times \CC$ are elements of $\Gamma(X_{\rm
  small},\O_{X})$.  We will prove that $Q,t$, and
$\lambda$ are pullbacks via the period map of $W$-invariant
functions. Moreover, the latter have some interesting property, which can be stated as
follows. Let us construct an open neighborhood of $\Omega_{\rm small}$
as the image of the map 
\ben
\Phi: \HH^2\times \CC^* \to \CC^3, \quad (\tau_1,\tau_2,y)\mapsto
(\tau_1\tau_2 y, -(\tau_1+\tau_2)y, y),
\een
where $\HH=\{\tau\in \CC\ |\ \operatorname{Im}(\tau)>0\}$ is the upper half-plane.
The image of $\Phi$ is the coarse moduli space for the orbifold quotient
$[\HH^2\times \CC^*/\mu_2]$, where $\mu_2$ is the cyclic group of
order 2 whose generator acts on $\HH^2$ by permutation
$(\tau_1,\tau_2)\mapsto (\tau_2,\tau_1)$. 
\begin{theorem}\label{t2}
a) There are $W$-invariant functions in $\Gamma(\Omega_{\rm
  small},\O_{\CC^3})^W $ of the form
\ben
Q(\tau_1,\tau_2,y) & = & \frac{8}{27}(2\pi/y)^6
\sum_{n=0}^\infty Q_n(\tau_{12}) (\tau_1-\tau_2)^{2n} ,\\
\lambda(\tau_1,\tau_2,y) & = & 2(2\pi/y)^2
\sum_{n=0}^\infty \lambda_n(\tau_{12}) (\tau_1-\tau_2)^{2n} ,\\
t(\tau_1,\tau_2,y) & =& -\frac{1}{32} (\tau_1-\tau_2)^2 y^2,
\een
where $\tau_{12}:=(\tau_1+\tau_2)/2$, such that their pullbacks via
the period map coincide with the coordinate functions $Q,\lambda$, and $t$. 

b) The coefficients $Q_n(\tau), \lambda_n(\tau)$ $(n\geq 0)$ are
quasi-modular forms, i.e., they are polynomials in the Eisenstein
series $E_i(\tau)$, $i=2,4,6$. 
\end{theorem}

Note that if $\tau_1=\tau_2$ then we recover the formulas from Theorem
\ref{t1}. In particular
\ben
Q_0(\tau) = E_4(\tau)^3-E_6(\tau)^2,\quad \lambda_0(\tau) =
E_4(\tau). 
\een
We have computed the quasi-modular forms  $\lambda_n$ and $Q_n$ for
$n=1,2,3$. The answer is the following
\ben
\lambda_1(\tau) & = &  \frac{1}{40}\partial_\tau^2 E_4,\\
\lambda_2(\tau) & = & \frac{1}{4480}\partial_\tau^4 E_4
-\frac{\pi^4}{2016}\Delta,\\
\lambda_3(\tau) & = &
\frac{1}{967680}\partial_\tau^6E_4-\frac{\pi^4}{209664} \partial_\tau^2\Delta-\frac{\pi^6}{101088}E_4\Delta,
\een
and 
\ben
Q_1(\tau) & = &  \frac{1}{104} \partial_\tau^2 \Delta +
\frac{\pi^2}{26} E_4\Delta,\\
Q_2(\tau) & = & \frac{1}{24960}\partial_\tau^4\Delta+
\frac{\pi^2}{2704}E_4 \partial_\tau^2 \Delta+
\frac{\pi^2}{1040} \Delta \partial_\tau^2 E_4+
\frac{17\pi^4}{20280}E_4^2 \Delta,\\
Q_3(\tau) &= & \frac{1}{10183680} \partial_\tau^6 \Delta +\frac{ 
  1611 \pi^2}{37856000}(E_4\partial_\tau^2E_4)^2 +\frac{ 3 \pi^2}{1514240}
  E_4\partial_\tau^4\Delta \\
&&+
\Big( \frac{3 \pi^2}{116480}\Delta - \frac{537
  \pi^2}{26499200}E_4^3\Big) \partial_\tau^4E_4 +
\frac{239 \pi^4}{2839200} E_4\Delta\partial_\tau^2E_4 \\
&&
-\frac {319 \pi^6}{
  26732160}\Delta^2 +\frac{3977 \pi^6}{202718880}
  E_4^3\Delta,
\een
where $\Delta:=E_4^3-E_6^2$.

\section{Reflection vectors in quantum cohomology of $\PP^2$}\label{sec:refl}
The main goal in this section is to prove Proposition
\ref{prop:refl-vec}.

\subsection{Mirror symmetry for the calibration}
In this section we recall an identity expressing the operator series
$S(t,z)$ in terms of an oscillatory integral. Let us denote by
$S(Q,z)$ the restriction of $S(t,z)$ to $t_1=t_3=0$.  Recall that
a Givental's mirror model of $\PP^2$ is given by the family of
functions
\ben
f(x,Q):=x_1+x_2+\frac{Q}{x_1x_2},\quad x=(x_1,x_2)\in (\CC^*)^2
\een
depending on the parameter $Q\in \CC^*$ and the holomorphic form on
$(\CC^*)^2$ defined by $\omega=\tfrac{dx_1}{x_1}\wedge
\tfrac{dx_2}{x_2}$.

Let us fix $Q=e^{2\pi\sqrt{-1} s}q$, where $q>0$ and $s$ are real
numbers and define the semi-infinite cycle $\alpha_Q:=e^{2\pi\sqrt{-1}
  s/3}(\RR_{>0})^2$, i.e., $\alpha_Q$ consists of all $x\in (\CC^*)^2$
such that $x_i = e^{2\pi\sqrt{-1} s/3}|x_i|$ ($i=1,2$). The following
result is due to Iritani (see \cite{Ir}, Theorem 4.14).
\begin{lemma}\label{le:osc-S}
  Suppose that $z=e^{2\pi\sqrt{-1} s/3} w$ where $w<0$ is a real
  number. Then the following formula holds
  \beq\label{central-charge}
  \int_{\alpha_Q} e^{f(x,Q)/z} \omega = \sqrt{2\pi} (-w) (S(q,w)
  (-w)^\theta (-w)^\rho \Psi(\O),\mathbf{1}) 
  \eeq
\end{lemma}
\proof
Note that the LHS is independent of $s$, because if we substitute
$x_i=e^{2\pi\sqrt{-1} s/3}y_i$, then the form $\omega$ is invariant,
the cycle $\alpha_Q$ is transformed into $\alpha_q$, and the function
\ben
f(x,Q) / z = \Big( y_1 + y_2 + \frac{q}{y_1y_2}\Big)/w.
\een
Therefore, we may assume that $s=0$ and that $Q=q>0$ and $z=w<0$ are real
numbers.

Let us simplify the RHS of \eqref{central-charge}. Recalling the definition of
$S$ we get
\ben
S(q,w)^T \mathbf{1} = 
\mathbf{1} +\sum_{k=0}^\infty
\sum_{i=1}^3\langle \phi_i\psi^k,\mathbf{1}\rangle_{0,2}(t) \phi^i
w^{-k-1}.
\een
Using the divisor and the string equations the above formula can
be transformed into
\ben
S(q,w)^T \mathbf{1}=
\Big(1+\sum_{k=0}^\infty\sum_{d=1}^\infty q^d\, 
\langle \phi_i\psi^k\rangle_{0,1,d} \, \phi^i
w^{-k-2}\Big)
e^{t/w}, 
\een
where recall that $t=t_2 p$ and $q=e^{t_2}$. 
The above formula is by definition Givental's J-function of
$\PP^2$. The J-functions of all projective spaces and certain classes
of complete intersections are computed explicitly in \cite{G4},
Theorem 9.1. We get
\ben
S(q,w)^T \mathbf{1}  = \Big(1+\sum_{d=1}^\infty
\frac{q^d}{(p+w)^3\cdots (p+d w)^3}\Big) e^{t/w}.
\een
Using the commutation relation $(-w)^{-\theta} p = (-w) p
(-w)^{-\theta}$ and $\theta(\mathbf{1})=\mathbf{1}$ we get 
\ben
(-w)(-w)^\rho (-w)^{-\theta} S(q,w)^T \mathbf{1} =
\Big(1+\sum_{d=1}^\infty
\frac{e^{(t_2-3\log(-w))d}}{(p-1)^3\cdots (p-d )^3}\Big) e^{-(t_2-3\log(-w))p} .
\een
Put $\tau:=t_2-3\log (-w)$. Then the RHS of \eqref{central-charge}
takes the form
\ben
\int_{[\PP^2]} e^{-\tau p} \Gamma(1+p)^3
\Big(1+\sum_{d=1}^\infty
\frac{e^{\tau d} }{(p-1)^3\cdots (p-d )^3}\Big) ,
\een
where $[\PP^2]$ is the fundamental class of $\PP^2$ and we used the
following formulas
$\Psi(\O)=\tfrac{1}{\sqrt{2\pi}}\Gamma(1+p)^3$,
$\theta^T=-\theta$, 
$\rho^T=\rho$.
The above integral can be written as a residue.  Namely
if $f(p)\in\CC[p]$ is a polynomial, then we have 
\ben
\int_{[\PP^2]} f(p) =
\operatorname{Res}_{p=0} f(p) \tfrac{dp}{p^3} .
\een
Using this fact and
the identity $\Gamma(1+p)=p(p-1)\cdots (p-d) \Gamma(p-d)$ we get that
the RHS of formula \eqref{central-charge} has the form
\ben
\sum_{d=0}^\infty \operatorname{Res}_{p=-d}
\Big(e^{-\tau p} \, \Gamma(p)^3\Big) dp .
\een
Let us transform the LHS of \eqref{central-charge}. 
Note that the oscillatory integral in \eqref{central-charge} can be
written as
\ben
I(\tau) = \int_{\RR^n} \exp\Big( -e^{y_1} -e^{y_2} 
-e^{\tau-y_1-y_2}\Big) dy_1dy_2,
\een
where $\tau=t_2-3\log (-w)$ and we used the substitution 
$x_i=-w e^{y_i}$ $(1\leq i\leq 2)$. If $\epsilon>0$ is a fixed real
number then $I(\epsilon,\tau):=e^{\epsilon \tau} I(\tau)\in L^1(\RR)$. Recalling the Fourier
inversion formula we get 
\ben
I(\epsilon,\tau) =
\frac{1}{2\pi} \lim_{a\to +\infty} \int_{-a}^a \int_{\RR}
e^{\sqrt{-1}(\tau-s)\xi} I(\epsilon,s)\, ds\,  d\xi.
\een 
Let us substitute in the above formula $I(\epsilon,\tau)=e^{\epsilon \tau}
I(\tau)$, $I(\epsilon,s)=e^{\epsilon s}
I(s)$, and make the substitution $p=\epsilon - \sqrt{-1}\xi$. We get 
\ben
I(\tau) = \frac{1}{2\pi\sqrt{-1} } \int_{\epsilon + \sqrt{-1}\RR} \int_\RR
  e^{-p(\tau-s)} I(s) \, ds\, dp,
\een
where the orientation of the contour $\epsilon + \sqrt{-1}\RR$ is from
$\epsilon - \sqrt{-1}\infty$ to $\epsilon + \sqrt{-1}\infty$.
The integral with respect to $s$ can be computed explicitly as follows:
\ben
\int_\RR e^{ps} I(s) d s = \prod_{i=1}^{3} \int_\RR e^{-e^{y_i} + p
  y_i} d y_i = \Gamma(p)^3,
\een
where the first equality is justified as follows.  The function 
\ben
\exp\Big( p s  -e^{y_1} - e^{y_2}-
e^{s-y_1-y_2}\Big) 
\een
with respect to $(y_1,y_2,s)\in \RR^3$ is of class $L^1$ for all $p$ satisfying
$\operatorname{Re}(p)>0$. Therefore we can use the Fubini's theorem
to transform the iterated integral $\int_\RR ds \int_{\RR^2} dy_1dy_2$ as a
multiple integral. It remains only to change the integration variables
$(y_1,y_2,s)$ into $(y_1,y_2,y_3)$ via the substitution
$s=y_1+ y_2+ y_3$. 
 
The oscillatory integral takes the form 
\ben
I(t) = \frac{1}{2\pi\sqrt{-1} } \int_{\epsilon + \sqrt{-1}\RR}  
e^{-p t} \Gamma(p)^{3} dp .
\een
The integrand has poles at $p=-d$ for $d=0,1,2,\dots$. Using the
Cauchy residue theorem and some standard estimates based on the
Stirling formula for the $\Gamma$-function we get 
\ben
I(t) = \sum_{d=0}^\infty \operatorname{Res}_{p=-d} \Big( e^{-p t}
\Gamma(p)^{3} \Big) dp.
\een
This completes the proof of \eqref{central-charge}.
\qed
\begin{remark}
The proof of Lemma \ref{le:osc-S} is inspired by the work of K. Hori
and M. Romo (see \cite{HR}).  The main idea is that the Laplace
transform with respect to the Novikov's variables of the oscillatory
integral is the integrand of a Mellin--Barnes integral. Therefore,
using the inverse Laplace transform we can identify the oscillatory
integral with a Mellin--Barnes integral. In fact, this method seems to be quite general and
although there are technical difficulties it will be interesting to
investigate more complicated targets.
\end{remark}

\subsection{Twisted thimble integrals}
Let us continue to work with $Q=e^{2\pi\sqrt{-1} s}q$ and
$z=e^{2\pi\sqrt{-1} s/3}w$ such that $q>0$, $w<0$, and $s$ are real
numbers. The key object that would allow us to prove Proposition
\ref{prop:refl-vec} is the following integral
\ben
\mathcal{I}^{(-m)}(Q,\lambda)=
\int_{\alpha_{Q,\lambda}}
\frac{(\lambda-f(x,Q))^{m-\tfrac{1}{2}} }{\Gamma(m+\tfrac{1}{2}) }\, \omega,
\een
where $m\geq 1$ is an integer, the parameter
$\lambda=e^{2\pi\sqrt{-1} s/3}r$ with $r\geq 3 q^{1/3}$ a real number,
and $\alpha_{Q,\lambda}\subset (\CC^*)^2$ is the 
submanifold with boundary defined by 
\ben
\{x\in \alpha_Q\ |\ |f(x,Q)|\leq |\lambda|\} .
\een
The integrand is a multivalued analytic function. We fix an analytic
branch as follows. Note that
\ben
\lambda-f(x,Q) = e^{2\pi\sqrt{-1} s/3} \Big(
r-|x_1|-|x_2| - \frac{q}{|x_1 x_2|}\Big)
\een
and that the expression in the brackets is non-negative due to the
mean arithmetic and mean geometric inequality and our assumption that
$r\geq 3 q^{1/3}$. We define
\ben
(\lambda-f(x,Q) )^{m-\tfrac{1}{2} }= e^{\pi\sqrt{-1} s(2m-1)/3} \Big(
r-|x_1|-|x_2| - \frac{q}{|x_1 x_2|}\Big)^{m-\tfrac{1}{2} }.
\een
Let us point out that $\alpha_{Q,\lambda}$ is a {\em Lefschetz
  thimble}, i.e., it is swept out by vanishing cycles in the following
sense. The function $f(x,Q)$ induces a  map
\ben
\alpha_{Q,\lambda} \to [ u_Q ,\lambda],\quad
x\mapsto f(x,Q),
\een
where $u_Q = 3 Q^{1/3} = 3e^{2\pi\sqrt{-1} s/3}q^{1/3}$ is the
critical value of $f(x,Q)$ corresponding to the critical point
$\xi_Q=(Q^{1/3},Q^{1/3})$. The fiber
over $\mu\in (u_Q,\lambda]$ is diffeomorphic to a circle, while the
fiber over $u_Q$ is the critical point $\xi_Q$.  
\begin{lemma}\label{le:thimble-osc}
Suppose that $Q=e^{2\pi\sqrt{-1} s/3}q$ and $z=e^{2\pi\sqrt{-1}
  s/3}w$, where $q>0$ and $w<0$ are real numbers. Then
\beq\label{thimble-osc}
\int_{u_Q}^\infty e^{\lambda/z} \mathcal{I}^{(-m)}(Q,\lambda)d\lambda
=
(-z)^{m+\tfrac{1}{2}} \, \int_{\alpha_Q} e^{f(x,Q)/z} \omega,
\eeq
where the integral on the LHS is along the ray
$\lambda=e^{2\pi\sqrt{-1} s/3}r$, $r\in [3q^{1/3},+\infty)$. 
\end{lemma}
\proof
Using Fubini's theorem we have
\beq\label{I-fubini}
\mathcal{I}^{(-m)}(Q,\lambda) =
\int_{u_Q}^\lambda
\frac{(\lambda-\mu)^{m-\tfrac{1}{2}} }{ \Gamma(m+\tfrac{1}{2})}
\, \int_{\alpha_{Q,\mu}} \frac{\omega}{df}\, d\mu.
\eeq
Therefore, the LHS of \eqref{thimble-osc} can be written as
\ben
\int_{u_Q}^\infty 
\int_{u_Q}^\lambda
\Big(
e^{\lambda/z}\,
\frac{(\lambda-\mu)^{m-\tfrac{1}{2}} }{ \Gamma(m+\tfrac{1}{2})} \,
\int_{\alpha_{Q,\mu}} \frac{\omega}{df}\Big)\, d\mu\,d\lambda.
\een
Using $e^{\lambda/z}=e^{(\lambda-\mu)/z} e^{\mu/z}$ and changing the
order of integration we get
\ben
\int_{u_Q}^\infty \, e^{\mu/z} \Big(
\int_{\mu}^\infty
e^{(\lambda-\mu)/z}\,
\frac{(\lambda-\mu)^{m-\tfrac{1}{2}} }{ \Gamma(m+\tfrac{1}{2})}\, d\lambda
\Big)\,
\int_{\alpha_{Q,\mu}} \frac{\omega}{df}\Big)\, d\mu.
\een
Let us change the integration variable $\lambda$ into $t$ via the
following substitution $\lambda-\mu=-t z$. Note that the integration
range for $t$ is $[0,+\infty)$ and that the integral with respect to
$\lambda$ turns into $(-z)^{m+\tfrac{1}{2} }$. Therefore, the LHS of
\eqref{thimble-osc} is
\ben
(-z)^{m+\tfrac{1}{2} }
\int_{u_Q}^\infty \, e^{\mu/z} \,
\int_{\alpha_{Q,\mu}} \frac{\omega}{df}\, d\mu.
\een
Recalling again the Fubini's theorem we get that the above integral
coincides with the RHS of \eqref{thimble-osc}.
\qed

\subsection{Mirror symmetry for the second structure connection}
\label{sec:ms-2nd_sc}
We will prove that the thimble integrals can be used to construct a
solution to the second structure connection.
\begin{lemma}\label{le:thimble-period}
  Suppose that $Q=e^{2\pi\sqrt{-1} s}q$ and $\lambda=
  e^{2\pi\sqrt{-1} s/3}r$, where $q>0$ and $r\geq 3 q^{1/3}$ are real
  numbers. Then there exists a constant vector $E\in H$ independent of
  $Q$ and $\lambda$ such that
  \beq\label{thimble-period}
  \left(
    I^{(-m-1)}(Q,\lambda) E,\phi_i
  \right) = (-Q\partial_Q)^{i-1} \mathcal{I}^{(-m-i+1)}(Q,\lambda),
  \quad 1\leq i\leq 3
  \eeq
  for all integers $m\geq 2$. 
\end{lemma}
\proof
Let us denote by $L(Q,\lambda)$ the $H$-valued function on $Q$ and
$\lambda$ defined uniquely in such a way that $(L(Q,\lambda),\phi_i)$
coincides with the RHS of \eqref{thimble-period} for all
$i=1,2,3$. Let us denote by $A(Q):=p\bullet_t\in \operatorname{End}(H)$ the linear
operator of quantum multiplication by $t$, where $t=t_2 p$ and
$t_2=\log Q = \tfrac{2\pi\sqrt{-1} s}{3}+\log q$. Since
$I^{(-m-1)}(t,\lambda)$ is a fundamental solution to the second structure
connection, the restriction of $I^{(-m-1)}(t,\lambda)$ to
$t_1=t_3=0$, $t_2=\log Q$ is a fundamental solution to the following system of
equations
\beqa\label{L-q_der}
Q\partial_Q \, I(Q,\lambda) & = &  -A(Q)\, \partial_\lambda\,
I(Q,\lambda) \\
\label{L-homog}
(\lambda\partial_\lambda + 3 Q\partial_Q) \, I(Q,\lambda) & = &
(\theta+m+\tfrac{1}{2})\, I(Q,\lambda).
\eeqa
In order to prove formula \eqref{thimble-period}, it is sufficient to
prove that the vector valued function $L(Q,\lambda)$ is a solution to
the above system of equations.

Recalling the definition of $\mathcal{I}^{(-m)}(Q,\lambda)$ we get
that it has the following scaling symmetry
\ben
\mathcal{I}^{(-m)}(Q c^3,\lambda c) =
c^{m-\tfrac{1}{2}}
\mathcal{I}^{(-m)}(Q,\lambda),
\een
for every real number $c>0$. Differentiating in $c$ and setting $c=1$
yields the following differential equation
\ben
(\lambda \partial_\lambda + 3 Q\partial_Q) \, \mathcal{I}^{(-m)}(Q,\lambda) =
\Big(m-\tfrac{1}{2} \Big) \, \mathcal{I}^{(-m)}(Q,\lambda).
\een
This implies that
\ben
(\lambda \partial_\lambda + 3 Q\partial_Q) \, (L(Q,\lambda),\phi_i) =
\Big(m+i-\tfrac{3}{2} \Big) \,  (L(Q,\lambda),\phi_i).
\een
On the other hand, since $\theta(\phi_i) = (2-i)\phi_i$ we have 
\ben
\Big(m+i-\tfrac{3}{2} \Big)\, \phi_i =
\Big(-\theta+m+\tfrac{1}{2} \Big) \, \phi_i  =
\Big(\theta+m+\tfrac{1}{2} \Big)^T \, \phi_i .
\een
It follows that $L(Q,\lambda)$ satisfies \eqref{L-homog}. Note that
$\partial_\lambda \mathcal{I}^{(-m-1)}(Q,\lambda) =
\mathcal{I}^{(-m)}(Q,\lambda)$ (see formula
\eqref{I-fubini}). Recalling the definition of $L(Q,\lambda)$ we get
\ben
-Q\partial_Q (L,\phi_i) = \partial_\lambda (L,\phi_{i+1}) =
\partial_\lambda (L,A(Q) \phi_{i}) = \partial_\lambda (A(Q)L,\phi_{i})  
\een
for $i=1,2$, where we used that $A(Q)$ is self-adjoint with respect to
the Poincare pairing. Since $A(Q)\phi_3 = Q \phi_1$ in order to prove
that $L(Q,\lambda)$ is a solution to \eqref{L-q_der} it remains only
to prove that
\beq\label{scalar-de}
(-Q\partial_Q)^3 \, \mathcal{I}^{(-m)}(Q,\lambda) =
Q \, \partial_\lambda^3 \, \mathcal{I}^{(-m)}(Q,\lambda),\quad m\geq 4.
\eeq
Note that $-Q\partial_Q \mathcal{I}^{(-m)}(Q,\lambda) $ is equal to 
\ben
\partial_\lambda \int_{\alpha_{Q,\lambda}}
\frac{(\lambda-f(x,Q))^{m-\tfrac{1}{2}} }{ \Gamma(m+\tfrac{1}{2} ) }\,
\frac{Q}{x_1x_2}\, \omega =
\partial_\lambda \int_{\alpha_{Q,\lambda}}
\frac{(\lambda-f(x,Q))^{m-\tfrac{1}{2}} }{ \Gamma(m+\tfrac{1}{2} ) }\,
x_1 \, \omega
\een
where we used the identity
\ben
d(\lambda-f(x,Q))\wedge \frac{d x_2}{x_2} = -x_1 \omega + \frac{Q}{x_1x_2}\, \omega 
\een
and integration by parts. Similarly,
\ben
(-Q\partial_Q)^2 \mathcal{I}^{(-m)}(Q,\lambda)  =
\partial^2_\lambda \int_{\alpha_{Q,\lambda}}
\frac{(\lambda-f(x,Q))^{m-\tfrac{1}{2}} }{ \Gamma(m+\tfrac{1}{2} ) }\,
x_1 x_2\, \omega.
\een
Differentiating the above formula with $-Q\partial_Q$ yields
\eqref{scalar-de}.
\qed
\subsection{Proof of Proposition \ref{prop:refl-vec}}
Let us fix an analytic
branch of $S(Q,z)$ for all $Q$ sufficiently close to $Q^\circ:=1$. Recall that 
$\lambda^\circ$ is a sufficiently large real number. Let us fix also
$m$ to be a sufficiently large positive integer number (e.g. $m\geq 2$ would
work).  If  $\lambda$ is sufficiently close to $\lambda^\circ$, then
we fix a branch of $I^{(-m-1)}(Q^\circ,\lambda) $ -- for example pick
the principal branch of $\log \lambda$, then this would determine
$\widetilde{I}^{(n)}(\lambda)$ for all $\lambda$ close to
$\lambda^\circ$ and hence $I^{(-m-1)}(Q^\circ,\lambda) $ is also
uniquely determined.

Recalling Lemma \ref{le:thimble-period} we get that there exist a
vector $E$ such that formula \eqref{thimble-period} holds with
$Q=Q^\circ$. Since the RHS of \eqref{thimble-period} is an integral
over a Lefschetz thimble, we get that it vanishes when $\lambda$ approaches
$u_{Q^\circ}$ along the line segment
$(u_{Q^\circ},\lambda^\circ]$. Therefore, since $u_{Q^\circ} = 3
=u_1^\circ$ and the line segment $(u_1^\circ,\lambda^\circ]$ is
precisely the path $\gamma_1$, we get that $E$ is proportional to the
reflection vector corresponding to the path $\gamma_1$. We claim that
$E=\sqrt{2\pi} \Psi(\mathcal{O})$. Since the intersection pairing
$(\Psi(\O)|\Psi(\O)) = 2$ this would prove that $\O$ is the reflection
vector corresponding to the path $\gamma_1$ as claimed.

Let us prove that $E=\sqrt{2\pi} \Psi(\O)$. Let us apply to formula
\eqref{thimble-period} the Laplace transform
$\int_{u_{Q^\circ}}^\infty d\lambda e^{\lambda/z}$. Recalling
\eqref{thimble-osc} we get
\ben
\int_{u_{Q}}^\infty e^{\lambda/z}
\left(I^{(-m-1)}(Q,\lambda)E,\phi_i
  \right)
  d\lambda =(-Q\partial_Q)^{i-1} (-z)^{m+i-\tfrac{1}{2}}
  \int_{\alpha_Q} e^{f(x,Q)/z}\omega, 
\een
where we allow $Q$ to be a real deformation of
$Q^\circ$. Recalling formula \eqref{central-charge} with $s=0$ and using the
quantum differential equation $-Q\partial_Q S(Q,z) = (-z)^{-1}
p{\bullet} S(Q,z)$ we get that the
RHS of the above identity is
\ben
\sqrt{2\pi} (-z)^{m+\tfrac{3}{2}} \left(
  S(Q,z) (-z)^\theta (-z)^\rho \Psi(\O),\phi_i\right).
\een
In other words we proved that
\beq\label{aux-2}
\int_{u_{Q}}^\infty e^{\lambda/z}
I^{(-m-1)}(Q,\lambda)E d\lambda =
\sqrt{2\pi} (-z)^{m+\tfrac{3}{2}} 
  S(Q,z) (-z)^\theta (-z)^\rho \Psi(\O).
\eeq
We would like to take the limit $Q\to 0$. Although this limit does not
exists, it is not hard to characterize the singularities of both sides at $Q=0$. Using
the divisor equation we have $S(Q,z)=T(Q,z) e^{P\log Q/z}$, where
$P:=p\cup$ denotes the operator of classical cup product multiplication by
$p$, $T(Q,z)$ is analytic at $Q=0$ and $T(0,z)=1$. Let us write
$T(Q,z)=\sum_{k=0}^\infty T_k(Q) z^{-k}$. We claim that
\beq\label{I-mon_Q}
I^{(-m-1)}(Q,\lambda) = \sum_{k=0}^\infty T_k(Q) (-1)^k
\widetilde{I}^{(-m-1+k)}(\lambda) e^{-P\log Q}.
\eeq
Indeed, by definition
\beq\label{aux-1}
I^{(-m-1)}(Q,\lambda)  = \sum_{k,l=0}^\infty
\frac{(-1)^{k+l} }{l!}T_k(Q) (\log Q)^l P^l \widetilde{I}^{(-m-1+k+l)}(\lambda).
\eeq
On the other hand, using that $P \theta =(\theta+1) P$ we get
\ben
P^l \widetilde{I}^{(-m-1+k+l)}(\lambda) =
e^{\rho\partial_\lambda\partial_m}\Big(
P^l \frac{\lambda^{\theta+m-k-l+\tfrac{1}{2}} }{
  \Gamma(\theta+m-k-l+\tfrac{3}{2})} \Big) =
\widetilde{I}^{(-m-1+k)}(\lambda) P^l.
\een
Substituting this formula in \eqref{aux-1} and summing over all $l$ we
get exactly \eqref{I-mon_Q}. Furthermore,
\ben
e^{P\log Q/z} (-z)^\theta (-z)^\rho  =
(-z)^\theta (-z)^{-\operatorname{ad_\theta}}(e^{P\log Q/z}) (-z)^\rho =
(-z)^\theta (-z)^\rho e^{-P\log Q},
\een
where we used that $-\operatorname{ad_\theta}(P)=P$. Now it is clear
that both sides are polynomials in $\log Q$ of degree 2, whose
coefficients take values in the ring of $H$-valued convergent power series in
$Q$. Comparing  the coefficients in front of
$(\log Q)^0$ and passing to the limt $Q\to 0$ we get
\beq\label{aux-3}
\int_{0}^\infty e^{\lambda/z}
\widetilde{I}^{(-m-1)}(\lambda)E d\lambda =
\sqrt{2\pi} (-z)^{m+\tfrac{3}{2}} 
  (-z)^\theta (-z)^\rho \Psi(\O).
\eeq
The LHS of \eqref{aux-3} is by definition
\ben
\int_0^\infty e^{\lambda/z} e^{\rho\partial_\lambda\partial_m} \Big(
\frac{\lambda^{\theta+m+\tfrac{1}{2}} }{\Gamma(\theta+m+\tfrac{3}{2}) }
\Big) E d\lambda = \sum_{k=0}^\infty \frac{\rho^k(-z)^{-k}\partial_m^k}{k!}
\int_0^\infty e^{\lambda/z}
\frac{\lambda^{\theta+m+\tfrac{1}{2}} }{\Gamma(\theta+m+\tfrac{3}{2}) }
d\lambda E.
\een
Using the substitution $\lambda=-t z$, we get that the integral
in the above formula is $(-z)^{\theta+m+\tfrac{3}{2}}$ and hence the
LHS of \eqref{aux-3} takes the form
\ben
(-z)^{m+\tfrac{3}{2} -\rho/z} (-z)^\theta E = (-z)^\theta
(-z)^{\rho+m+\tfrac{3}{2} } E.
\een
Comparing with the RHS of \eqref{aux-3} we get $E=\sqrt{2\pi}
\Psi(\O)$.

The rest of the proposition is easy to complete. Let us look at
formula \eqref{thimble-period} for $Q=e^{2\pi\sqrt{-1}s}Q^\circ$,
$\lambda=e^{2\pi\sqrt{-1} s/3}\lambda^\circ$ and decrease $s$
continuously from $s=0$ to $s=-1$. The vector $E=\sqrt{2\pi} \Psi(\O)$
does not change, so according to formula \eqref{I-mon_Q} the LHS will
be transformed into $(I^{(-m-1)}(Q^\circ,\lambda^\circ \zeta^{-1})
e^{2\pi\sqrt{-1} p}E,\phi_i)$, which is the same as the analytic
continuation of $(I^{(-m-1)}(Q^\circ,\lambda^\circ)
e^{2\pi\sqrt{-1} p}E,\phi_i)$ along the arc $\lambda=e^{2\pi\sqrt{-1}
  s/3}\lambda^\circ$, $-1\leq s\leq 0$. On the other hand, on the RHS of
\eqref{thimble-period} the only change will be that in the corresponding
thimble integrals the cycle $\alpha_{Q^\circ,\lambda^\circ}$ will be
transformed to the Lefschetz thimble $\alpha_{Q^\circ
  e^{-2\pi\sqrt{-1}},\lambda^\circ \zeta^{-1}}$. The conclusion is
that $(I^{(-m-1)}(Q^\circ,\lambda) e^{2\pi\sqrt{-1} p}E,\phi_i)$
vanishes as $\lambda$ approaches $u_{Q^\circ
  e^{-2\pi\sqrt{-1}}}=u_{Q^\circ}\zeta^{-1} = u^\circ_2$ first along
the arc $\lambda=e^{2\pi\sqrt{-1}
  s/3}\lambda^\circ$ ($-1\leq s\leq 0$) and then along the line
segment $[\lambda^\circ\zeta^{-1}, u^\circ_2]$, i.e., as $\lambda$
travels along the path $\gamma_2$. Therefore, $e^{2\pi\sqrt{-1}
  p}E=\sqrt{2\pi} e^{2\pi\sqrt{-1}p} \Psi(\O) =\sqrt{2\pi}\Psi(\O(1)) $
is proportional to the reflection vector corresponding to
the path $\gamma_2$. The intersection pairing
$(\Psi(\O(1))|\Psi(\O(1)) )=2$, so  $\Psi(\O(1))$ is a reflection
vector. The argument that $\Psi(\O(2))$ is the reflection vector
corresponding to the path $\gamma_3$ is similar -- one just has to
decrease further $s$ from $-1$ to $-2$. 
\qed

\medskip

Let us sketch the proof of Corollary \ref{cor:refl}. First, let us
point out that our proof of Proposition \ref{prop:refl-vec} implies
that $L^n=\O(n)$ is a reflection vector for all $n\in \ZZ$. Let $\ell_i\in
\pi_1(F^\circ,\lambda^\circ)$ be a  simple loop around $u_i^\circ$
corresponding to the path $\gamma_i$ and $M_i\in W$ be the monodromy
transformation representing $\ell_i$. Let $c\in \pi_1((M\times
\CC)',(t^\circ,\lambda^\circ))$ be the loop that under the projection  $(M\times
\CC)' \to M$ maps to the loop in 
$\pi_1(M,t^\circ)\cong \pi_1(\CC^*,Q^\circ)\cong \ZZ$ that goes once
around $Q=0$ in clockwise direction. Formula \eqref{I-mon_Q} implies
that the monodromy transformation $K\in W$ representing the loop $c$
is multiplication by $e^{2\pi\sqrt{-1} p}$ which via Iritani's map
$\Psi$ corresponds to K-theoretic multiplication by $\O(1)$. 
By definition $W_R=\langle M_1,M_2,M_3\rangle$ and $W=\langle
M_1,M_2,M_3, K\rangle.$

Every path $\gamma\subset F^\circ$ from $\lambda^\circ$ to one of the punctures $u_j^\circ$ $(1\leq j\leq 3)$
determines uniquely up to a sign a reflection vector $E\in R$. Since
$M_i^2=1$ we get that $M_i(E)$ is the reflection
vector corresponding to the path $\gamma\circ \ell_i$. We may assume that
the paths $\gamma_j$ and $\gamma$ coincide in a sufficiently small
neighborhood of the puncture $u_j^\circ$. The composition
$\gamma_j^{-1}\circ \gamma=\ell_{i_1}\circ\cdots \circ
\ell_{i_r}$. Therefore, the reflection vector corresponding to the
path $\gamma = \gamma_j\circ  \ell_{i_1}\circ\cdots \circ
\ell_{i_r}$ is $M_{i_r}\cdots M_{i_1}(L^{j-1})$. Since the path
$\gamma$ is arbitrary we get that $R=W_R\cdot 1\cup W_R\cdot L\cup
W_R\cdot L^2$. Finally, in order to complete the proof of Corollary
\ref{cor:refl} we need only to use that  $W_R$ is a normal subgroup of $W$,
the monodromy transformation $K$ and $W_R$ generate $W$, and
$K(L^{j-1})=L^j\in R$. 

\begin{remark}
One can prove that $R=\{E\in K^0(\PP^2)\ |\ \chi(E)=1\}$. In particular,
the set of reflection vectors is bigger than the set formed by the
K-theoretic classes of the exceptional objects in the bounded derived
category $\D^b\operatorname{Coh}(\PP^2)$. 
\end{remark}

\section{The period map for small quantum cohomology}

The goal of this section is to prove Theorem \ref{t1}. Let us assume
that $t_3=0$ and denote by $Z(Q,\lambda)$ the value of the period map
at the point $(Q,\lambda)\in X_{\rm small}$. We do not use an explicit
notation, but we will always keep in mind that $Z(Q,\lambda)$ depends
on the choice of a reference path. 

\subsection{The monodromy group of the second structure connection}
\label{sec:mon-2nd-con}
Let us sketch the main steps in computing the monodromy group $W$. 
The matrix of the intersection form in the basis (over $\mathbb{Q}$)
$E_1,E_2,E_3$ takes the form 
\beq\label{E-im}
\begin{bmatrix}
0 & 0 & 4\\
0 & -8 & 0 \\
4 & 0 & 0
\end{bmatrix}
\eeq
In other words the only non-vanishing pairings are $(E_1|E_3)=4$ and
$(E_2|E_2)=-8$.
Let us denote by $R_i$ the monodromy transformation of the basis
$E=(E_1,E_2,E_3)$ corresponding to analytic continuation along the
path $\gamma_i$ (i.e. the path that turns $L^{i-1}$ into a reflection
vector). We represent $R_i$ by a matrix such that the monodromy
transformation of the row $E$ is $E\,R_i$. A direct computation yields
\ben
R_1=
\begin{bmatrix}
 0& 0 &-1 \\
 0& 1 & 0 \\
-1& 0 & 0
\end{bmatrix},
\quad 
R_2=
\begin{bmatrix}
-1 & 2 &-1 \\
-2& 3 & -1 \\
-4& 4 & -1
\end{bmatrix},
\quad
R_3=
\begin{bmatrix}
-4 & 4 & -1 \\
-10 & 9 & -2 \\
-25 & 20 & -4
\end{bmatrix}.
\een
The matrices $R_i$ ($1\leq
i\leq 3)$ generate reflection group $W_R$ that can be embedded as a finite index
subgroup of the modular groups as follows. Let 
\beq\label{def:phi}
\phi: \CC^3 \to \operatorname{Sym}^2(\CC^2)
\eeq
be the isomorphism identifying $\CC^3$ with the space of symmetric
quadratic forms on $\CC^2$. More precisely
\ben
\phi(z)(u_1,u_2) = z_1 u_1^2 +z_2 u_1 u_2 + z_3 u_2^2.
\een
The modular group $\Gamma:=\operatorname{PSL}_2(\ZZ)$ acts naturally
on the space of quadratic forms 
\ben
q(u_1,u_2)\mapsto (q\cdot g)(u_1,u_2):= q( au_1+b u_2,c u_1+d u_2),\quad 
g=\begin{bmatrix}
a & b \\
c & d
\end{bmatrix}.
\een
Note that the above action is a {\em right} action: 
$q\cdot (g_1g_2)=(q\cdot g_1)\cdot g_2$ . Let us define a
group homomorphism 
\ben
\rho: \operatorname{PSL}_2(\ZZ) \to \operatorname{SL}_3(\ZZ)
\een
such that $\phi(z)\cdot g = \phi(z\rho(g))$, where $z=(z_1,z_2,z_3)$
is a row vector and the matrix $\rho(g)$ acts on $z$ via matrix
multiplication from the right. Explicitly
\ben
\rho(g) = 
\begin{bmatrix}
a^2 & 2ab & b^2 \\
ac & ad+bc & bd \\
c^2 & 2cd & d^2
\end{bmatrix},
\quad 
g=\begin{bmatrix}
a & b \\
c & d
\end{bmatrix}
\een
Note that $R_i=-\rho(g_i)$, where
\ben
g_1=
\begin{bmatrix}
0 & 1 \\
-1 & 0
\end{bmatrix},
\quad
g_2=
\begin{bmatrix}
1 & -1 \\
2 & -1
\end{bmatrix},
\quad
g_3=
\begin{bmatrix}
2 & -1 \\
5 & -2
\end{bmatrix}.
\een
The monodromy group is generated by $R_1,R_2,R_3$, and $K$, where $K$
is the monodromy transformation corresponding to the analytic
continuation in the $Q$-plane along a loop around $Q=0$ in clockwise
direction. Recall that $K$ coincides with the operator of
K-theoretic multiplication by $\O(1)$ (see the discussion after the
proof of Proposition \ref{prop:refl-vec} in Section \ref{sec:ms-2nd_sc}).
Therefore $E=(E_1,E_2,E_3)$ transforms into $E\,K$ where  
\beq\label{matr-K}
K=\begin{bmatrix}
1 & 0& 0 \\
1 & 1 & 0 \\
1 & 2 &1 
\end{bmatrix} =\rho(\kappa),\quad 
\kappa=
\begin{bmatrix}
1 & 0 \\
1 & 1
\end{bmatrix} .
\eeq
The relation between $W_R$, $W$, and the modular group can be
described as follows. The matrices $g_1\kappa$ and $g_1$ have orders respectively $3$ and
$2$ and we have $\operatorname{PSL}_2(\ZZ) = \langle g_1\kappa\rangle * \langle
g_1\rangle. $ Using this presentation of the modular group we define the characters
\ben
\chi_2,\chi_3 :\operatorname{PSL}_2(\ZZ) \to \CC^*,
\een
such that 
\ben
\chi_2(g_1\kappa)=\chi_3(g_1)=1,\quad \chi_2(g_1)=-1,\quad
\chi_3(g_1\kappa)=\zeta. 
\een
Let us define a group homomorphism 
\beq\label{mod-mon}
\operatorname{PSL}_2(\ZZ)\times \{\pm 1\} \to \operatorname{GL}(\CC^3),
\quad 
(g,\sigma)\mapsto \sigma\, \chi_2(g) \,\rho(g).
\eeq
Using that 
\ben
R_2 = K R_1 K^{-1}, \quad R_3= K^2 R_1 K^{-2}.
\een
we get that the image of the map \eqref{mod-mon} is the monodromy
group $W$. It is not very difficult to check that the map is also
injective, so it gives a group isomorphism $W\cong
\operatorname{PSL}_2(\ZZ) \times \{\pm 1\} $.  Finally, using that the
reflection group $W_R$ is generated by $R_1,R_2,R_3$ we get that the map
\ben
\operatorname{Ker}(\chi_3)\to W_R,\quad g\mapsto \chi_2(g) \,\rho(g)
\een
is a group isomorphism.

\subsection{Quadratic relation}

The second structure connection for small quantum cohomology takes the form
\ben
Q\partial_Q I^{(-1)} (Q,\lambda) = -\partial_\lambda (P\bullet  I^{(-1)} (Q,\lambda) )
\een
and 
\ben
(\lambda\partial_\lambda + 3 Q\partial_Q)  I^{(-1)} (Q,\lambda)  =
(\theta+1/2)  I^{(-1)} (Q,\lambda) .
\een
Using these equations, we get that the period map satisfies the
following differential equations
\ben
\Big( (Q\partial_Q)^3 - Q\partial_\lambda^3 \Big) Z(Q,\lambda)=0
\een
and
\ben
(\lambda\partial_\lambda + 3 Q\partial_Q) Z(Q,\lambda) = -\frac{1}{2} Z(Q,\lambda),
\een
where we used that in small quantum cohomology $P\bullet P\bullet P=
Q$ and that $\theta(1)=1$. 

The second equation implies that the period map has the form 
\ben
Z(Q,\lambda) = Q^{-1/6} z(x),\quad x:=\frac{\lambda^3}{27 Q}
\een
while the first equation implies that the vector valued function $z(x)$ is a solution
to the hypergeometric equation of type $(3,2)$ defined by the
differential operator
\beq\label{hge:3F2}
D (D-\rho_1)(D-\rho_2) - x (D+\alpha_1) (D+\alpha_2) (D+\alpha_3),
\eeq
where $D=x\partial_x$, $\rho_1=\frac{1}{3}$, $\rho_2=\frac{2}{3}$, and
$\alpha_1=\alpha_2=\alpha_3=\frac{1}{6}.$
\begin{lemma}\label{le:quadr-rel}
The image of the period map 
$
Z_{\rm small}: \widetilde{X}_{\rm small} \to \CC^3
$
is contained in the quadratic cone $Z_2^2=4Z_1Z_3$.
\end{lemma}
\proof
The equation of the quadratic cone coincides with 
\ben
\sum_{i,j=1}^3 \eta^{ij} Z_i Z_j =0
\een
where $\eta^{ij}$ are the entries of the matrix inverse to the
matrix $\eta$ whose entries are the intersection numbers
$\eta_{ij}:=(E_i|E_j)$ (see formula \eqref{E-im}). Let us denote by
$I^{(n)}$ the matrix whose $(i,j)$-entry is given by
$(I^{(n)}_{E_i}(t,\lambda), \phi_j)$. Using the equations of the
second structure connection we get 
\ben
I^{(0)} (\lambda-E\bullet) = I^{(-1)} \Big(-\theta+\frac{1}{2}\Big),
\een
where 
\ben
E\bullet=\begin{bmatrix}
0 & 0 & 3 Q \\
3 & 0 & 0\\
0 & 3 & 0
\end{bmatrix}.
\een
On the other hand, the entries of the intersection pairing are 
\ben
\eta_{ij}= (I^{(0)}_{E_i},(\lambda-E\bullet) I^{(0)}_{E_j}) = 
\sum_{k,\ell=1}^3 (I^{(0)}_{E_i},\phi_k)g^{k\ell} (
I^{(0)}_{E_j},(\lambda-E\bullet)\phi_\ell), 
\een
where $g^{k\ell}$ are the entries of the matrix
inverse to the matrix of the Poincare pairing.
Therefore
\ben
\eta = I^{(0)} g^{-1} (\lambda-E\bullet)^T (I^{(0)})^T.
\een
Using that $Z=(Z_1,Z_2,Z_3)= (I^{(-1)} e_1)^T$ we get 
\ben
Z\eta^{-1} Z^T  = e_1^T (I^{(-1)})^T \eta^{-1} I^{(-1)}  e_1 .
\een
Recalling the formula for $\eta$ from above we get 
\ben
(I^{(-1)})^T \eta^{-1} I^{(-1)} = 
((I^{(0)})^{-1} I^{(-1)})^T   ((\lambda-E\bullet)^T)^{-1} g 
(I^{(0)})^{-1} I^{(-1)} .
\een
On the other hand
\ben
(I^{(0)})^{-1} I^{(-1)} = (\lambda-E\bullet) (-\theta+1/2)^{-1}.
\een
Therefore $Z\eta^{-1} Z^T$ is the (1,1)-entry of the matrix 
\ben
(-\theta+1/2)^{-1} \, g\, (\lambda-E\bullet) (-\theta+1/2)^{-1} = 
\begin{bmatrix}
0 & 12 & -4\lambda/3 \\
12 & 4\lambda & 0 \\
-4\lambda/3 & 0 & -4Q/3
\end{bmatrix}.
\qed
\een

\subsection{Connection Formula}\label{sec:conn-formula} 
The differential equation \eqref{hge:3F2} has the following basis of
solutions near $x=\infty$.
\ben
z^\infty_1(x):= \sum_{n=0}^\infty a_n x^{-n-\frac{1}{6}},
\een
where the coefficients $a_n$ are defined by 
\ben
a_0=1,\quad a_{n+1}=\frac{(n+1/6)(n+3/6)(n+5/6)}{(n+1)^3},\quad n\geq 0.
\een
Note that $z^\infty_1(x)$ coincides with the generalized hypergeometric function
\ben
{}_3 F_2\left[
\begin{matrix}
1/6\quad 3/6\quad 5/6\\
1\quad 1
\end{matrix}
\, ;x^{-1}\right]\, x^{-1/6}.
\een
The second solution is 
\ben
z^\infty_2(x)=\sum_{n=0}^\infty a_n x^{-n-\frac{1}{6}} \left( \log x - b_n\right),
\een
where the constants $b_n$ $(n\geq 0)$ are defined by $b_0=0$ and
\ben
b_{n+1} = b_n + \frac{1}{n+\frac{1}{6}} +
\frac{1}{n+\frac{3}{6}}+\frac{1}{n+\frac{5}{6}} -
\frac{3}{n+1}.
\een
Finally, the third solution is given by 
\ben
z^\infty_3(x) = \sum_{n=0}^\infty 
a_n x^{-n-\frac{1}{6}} \left( (\log x - b_n)^2+c_n\right),
\een
where the constants $c_n$ $(n\geq 0)$ are defined by $c_0=0$ and 
\ben
c_{n+1}= c_n + \frac{3}{(n+1)^2} - 
\frac{1}{\left(n+\frac{3}{6}\right) \left(n+\frac{5}{6}\right)} -
\frac{1}{\left(n+\frac{1}{6}\right) \left(n+\frac{5}{6}\right)   } -
\frac{1}{\left(n+\frac{1}{6}\right)\left(n+\frac{3}{6}\right)}.
\een
Let us find the transition matrix $C^\infty$ defined by 
\beq\label{conn-formula}
z(x)=z^\infty(x) \, C^\infty,
\eeq
where
\ben
z(x)=(z_1(x),z_2(x),z_3(x)) ,
\quad
z^\infty(x)=(z^\infty_1(x),z^\infty_2(x),z^\infty_3(x)).
\een
To begin with, note that the analytic continuation of $Z(Q,\lambda)$
along a clock-wise loop around $0$ in the $Q$-plane corresponds to
analytic continuation of $z(x)$ along an anticlock-wise loop around $0$
and $1$ in the $x$-plane. Since the analytic continuation transforms
$Z(Q,\lambda)$ into $Z(Q,\lambda)\, K$ where $K$ is the matrix
\eqref{matr-K} we get that $z(x)=Q^{1/6} Z(Q,\lambda)$ transforms into
\ben
\zeta^{-1/2}\, z(x)\, K.
\een
On the other hand, the analytic continuation of $z^\infty(x)$ is 
\ben
\zeta^{-1/2} \, z^\infty(x)\, K^\infty,\quad
K^\infty:=
\begin{bmatrix}
1 & 2\pi\sqrt{-1} & (2\pi\sqrt{-1})^2 \\
0 & 1 & 4\pi\sqrt{-1}\\
0 & 0 & 1 
\end{bmatrix}.
\een
Therefore we have the following relation 
\ben
C^\infty K = K^\infty C^\infty.
\een
Let us denote by $C^\infty_i$ the $i$-th column of $C^\infty$. Then
comparing the columns in the above relation we get
\ben
C^\infty_2=\Big(K^\infty-1 -\frac{1}{2} (K^\infty-1)^2\Big)\, C^\infty_1
= \begin{bmatrix}
0 & 2\pi\sqrt{-1} & 0 \\
0 & 0 & 4\pi\sqrt{-1}\\
0 & 0 & 0 
\end{bmatrix}\, C^\infty_1
\een
and 
\ben
C^\infty_3=\frac{1}{2} (K^\infty-1)^2\, C^\infty_1 = 
\begin{bmatrix}
0 & 0 & (2\pi\sqrt{-1})^2 \\
0 & 0 & 0\\
0 & 0 & 0 
\end{bmatrix}\, C^\infty_1
\een
Therefore, we need just to find the first column of $C^\infty_1$, i.e.,
the coefficients in the relation 
\beq\label{z1}
z_1(x) = C^\infty_{11} z^\infty_1(x)+C^\infty_{21}
z^\infty_2(x)+C^\infty_{31} z^\infty_3(x) .
\eeq
The leading order term in the expansion of the RHS of \eqref{z1} at $x=\infty$ is
precisely
\ben
C^\infty_{11} x^{-1/6}+
C^\infty_{21}x^{-1/6}\log x+
C^\infty_{31} x^{-1/6}(\log x)^2.
\een
By definition
$z_1(x)=Q^{1/6} (I^{(-1)}_{E_1}(Q,\lambda),\mathbf{1})$,
where $E_1=1+2L-L^2 $. Recalling formula \eqref{I-mon_Q} we get that
the leading order term of the expansion of the LHS of \eqref{z1}
coincides with the leading order term of 
\ben
Q^{1/6} \left(
  \widetilde{I}^{(-1)}(\lambda)\, e^{-p\log Q}\,  \Psi(1+2L-L^2),
  \mathbf{1}\right).
\een
Recalling the definitions of the period
$\widetilde{I}^{(-1)}(\lambda)$ and the Iritani's map $\Psi$
we get that the above expression is precisely
\ben
\frac{1}{\sqrt{6}\pi} \left( \log x + \log 1728\right)^2 x^{-1/6}.
\een
Therefore,
\ben
C^\infty_{11}=\frac{(\log 1728)^2}{\sqrt{6}\pi},\quad 
C^\infty_{21}=\frac{2\log 1728}{\sqrt{6}\pi},\quad 
C^\infty_{31}=\frac{1}{\sqrt{6}\pi}
\een
and we get
\ben
C^\infty=
\frac{1}{\sqrt{6}\pi}\,
\begin{bmatrix}
(\log 1728)^2 & 2(2\pi\sqrt{-1})\,\log 1728 & (2\pi\sqrt{-1})^2 \\
2\log 1728  & 4\pi\sqrt{-1} & 0\\
1 & 0 & 0 
\end{bmatrix}.
\een

\subsection{Symmetric square of a hypergeometric equation}

We are going to prove that the solution space of the generalized hypergeometric
equation \eqref{hge:3F2} has a basis of the form
\beq\label{sym2}
z_1=\frac{u^2}{2} ,\quad z_2=uv,\quad z_3=\frac{v^2}{2},
\eeq
where $u$, and $v$ is a basis of solutions for the differential
equation defined by the differential operator
\beq\label{hge:2F1}
D^2-\frac{2+x}{6(1-x)} \, D -\frac{x}{144(1-x)},\quad D:=x\partial_x.
\eeq
Note that the above operator defines a differential equation
equivalent to the classical hypergeometric equation defined by the
differential operator
\ben
(1-x) x\partial_x^2 +(c-(a+b+1)x)\partial_x -ab
\een
with $a=b=\frac{1}{12}$ and $c=\frac{2}{3}$.

In order to compute the symmetric square of \eqref{hge:2F1} we have to
find 3 functions $p_i(x)$ $(1\leq i\leq 3)$ such that 
\ben
D^3 z_i = p_2(x) D^2 z_i +p_1(x) D z_i +p_0(x) z_i=0,\quad 1\leq i\leq 3,
\een
where $z_i$ are given by \eqref{sym2} with $\{u,v\}$ a basis of
solutions to \eqref{hge:2F1}. Using that $u$ and $v$ are solutions to
\eqref{hge:2F1} we can express the equations of the above linear
system for $p_0,p_1$, and $p_2$ as differential polynomials in $u$ and
$v$. After a direct computation we get
\ben
c_1(Du)^2 +c_2 (u Du) + c_3 u^2 & = & 0,\\
2c_1 Du Dv+c_2 (uDv+vDu) +2c_3uv & = & 0,\\
c_1 (Dv)^2 +c_2 (v Dv)+ c_3 v^2 =0,
\een
where
\ben
c_1 & =&p_2-3\alpha\\
c_2 & =&p_1+\alpha p_2-4\beta-\alpha^2-D\alpha\\
c_3 & =&p_0/2-\alpha\beta-D\beta +p_2\beta,
\een
and
\ben
\alpha:=\frac{2+x}{6(1-x)},\quad \beta:= \frac{x}{144(1-x)}.
\een
From here we get that $c_1=c_2=c_3=0$, so 
\ben
p_2=3\alpha = \frac{2+x}{2(1-x)},
\een
\ben
p_1=4\beta-2\alpha^2 + D\alpha = -\frac{8-3x}{36(1-x)},
\een
and 
\ben
p_0=-4\alpha\beta +2D\beta=\frac{x}{216(1-x)}.
\een
Finally, it remains only to verify that the differential operator
$(1-x)^{-1}(D^3-p_2 D^2-p_1 D-p_0)$ coincides with \eqref{hge:3F2}.  

The classical hypergeometric equation \eqref{hge:2F1} has a basis of
solutions near $x=\infty$ of the following form
\ben
v^\infty(x)= \sum_{n=0}^\infty v_n x^{-n-\frac{1}{12}}
\een
and 
\ben
u^\infty(x) = \sum_{n=0}^\infty v_n \left(\log x - u_n\right) x^{-n-\frac{1}{12}},
\een
where the constants $u_n,v_n$ $(n\geq 0)$ are defined by 
\ben
v_0=1,\quad v_{n+1} = \frac{\left(n+\frac{1}{12}\right)
  \left(n+\frac{5}{12}\right) }{(n+1)^2} v_n,
\een
and 
\ben
u_0=0,\quad u_{n+1} = u_n+\frac{1}{n+\frac{1}{12} } +
\frac{1}{n+\frac{5}{12} } - \frac{2}{n+1}. 
\een
Comparing the leading coefficients in the Laurent series expansion
near $x=\infty$ we get 
\ben
z^\infty_1(x) = (v^\infty(x))^2,\quad 
z^\infty_2(x) = u^\infty(x)v^\infty(x),\quad
z^\infty_3(x) = (u^\infty(x))^2.
\een
Recalling the quadratic relation $Z_2^2=4Z_1Z_3$, we get that the
period map $Z_{\rm small}$ can be expressed in terms of the 
hypergeometric functions $u^\infty$ and $v^\infty$ as follows:
\beq\label{z-hge}
Z_{\rm small}(Q,\lambda) =
(\tau(x)^2 z_3(x) Q^{-1/6}, -2\tau(x) z_3(x) Q^{-1/6}, z_3(x) Q^{-1/6}),
\eeq
where $x=\tfrac{\lambda^3}{27Q}$,
\beq\label{z3}
z_3(x)= \frac{1}{\sqrt{6}\pi} \, (2\pi\sqrt{-1})^2 (v^\infty(x))^2
\eeq
and
\beq\label{def:tau}
\tau(x):=-\frac{z_2(x)}{2 z_3(x)} = -\frac{1}{2\pi\sqrt{-1}} \left( 
\log 1728 +\frac{u^\infty(x)}{v^\infty(x)}\right),
\eeq
where \eqref{z3} and the second equality in \eqref{def:tau} are
obtained by formula \eqref{conn-formula}.

\subsection{The Schwarz map}
The map \eqref{def:tau} is known as the Schwarz map for the hypergeometric
equation \eqref{hge:2F1} (with $a=b=\tfrac{1}{12}$,
$c=\tfrac{2}{3}$). The exponents of \eqref{hge:2F1} at $x=0,1$, 
and $\infty$ are respectively $|c-a-b|=\tfrac{1}{2}$,
$|1-c|=\tfrac{1}{3}$, and $|a-b|=0=\tfrac{1}{\infty}$. Therefore, the
holomorphic branches of $\tau(x)$ define maps whose images are
hyperbolic triangles that define a triangulation of the upper-half
plane. Let us work out the precise structure of the triangulation.

Let us fix $x^\circ=\tfrac{(\lambda^\circ)^3}{27}$ as a reference point
and define the holomorphic branch of $\tau(x)$ near $x=x^\circ$ to be
the one for which $\log x^\circ$ is a real number. If $x$ is in a
neighborhood of $0$, then the hypergeometric equation admits the basis of
solutions
\ben
v^0(x) & = &  {}_2F_1(\tfrac{1}{12},\tfrac{1}{12};\tfrac{2}{3}; x)\\
u^0(x) & = & x^{1/3} {}_2F_1(\tfrac{5}{12},\tfrac{5}{12};\tfrac{4}{3}; x).
\een
If $x$ is in a neighborhood of $x=1$, then a basis of solution is given
by
\ben
v^1(x) & = &  {}_2F_1(\tfrac{1}{12},\tfrac{1}{12};\tfrac{1}{2};1-x)\\
u^1(x) & = & (1-x)^{1/2} {}_2F_1(\tfrac{7}{12},\tfrac{7}{12};\tfrac{3}{2};1-x).
\een
The values of the above functions are defined to be real for all real
$x\in (0,1)$ while for other $x\in \CC\setminus{\{0,1\}}$ the values
are specified only after we choose a reference path from $x$ to the
interval $(0,1)$.

The general theory of the Gauss hypergeometric equations provides
explicit formulas for the two linear transformations that relate the
bases $(u^\infty,v^\infty)$ and $(u^0,v^0)$ and the bases
$(u^0,v^0)$ and $(u^1,v^1)$ (see \cite{BE}). The first linear
transformation takes the form 
\ben
v^0(x) & = & 
\frac{
  \Gamma(\tfrac{2}{3}) e^{  \tfrac{\pi\sqrt{-1}}{12} }} {
    \Gamma(\tfrac{1}{12})\Gamma(\tfrac{7}{12})}
  \left(
    u^\infty(x) + v^\infty(x) \Big(-\pi\sqrt{-1} +\log
    (1728)+\pi\sqrt{3}\Big)\right)\\
  u^0(x) & = & 
\frac{
  \Gamma(\tfrac{4}{3}) e^{  \tfrac{5\pi\sqrt{-1}}{12} }} {
    \Gamma(\tfrac{5}{12})\Gamma(\tfrac{11}{12})}
  \left(
    u^\infty(x) + v^\infty(x) \Big(-\pi\sqrt{-1} +\log
    (1728)-\pi\sqrt{3}\Big)\right),
\een
where $u^0(x)$ and $v^0(x)$ are analytically extended to a
neighborhood of $x=x^\circ$ along a path consisting of the arc
$x=x^\circ e^{\pi\sqrt{-1} s}$ ($0\leq s\leq 1$), an interval
$[-x^\circ,-\epsilon]$ such that $0<\epsilon<1$, and the arc
$x=\epsilon e^{\pi\sqrt{-1} s}$ ($-1\leq s\leq 0$). The derivation of
the above formulas uses formulas (18) and (23) from Section 2.1.4 in \cite{BE}
and the following identities
\ben
u_n=\psi(n+\tfrac{1}{12})+ \psi(\tfrac{7}{12}-n) -2\psi(n+1)-
\Big(\psi(\tfrac{1}{12})+\psi(\tfrac{7}{12})-2 \psi(1)\Big),
\een
\ben
u_n=\psi(n+\tfrac{5}{12})+ \psi(\tfrac{11}{12}-n) -2\psi(n+1)-
\Big(\psi(\tfrac{5}{12})+\psi(\tfrac{11}{12})-2 \psi(1)\Big),
\een
\ben
\psi(\tfrac{1}{12})+\psi(\tfrac{7}{12})-2 \psi(1) =-\log(1728)-\pi\sqrt{3},
\een
\ben
\psi(\tfrac{5}{12})+\psi(\tfrac{11}{12})-2 \psi(1) =-\log(1728)+\pi\sqrt{3},
\een
where $\psi(z)=\tfrac{d}{dz} \log \Gamma(z)$ is the logarithmic
derivative of the Gamma function. For the last two formulas we used a
formula due to Gauss that expresses the rational values $\psi(p/q)$ in
terms of elementary transcendental functions (see formula (29) from
Section 1.7.4 in \cite{BE}). 

Similarly, using formula (1) from Section 2.10 and formulas (5) and (6)
from Section 2.9 in \cite{BE} we get 
\ben
v^0(x) & = & 
\frac{
  \Gamma(\tfrac{1}{2}) \Gamma(\tfrac{2}{3})   } {
  \Gamma(\tfrac{7}{12})^2} \, v^1(x) +
\frac{
  \Gamma(\tfrac{-1}{2}) \Gamma(\tfrac{2}{3})   } {
  \Gamma(\tfrac{1}{12})^2} \, u^1(x) \\
u^0(x) & = & 
\frac{
  \Gamma(\tfrac{1}{2}) \Gamma(\tfrac{4}{3})   } {
  \Gamma(\tfrac{11}{12})^2} \, v^1(x) +
\frac{
  \Gamma(\tfrac{-1}{2}) \Gamma(\tfrac{4}{3})   } {
  \Gamma(\tfrac{5}{12})^2} \, u^1(x) \, ,
\een
for all $x$ in a contractible open neighborhood of the open interval $(0,1)$.  

Put $\tau^i(x)=\tfrac{u^i(x)}{v^i(x)}$ for $i=0,1,\infty$. By
definition $\tau^\infty(x) +\log(1728)= -2\pi\sqrt{-1}
\tau(x)$. Therefore, using the linear relation between
$(u^\infty,v^\infty)$ and $(u^0,v^0)$ we get
\beq\label{tau-tau0}
\tau^0 = k \, \frac{\rho \tau+1}{\tau+\rho},
\quad k :=\frac{
  \Gamma(\tfrac{4}{3}) \Gamma(\tfrac{1}{12})
  \Gamma(\tfrac{7}{12}) }{
  \Gamma(\tfrac{2}{3}) \Gamma(\tfrac{5}{12})
  \Gamma(\tfrac{11}{12}) },
\eeq
where $\rho:=e^{\tfrac{\pi\sqrt{-1}}{3} }=\tfrac{1}{2}
+\tfrac{\sqrt{-3}}{2}$. This formula allows us to determine the image
of the ray $(-\infty,0)$ under the Schwarz map
\eqref{def:tau}. Indeed, since
\beq\label{tau-x}
-2\pi\sqrt{-1}\, \tau(x) = \log(1728\, x) -
\frac{\sum_{n=1}^\infty  u_n v_n x^{-n}   } {
  1+ \sum_{n=1}^\infty v_n x^{-n} }
\eeq
we get that $\operatorname{Re}(\tau(x))=-\tfrac{1}{2}$ and
$\lim_{x\to -\infty} \tau(x)=-\tfrac{1}{2} +\infty \,\sqrt{-1}$.
Using the relation \eqref{tau-tau0}, we get
$\lim_{x\to 0} \tau(x)=-\rho^{-1}$. In other words, the image of
$(-\infty,0)$ is the geodesic line in $\HH$ from
$-\tfrac{1}{2}+\infty\,\sqrt{-1}$ to $-\rho^{-1}$.

Furthermore, let us determine the image of the interval $(1,+\infty)$
under the Schwarz map \eqref{def:tau}. The linear relation between the
bases $(u^0,v^0)$ and $(u^1,v^1)$ yields 
\beq\label{tau-tau1}
\frac{\rho\tau+1}{\tau+\rho} = 
k^{-1} \tau^0(x) =  \frac{
  \alpha- 2\alpha^{-1}\tau^1}{\beta-2\beta^{-1} \tau^1},
\eeq
where
$\alpha=\tfrac{\Gamma(5/12)}{\Gamma(11/12)}$ and 
$\beta=\tfrac{\Gamma(1/12)}{\Gamma(7/12)}$. If $x\to 1$ with real
values, then $\tau^1(x)\to 0$, so using relation \eqref{tau-tau1} we
find that the limit $\lim_{x\to 1}\tau(x)=\sqrt{-1}$. Using
\eqref{tau-x} we get that if $x\in (1,+\infty)$, then
$\operatorname{Re}(\tau(x))=0$ and $\lim_{x\to +\infty} \tau(x) =
+\infty \sqrt{-1}$. The conclusion is that the image of the interval
$(1,\infty)$ is the geodesic line in $\HH$ from $\sqrt{-1}$ to
$+\infty\,\sqrt{-1}$.

Finally, the image of the interval $(0,1)$ under the Schwarz map
\eqref{def:tau} must be a geodesic in $\HH$ connecting $-\rho^{-1}$
and $\sqrt{-1}$. Therefore it must be the arc on the unit circle
connecting $-\rho^{-1}$ and $\sqrt{-1}$. Using again formula
\eqref{tau-x} we fined that the real part of $\tau(\sqrt{-1} x^\circ)$
is $-\tfrac{1}{4}$. Therefore, the Schwarz map \eqref{def:tau} maps the
upper-half plane $\operatorname{Im}(x)>0$ into the interior of the
hyperbolic triangle $\Delta_+$ in $\HH$ with vertices $\infty$,
$-\rho^{-1}$ and $\sqrt{-1}$.

The image of the lower half-plane
$\operatorname{Im}(x)<0$ can be determined by monodromy
reasons. Namely, the values of the Schwarz map above and below the
interval $(1,+\infty)$ vary continuously as we cross the interval,
while the values above and below the intervals $(0,1)$ and
$(-\infty,0)$ make a jump that can be computed if we new the monodromy
of the Schwarz map around respectively the loop around $x=1$ and the
loop around both $x=1$ and $x=0$. Let us compute the monodromy of the
Schwarz map around these two loops.
Using formula \eqref{tau-x} we get that the analytic
continuation of $\tau(x)$ along the circle $x=x^\circ
e^{-2\pi\sqrt{-1} s }$ ($0\leq s\leq 1$) transforms $\tau(x^\circ)$
into $\tau(x^\circ)+1$, i.e., changing the reference path by
pre-composing it with a big loop $\ell_\infty$ around $x=0$ and $x=1$ (with
clock-wise orientation) transforms the holomorphic branch $\tau\mapsto
\tau+1$. Similarly, the analytic continuation along the loop $\ell_1$ that
starts at $x^\circ$ approaches $x=1$ along the real axes, goes around
$x=1$ anti-clockwise along a small circle, and returns back to
$x^\circ$ can be determined by formula \eqref{tau-tau1}. Indeed, such
analytic continuation transforms $\tau^1\mapsto -\tau^1$. Using formula
\eqref{tau-tau1} after a direct computation we find that $\tau\mapsto
-\tau^{-1}$.  Let us summarize the conclusions of our computations.
\begin{proposition}\label{prop:schw}
a) The analytic continuation of the Schwarz map \eqref{def:tau}
defines a representation 
\beq\label{hge-mon}
\pi_1(\CC\setminus{\{0,1\}},x^\circ)\to \operatorname{Aut}(\HH)
\eeq
such that the generators $\ell_1$ and $\ell_\infty$ are mapped
respectively to the transformations $\tau\mapsto -\tau^{-1}$ and
$\tau\mapsto \tau+1$.

b) The image of the upper (resp. lower) half-plane
$\operatorname{Im}(x)\geq 0$ (resp. $\operatorname{Im}(x)\leq 0$) under
the Schwarz map is the hyperbolic triangle $\Delta_+$
(resp. $\Delta_-$) with vertices
$\infty$, $-\rho^{-1}$, $\sqrt{-1}$ (resp. $\infty,\sqrt{-1},\rho$). 
\end{proposition}

Note that the image of the monodromy representation \eqref{hge-mon} is
the modular group $\operatorname{PSL}_2(\ZZ)$ 
and that $\Delta_+\cup \Delta_-$ is its fundamental domain. Therefore
all Schwarz triangles are obtained from $\Delta_+$ and $\Delta_-$  by
the action of $\operatorname{PSL}_2(\ZZ)$. The vertices of the
triangulation that are inside $\HH$ are the two orbits
$\operatorname{PSL}_2(\ZZ)\cdot \sqrt{-1} = \{\tau\in \HH\ |\
E_4(\tau)=0\}$ and $\operatorname{PSL}_2(\ZZ)\cdot \rho = \{\tau\in \HH\ |\
E_6(\tau)=0\}$. The Schwarz map allows us to factorize the universal
covering map $(\CC\setminus{\{0,1\}})^\sim \to \CC\setminus{\{0,1\}}$
as follows
\begin{diagram}
(\CC\setminus{\{0,1\}})^\sim & \rTo^{\tau} &
\HH\setminus{\{E_4(\tau)=0\}\cup \{E_6(\tau)=0\} } & \rTo^{J} &
\CC\setminus{\{0,1\}},
\end{diagram}
where $\tau$ is the Schwarz map and $J$ is the regular covering
corresponding to the kernel of the monodromy representation
\eqref{hge-mon}. 
\begin{lemma}
a) The function $J$ extends to a holomorphic function on $\HH$ and it
coincides with the $J$-invariant, i.e., the unique
$\operatorname{PSL}_2(\ZZ)$-invariant holomorphic function on $\HH$
such that  $J(\sqrt{-1})=1$ and the Fourier series of $J(\tau)$ has a pole
of order 1 at $q=0$, where $q=e^{2\pi\sqrt{-1} \tau}$.

b) The pullback of any holomorphic branch of $(v^\infty(x))^{12}$ to
$\HH$ via the map $J$ extends to an analytic function on $\HH$ and it 
coincides with the modular form $E_4(\tau)^3-E_6(\tau)^2$.
\end{lemma}
\proof
a) Suppose that $\tau\in \HH$ is not a vertex of a Schwarz
triangle and let us put $x=J(\tau)$. Recall that the Schwarz triangles are
biholomorphic to the upper or the lower half-plane and define a
triangulation of $\HH$. Therefore, there exists a reference
path $C$ in $\CC\setminus{\{0,1\}}$ that defines a branch $\tau_C$ of
the Schwarz map in a neighborhood of $x$ such that  $\tau_C(x):=\tau$. The
choice of such a path is not unique, but the branch $\tau_C$ is
independent of the choice of $C$. In other words, locally the map 
$J$ is defined by inverting the Schwarz map.

Let us prove that $J$ is
$\operatorname{PSL}_2(\ZZ)$-invariant. Suppose that $\tau''=g(\tau')$
for some $g\in \operatorname{PSL}_2(\ZZ)$. Pick a point $x'$ and a
reference path $C'$ such that $\tau_{C'}(x')=\tau'$. Let $C$ be a loop
based at $x^\circ$ such that the monodromy representation
\eqref{hge-mon} maps the 
homotopy class of $C$ to $g$. Then $\tau_{C\cdot
  C'}(x')=g(\tau_{C'}(x'))=\tau''$. In particular, $x'=J(\tau'')$. On
the other hand, the identity $\tau_{C'}(x')=\tau' $ implies
$J(\tau')=x'$, so $J(\tau')=J(\tau'')=J(g(\tau'))$. 

Since $J$ is $\operatorname{PSL}_2(\ZZ)$-invariant, in order to prove
that $J$ extends to a holomorphic function on $\HH$, it is enough to
prove that it extends holomorphically in a neighborhood of
$\tau=\sqrt{-1}$ and $\tau=-\rho^{-1}$. Let us give the argument for
$\tau=\sqrt{-1}$. The other case is similar. If $x$ is
sufficiently close to $x=1$ then substituting 
$\tau^1(x)=(1-x)^{1/2}(1+O(1-x))$ in \eqref{tau-tau1} and solving for
$1-x$ in terms of $\tau$, we get that $1-x = f(\tau(x))$ for some function
$f(\tau)=c_2(\tau-\sqrt{-1})^2 + c_3(\tau-\sqrt{-1})^3+\cdots $
holomorphic near $\tau=\sqrt{-1}$.  Since $x=J(\tau(x))=1-f(\tau(x))$
for all $x$ in a neighborhood of $x=1$, we conclude that the
holomorphic functions $1-f(\tau)$ and $J(\tau)$ agree on a dense
subset of a neighborhood of $\tau=\sqrt{-1}$, so they must coincide
identically. This proves that $J(\tau)$ is holomorphic at
$\tau=\sqrt{-1}$ and moreover that $J(\sqrt{-1})=1$.

It remains only to prove that the Fourier series of $J(\tau)$ has a
pole of order 1.  If $x$ is
sufficiently close to $\infty$ in $\CC\setminus{\{0,1\}}$, then let us
exponentiate relation \eqref{tau-x}. We get the following formula:
\beq\label{q-x}
q=\frac{x^{-1}}{1728}\, 
\exp\Big(
\frac{
\sum_{m=1}^\infty u_m v_m x^{-m} 
}{
\sum_{m=1}^\infty u_m x^{-m}}\Big),\quad q:=e^{2\pi\sqrt{-1}\tau}.
\eeq
Formula \eqref{q-x} can be solved for $x^{-1}$ in terms of $q$. We get
that $x= 1728 q^{-1} + O(1)$. This completes the proof of part a).

b) Let us first check that the pullback of $(v^\infty)^{12}$ is 
holomorphic on the complement in $\HH$ of the vertices of the Schwarz
triangles. The pullback of $v^\infty$ to the universal cover is an
analytic function and the complement of the vertices of the Schwarz
triangles is a quotient of the universal cover by the kernel of the
monodromy representation \eqref{hge-mon}. Therefore, we have to verify
that $(v^\infty(x))^{12}$ is invariant under the analytic continuation
of every loop $C$ whose homotopy class is in the kernel of the
monodromy representation. On the other hand, according to Proposition
\ref{prop:schw}, a), the monodromy representation \eqref{hge-mon} maps
the loops $\ell_1$ and $\ell_\infty$ to the standard generators of the modular group
$\operatorname{PSL}_2(\ZZ)$. The relations between the standard
generators are well known, namely the kernel of the monodromy
representation is the normal subgroup of
$\pi_1(\CC\setminus{\{0,1\}},x^\circ)$ generated by $\ell_1^2$ and
$(\ell_1\ell_\infty)^3$.

Let us work out the analytic transformation of $v^\infty(x)$ along the
loops $\ell_1$ and $\ell_\infty$. The analytic continuation along
$\ell_\infty$  acts on the vector column with entries $u^\infty$ and
$v^\infty$ as multiplication by the matrix
\ben
M_\infty = e^{2\pi\sqrt{-1}/12}\, \begin{bmatrix}
  1 & 2\pi\sqrt{-1} \\
  0 & 1 
\end{bmatrix}.
\een
In particular, $\tau(x)\mapsto \tau(x)+1$ and $v^\infty(x)\mapsto
e^{2\pi\sqrt{-1}/12}\, v^\infty(x).$ The analytic continuation along
$\ell_1$ requires a long but straightforward computation. Let us point
out the main steps leaving the details as an exercise. Using the
linear transformation formulas between the bases
$(u^\infty,v^\infty)$, $(u^0,v^0)$, and $(u^1,v^1)$ we find
\ben
\begin{bmatrix}
  e^{5\pi\sqrt{-1}/12} & A e^{5\pi\sqrt{-1}/12} \\
  e^{\pi\sqrt{-1}/12} & B e^{\pi\sqrt{-1}/12}
\end{bmatrix}\,
\begin{bmatrix} u^\infty \\ v^\infty \end{bmatrix} =
\begin{bmatrix}
  -2\alpha^{-1} & \alpha \\
  -2\beta^{-1}  & \beta
\end{bmatrix}\,
\begin{bmatrix} \Gamma(1/2) u^1 \\ \Gamma(1/2) v^1 \end{bmatrix},
\een
where
$A=-\pi\sqrt{-1} +\log(1728) -\pi\sqrt{3}$ and
$B=-\pi\sqrt{-1} +\log(1728) +\pi\sqrt{3}.$
From this formula, we find that the analytic continuation along
$\ell_1$ acts on the column vector with entries $u^\infty$ and
$v^\infty$ by multiplication by the matrix
\ben
M_1=
\frac{1}{2\pi}
\begin{bmatrix}
  -\log (1728) & -\log^2 (1728) + 4\pi^2\\
  1 & \log (1728) 
\end{bmatrix}\, .
\een
Note that 
\ben
\frac{1}{2\pi}\left(
  u^\infty(x) + \log(1728) v^\infty(x)\right) =
\frac{1}{2\pi}\left(
  \tau^\infty(x) + \log(1728) \right) v^\infty(x) = -\sqrt{-1} \,
\tau(x)\, v^\infty(x).
\een
Therefore, the analytic continuation along $\ell_1$ transforms
\ben
\tau(x)\mapsto -1/\tau(x),\quad
v^\infty(x)\mapsto -\sqrt{-1} \, \tau(x)\, v^\infty(x). 
\een
The above formula implies that $v^\infty(x)$ is invariant under the
analytic continuation along $\ell_1^2$. The analytic continuation
along the path $\ell_1\ell_\infty$ can be represented by the following
diagrams 
\begin{diagram}
  v^\infty(x) & \rMapsto^{\ell_1} & 
  -\sqrt{-1} \, \tau(x)\, v^\infty(x) & \rMapsto^{\ell_\infty} &
  e^{-\tfrac{\pi\sqrt{-1}}{3} }\, (\tau(x)+1)\, v^\infty(x)
\end{diagram}
and 
\begin{diagram}
  \tau(x) & \rMapsto^{\ell_1} & 
  -1/\tau(x)& \rMapsto^{\ell_\infty} &
  -1/(\tau(x)+1).
\end{diagram}
Therefore, the analytic continuation along $(\ell_1\ell_\infty)^3$
transforms $v^\infty(x)$ into $-v^\infty(x)$. We get that
$(v^\infty(x))^{2k}$ for all integers $k\geq 0$ is invariant under the
analytic  continuation along $\ell_1^2$ and $(\ell_1\ell_\infty)^3$.

Let us denote by $V(\tau)$ the pullback of $(v^\infty(x))^{12}$, i.e.,
$V$ is a holomorphic function on
$\HH\setminus{\{E_4(\tau)=0\}\cup\{E_6(\tau)=0\}}$ such that
$V(\tau(x))=(v^\infty(x))^{12}$, where the reference path defining the
branch of the Schwarz map $\tau(x)$ is chosen to be the same as the reference path
specifying the branch of  $v^\infty(x)$. The transformation formulas
of $\tau(x)$ and $v^\infty(x)$ under the analytic continuation along
$\ell_1$ and $\ell_\infty$ yield the following symmetries
\ben
V(\tau+1)=V(\tau),\quad V(-1/\tau) = \tau^{12} V(\tau). 
\een
This implies that $V$ transforms as a modular form of weight 12. Let
us check that $V$ is analytic at the cusp $\tau=\infty$. Suppose that
$\tau$ belongs to a Schwarz triangle $\Delta$ with vertex at
infinity. According to part a), the $J$-invariant defines a map that
is inverse to the analytic branch of the Schwarz map whose image is 
the Schwarz triangle $\Delta$. Then $V(\tau) = v^\infty(J(\tau))^{12}$, so
by substituting the Fourier series of the $J$-invariant
$x=J(\tau)=\tfrac{1}{1728} (q^{-1} +744+O(q))$ in   $v^\infty(x)^{12}$
we get that $V(\tau)$ has a Fourier expansion in terms of $q$ that has
only positive powers. The leading order term is $V(\tau) = 1728 q +O(q^2)$. 
Similar argument shows that $V(\tau)$ extends analytically across the
vertices of the Schwarz triangles. Therefore, $V(\tau)$ is a modular
form of weight 12. The space of modular forms of weight 12 has
dimension 2 and it is spanned by $E_4(\tau)^3$ and $E_6(\tau)^2$. We
get $V(\tau)=c_1 E_4(\tau)^3 + c_2 E_6(\tau)^2$. Comparing the
coefficients in the corresponding Fourier series in front of $q^0$ and $q^1$ we get
$c_1=1$ and $c_2=-1$. 
\qed

\subsection{Proof of Theorem \ref{t1}, a)}\label{sec:pd}
Let us recall formula \eqref{z-hge}. If $\lambda\neq 0$, then
$-Z_2(Q,\lambda)/(2Z_3(Q,\lambda))=\tau(x)$ is a point in the
uper-half plane satisfying the condition
$J(\tau(x))=x=\tfrac{\lambda^2}{27Q}\neq 1$, i.e., $Z_{\rm
  small}(Q,\lambda) \in \Omega_{\rm small}$. Since the period map is
locally analytic near $\lambda=0$, its value $Z_{\rm
  small}(Q,0)=\lim_{n\to \infty} Z_{\rm small}(Q,\lambda_n)$ can be
computed by choosing any convergent sequence $\lambda_n\to 0$. For
example, let us take $\lambda_n$ to be real. The Schwarz map
$\tau(\lambda_n)$ has a limit $\tau_0$ whose value depends on the reference
path, but in any case $\tau_0$ is a vertex of a Schwarz triangle and
$J(\tau_0)=\lim_n J(\tau(\lambda_n)) = \lim \lambda_n=0$, i.e.,
$\tau_0$ belongs to the orbit $\operatorname{PSL}_2(\ZZ)\cdot \rho$.
This completes the proof that $Z_{\rm small}: \widetilde{X}_{\rm
  small}\to \Omega_{\rm small}\setminus{\{E_6(-z_2/(2z_3))=0\}}$.
We have to check that this map is surjective.

Let us define the analytic isomorphism 
\beq\label{def:small-Phi}
\Phi_{\rm small}: \HH\times \CC^*\to \Omega_{\rm small},\quad (\tau,y)\mapsto (\tau^2\,y,
-2\tau\, y,y).
\eeq
Then we have to prove that the following map is surjective 
\ben
\widetilde{X}_{\rm small}\to (\HH\setminus{\{E_6(\tau)=0\}})\times \CC^*,\quad
(Q,\lambda)\mapsto (-Z_2(Q,\lambda)/(2Z_3(Q,\lambda)), Z_3(Q,\lambda)).
\een
Given a point $(\tau,y)\in \HH\times \CC^*$ we first pick
$(Q,\lambda)$ such that $\lambda^3/(27Q) = J(\tau)$ and then we fix
$Q$ in such a way that $Q^{1/6} z_3(\lambda^3/27Q) = y$. The
surjectivity follows.
\qed

\subsection{Proof of Theorem \ref{t1}, b)}
Recall that the period map $Z_{\rm small}(Q,\lambda)$ has the form
\eqref{z-hge}. According to Proposition \ref{hge-mon} if we put
$\tau=-\tfrac{Z_2(Q,\lambda)}{ 2 Z_3(Q,\lambda)}$, then we have  
\ben
Q Z_3(Q,\lambda)^6 = z_3(x)^6 =  \frac{8}{27}\, (2\pi)^{6}\,
(E_4(\tau)^3 -E_6(\tau)^2), 
\een
where $x:=\tfrac{\lambda^3}{27Q}$. Therefore,
\ben
Q = \frac{8}{27}\, (2\pi)^{6}\, Z_3(Q,\lambda)^{-6}\,  (E_4(\tau)^3 -E_6(\tau)^2).
\een
The formula for $\lambda$ follows from the relation 
\ben
\frac{\lambda^3}{27 Q} = x = J(\tau) = \frac{E_4(\tau)^3}{E_4(\tau)^3-E_6(\tau)^2},
\een
which implies that 
\ben
\lambda=2\, (2\pi)^{2}\, Z_3(Q,\lambda)^{-2}\, E_4(\tau)\, \xi,
\een
where $\xi^3=1$. 
On the other hand, $Z_3(Q,\lambda)^2 \lambda = 8\pi^2
+O(\lambda^{-1})$, so $\xi=1$.

\subsection{Proof of Theorem \ref{t1}, c)}
Let us first identify the period domain $\Omega_{\rm small}$ with
$\HH\times \CC^*$ via the isomorphism \eqref{def:small-Phi}. This
identification will induce an action of the 
monodromy group $W$ of quantum cohomology on $\HH\times \CC^*$. Let us
work out this action explicitly. Let us define a {\em left} action of
$\operatorname{PSL}_2(\ZZ)\times\{\pm 1\}$ on
$\HH\times \CC^*$ by 
\ben
(g,\sigma)\cdot(\tau,y) := \Big(\frac{a\tau+b}{c\tau+d}\, ,\,
\sigma \chi_2(g) (c\tau+d)^2 y\Big),\quad 
g=
\begin{bmatrix}
a& b\\
c& d
\end{bmatrix}\ \in\ \operatorname{PSL}_2(\ZZ).
\een
\begin{lemma}\label{le:mod-action}
Let $w\in W\subset \operatorname{GL}(\CC^3)$ be the monodromy transformation
corresponding to an element $(g,\sigma)\in
\operatorname{PSL}_2(\ZZ)\times\{\pm 1\}$  via the group homomorphism
\eqref{mod-mon}. Then 
\ben
\Phi_{\rm small}(\tau,y)\cdot w= \Phi_{\rm small}((s g^{-1} s^{-1},\sigma)\cdot (\tau, y)).
\een
where 
$
s=\begin{bmatrix}
0& 1\\
1& 0
\end{bmatrix}
$ and $(\tau,y)\in \HH\times \CC^*.$
\end{lemma}
\proof
Put $z:=\Phi_{\rm small}(\tau,y)$, then by definition
$\tau=\tau(z):=-z_2/(2z_3)$
and 
\ben
w=\sigma \chi_2(g) 
\begin{bmatrix}
a^2 & 2ab& b^2 \\
ac & (ad+bc) & bd \\
c^2 & 2cd & d^2  
\end{bmatrix}.
\een
After a straightforward computation we get
\ben
\tau(z\cdot w) = 
-\frac{2ab z_1 +(ad+bc)z_2+2cd z_3}{2(b^2z_1+bd z_2 +d^2 z_3)} = \frac{a\tau(z)-c}{-b\tau(z)+d},
\een
where we used that $z_1/z_3=z_1z_3/z_3^2 = \tau^2$. It remains only to
use that
\ben
\begin{bmatrix}
a& -c\\
-b& d
\end{bmatrix} = 
s\, \begin{bmatrix}
a& b\\
c& d
\end{bmatrix}^{-1}\, s^{-1}
\een
and $\chi_2(sg^{-1}s^{-1}) = \chi_2(g)$. \qed

The proof of part c) can be completed as follows.  Suppose that 
$\pi_{\rm small}(\tau,y) = \pi_{\rm small}(\tau',y')$. Since $J(\tau)=J(\tau')$, there exists
$g\in \operatorname{PSL}_2(\ZZ)$ such that $\tau'=g(\tau)$. There are
2 cases. First, if $J(\tau)\neq 0$, then $E_4(\tau)\neq 0$ and we get
\ben
y^{-2} E_4(\tau) = (y')^{-2} E_4(\tau')= (y')^{-2} (c\tau+d)^4 E_4(\tau),
\een 
so
\ben
y' = \sigma' (c\tau+d)^2 y
\een
for some $\sigma'\in \{\pm 1\}$. Defining $\sigma = \sigma'
\chi_2(g)$ we get $(\tau',y') = (g,\sigma) \cdot (\tau,y)$. 

The second case is the case when $J(\tau)=0$.  Such a $\tau$ is a
vertex in the triangulation by Schwarz triangles and it is in the same
$\operatorname{PSL}_2(\ZZ)$-orbit as
$\zeta=e^{\tfrac{2\pi\sqrt{-1}}{3}}$.  We may assume that
$\tau=\tau'=\zeta$, because the point $(\tau,y)$ (resp. $(\tau',y')$) is in the 
$W$-orbit of a point of the form $(\zeta,\widetilde{y})$
(resp. $(\zeta,\widetilde{y'})$). Using that
\ben
y^{-6} E_6(\tau)^2 = (y')^{-6} E_6(\tau')^2 = (y')^{-6} E_6(\tau)^2
\een
we get that $y'= \xi y$, for some $\xi$ such that
$\xi^6=1$. Note that the stabilizer of $\zeta\in \HH$ in the modular
group is a cyclic group of order 3 generated by $g_1\kappa$, where the
matrices $g_1,\kappa\in \operatorname{PSL}_2(\ZZ)$ are defined in Section 
\ref{sec:mon-2nd-con}.   We have
\ben
((g_1\kappa)^i,\sigma)\cdot (\tau,y) = (\tau,\sigma^i \zeta^{2i} y),
\een
where $i$ is an integer.  Since we can always choose $\sigma\in \{\pm
1\}$ and $i$ such that $\sigma^i \zeta^{2i}=\xi$ we get that
$(\tau,y)$ and $(\tau,y')$ are in the same $W$-orbit. 
\qed

\section{Holomorphic thickening}
Let us return to the general case of the period map for the big
quantum cohomology. Recall that 
\ben
X=(B\times \CC)' \subset \CC^*\times \CC\times \CC
\een
and
\ben
X_{\rm small}=\{t_3=0\}\subset X.
\een
Let us introduce coordinates $(Q,t,\lambda)$ on $X$ such that $Q:=e^{t_2}$,
$t:=t_3$. We fix a
base point $y^\circ=(Q^\circ,\lambda^\circ)$ in $X_{\rm small}$, which will
be used as a base point of $X$ as well. Given a point $(Q,\lambda)\in
X_{\rm small}$ we will be interested in the Taylor's series expansion
\ben
Z(Q,t,\lambda) = \sum_{n=0}^\infty Z^{(n)}(Q,\lambda)\frac{t^n}{n!}.
\een
We are going to construct a covering of $X$ such that the pullback of
$Z^{(n)}(Q,\lambda)$ is a quasi-modular form. This would allow us to find
the inverse of the period map and generalize to some extend the
statement of Theorem \ref{t1}. 

\subsection{Auxiliary covering}
Let $\U\subset \HH\times \CC^* \times \CC$ be an
open subset and let 
\ben
\pi^{\rm aux}: \mathcal{U} \to \CC^*\times \CC\times \CC,
\quad (\tau,x,s)\mapsto (Q,t,\lambda)
\een
be the holomorphic map defined by 
\ben
Q & = &\frac{8}{27}(2\pi/x)^6 (E_4(\tau)^3 - E_6(\tau)^2),\\
\lambda & = & 2(2\pi/x)^2 E_4(\tau), \\
t &= & s E_6(\tau)^2/x^6.
\een
\begin{remark}
So far the variable $x$ was used to denote the coordinate on
the domain $\CC\setminus{\{0,1\}}$ of the Schwarz map. We will no
longer deal with the Schwarz map. From now on we will use $x$ to
denote the coordinate on the complex circle $\CC^*\subset \U$,  except
for some local proofs where the demand for letters is hard to meet
without using $x$.
\end{remark}

We choose $\U$ to be the trivial disk bundle 
\ben
\U=\{(\tau,x,s)\ |\ |s|<\delta(\tau,x)\}
\een
where
\ben
\delta:\HH\times \CC^* \to \RR_{>0}
\een
is a smooth function defined as follows. We choose $\delta$ in such a way that the
preimage under $\pi^{\rm aux}$ of the discriminant is the analytic hypersurface
$E_6(\tau)=0$. More precisely, the
equation of the discriminant has the form
\ben
0=\operatorname{det}(\lambda-E\bullet) = \lambda^3-27Q + t\,g(t,Q,\lambda),
\een
where $g\in \O(B)$ is some holomorphic function. If
$(t,Q,\lambda)=\pi^{\rm aux}(\tau,x,s)$, then the above equation becomes
\ben
E_6(\tau)^2 \, \Big(8(2\pi)^6+ s g\circ \pi^{\rm aux}(\tau,x,s)\Big) \, x^{-6}=0.
\een
For fixed $(\tau,x)\in \HH\times \CC^*$ we choose $\delta(\tau,x)$
such that $\pi^{\rm aux}(\tau,x,s)\in B$ and  $|s\, g\circ\pi^{\rm aux}(\tau,x,s)|<8(2\pi)^6$ for all
$|s|<\delta(\tau,x)$. 

\begin{lemma}\label{X-retract-Y}
If we choose the constant $\epsilon$ in the definition of the domain
$B$ sufficiently small, then the subvariety $X_{\rm small}$ is a deformation retract of $X$.  
\end{lemma}
\proof
A deformation retraction 
\ben
\Psi: X\times [0,1]\to X
\een
can be taken in the form
\beq\label{retr-XtoY}
\Psi(Q,t,\lambda,s):= (Q,\lambda + \sum_{i=1}^3 \rho_i(\lambda
Q^{-1/3}) (u_i(Q,(1-s)t)-u_i(Q,t)), (1-s)t)
\eeq
where $u_i(Q,t)$ ($1\leq i\leq 3$) are the eigenvalues of the quantum
multiplication by $E\bullet_{Q,t}$, i.e., the {\em canonical
  coordinates}. Note that $u_i(Q,t) =Q^{1/3} u_i(1,t Q^{1/3})$ and
$u_i(Q,0)=3\zeta^iQ^{1/3}$, where $\zeta=e^{2\pi\sqrt{-1}/3}$. Given a real number $\delta>0$ we can always
choose $\epsilon$ sufficiently small so that
$|u_i(1,tQ^{1/3})-u_i(1,0)|<\delta$ for all $(Q,t)$, s.t.,
$|tQ^{1/3}|<\epsilon$.  We claim that if we choose $\delta<\frac{3}{16}|1-\zeta|$
and $\rho_i$ ($1\leq i\leq 3$) to be smooth functions such that 
$\rho_i(x)=1$ for all $|x-3\zeta^i| <4\delta$ and $\rho_i(x)=0$ for all
$|x-3\zeta^i|>8\delta$, then formula \eqref{retr-XtoY} defines a deformation
retract, i.e., a homotopy between the identity map and a retraction 
$X\to X_{\rm small}$.

Clearly we have $\Psi(x,0)=x$ for all $x\in X$ and $\Psi(y,s)=y$ for
all $y\in X_{\rm small}$ and for all $s$. We have to verify that
$\Psi(Q,t,\lambda,s)$ is not a point on the discriminant. There are
two cases. First, if $|\lambda Q^{-1/3}-3\zeta^i|>8\delta$ for all $i$,
then $\Psi(Q,t,\lambda,s)=(Q,\lambda,(1-s)t)$. We have 
\ben
&&
|\lambda Q^{-1/3}-u_j(1,(1-s)tQ^{1/3})|\geq \\
&&
|\lambda Q^{-1/3}-u_j(1,0)|-
|u_j(1,(1-s)tQ^{1/3})-u_j(1,0)|> 8\delta -\delta =7\delta>0,
\een
so $\lambda\neq u_j(Q,(1-s)t)$ for all $j$, i.e.,
$(Q,\lambda,(1-s)t)\in X$.  The second case is if $|\lambda
Q^{-1/3}-3\zeta^i|\leq 8\delta$ for some $i$. Note that if $j\neq i$,
then 
\ben
|\lambda Q^{-1/3}-3\zeta^j|\geq |3\zeta^i-3\zeta^j|-|\lambda
Q^{-1/3}-3\zeta^i|=3|1-\zeta|-|\lambda
Q^{-1/3}-3\zeta^i| > 8\delta.
\een
Therefore the second component of $\Psi(Q,t,\lambda,s)$ is 
\ben
\lambda + \rho_i(\lambda Q^{-1/3}) (u_i(Q,(1-s)t)-u_i(Q,t)).
\een
We have to prove that the above number does not coincide with $u_j(Q,(1-s)t)$ for
all $j$. Let us assume that this is not the case, i.e., the number
coincides with $u_j(Q,(1-s)t)$ for some $j$. Using the estimate
\ben
|\lambda Q^{-1/3} - 3 \zeta^j|\leq
|u_j(1,(1-s)t Q^{1/3})-3\zeta^j|+|u_i(1,(1-s)t Q^{1/3})-u_i(1,t Q^{1/3})|<3\delta
\een
we get that we must have $j=i$ and $\rho_i(\lambda
Q^{-1/3})=1$. Therefore our assumption implies that 
\ben
\lambda + u_i(Q,(1-s)t)-u_i(Q,t) = u_i(Q,(1-s)t)
\quad 
\Rightarrow
\quad 
\lambda = u_i(Q,t).
\een
This however contradicts the fact that $(Q,t,\lambda)\in X$. 
\qed

\begin{proposition}\label{prop:hol-lift}
a)
Let $\pi':\U'=(\pi^{\rm aux})^{-1}(X)\to X$ be the map induced from $\pi^{\rm aux}$. Then
the period map admits a holomorphic lift
$Z^{\rm aux}:\U'\to \CC^3$.

b) The map $Z^{\rm aux}$ extends holomorphically  on the entire domain $\U$. 
\end{proposition}
\proof
a)
Note that our definition of $\U$ implies that
$\U'=\U-\{E_6(\tau)=0\}$. 
Let us define an action of the monodromy group
$W=\operatorname{PSL}_2(\ZZ)\times \{\pm 1\}$ on $\U$ 
\ben
(g,\sigma)\cdot (\tau,x,s) = (g(\tau),\sigma\,\chi_2(g)\, x,s).
\een
Note that the points with non-trivial stabilizers are given by the
analytic hypersurfaces
$
\{E_4(\tau)=0\}
$
and $\{E_6(\tau)=0\}$. Let $u^\circ=(\tau^\circ,x^\circ,0)\in \U'$ be a reference point, such that
$E_4(\tau^\circ)\neq 0$ and $\pi^{\rm aux}(u^\circ)=y^\circ$.

Let us construct a lift $Z'$ of the period map on $\U'$. If $u=(\tau,x,t)\in \U'$, then
we pick a reference path $\gamma\subset \U'$ and define  
$
Z'(u)=Z(\pi^{\rm aux}(u))
$
where the value of $Z(\pi^{\rm aux}(u))$ is defined via the reference path
$\pi^{\rm aux}(\gamma)$. We claim that choosing a different reference path $\gamma'\subset
\U'$ does not change the value of $Z'(u)$. In other words we claim
that if $L\in \pi_1( \U',u^\circ)$ is a loop based at $u^\circ$, then the image
$\pi^{\rm aux}(L)$ is in the kernel of the monodromy representation ${\rm
  per}:\pi_1(X)\to W$. By making a small perturbation
(without changing the homotopy class) we can arrange that $L$ is a
loop in $\U'-\{E_4=0\}$.  Note that the
projections $r':(\tau,x,s)\mapsto (\tau,x)$ and $r'':(Q,t,\lambda)\mapsto
(Q,\lambda)$ give rise to a commutative diagram 
\ben
\begin{CD}
\U' @>r'>> (\HH\times \CC^*)'& := &  \HH\times \CC^* -\{E_6=0\} \\
@V\pi' VV @VV\pi'' V & & \\
X' @>r''>> X_{\rm small}&  &  
\end{CD}
\een
where $\pi'$ and $\pi''$ are the maps induced from $\pi^{\rm aux}$, 
\ben
X' := \{(Q,t,\lambda)\in X\,:\, \lambda^3-27Q\neq 0\},
\een
and the horizontal arrows are deformation retractions. The map $\pi''$
induces a covering
\beq\label{covering}
\HH\times \CC^* -\{E_6=0\} \cup\{E_4=0\} \to X_{\rm small}':=X_{\rm small}-\{\lambda=0\}.
\eeq
According to Theorem \ref{t1}, locally the inverse of the covering map
\eqref{covering} is given explicitly by the following formulas
\ben
\tau=-Z_2(Q,\lambda)/(2Z_3(Q,\lambda)),\quad x=Z_3(Q,\lambda).
\een
Therefore  we have  a commutative diagram 
\ben
\begin{diagram}
\pi_1(X_{\rm small}',y^\circ)& & \rTo & &\pi_1(X_{\rm small},y^\circ)  \\
&\rdTo_{\rm cov} &   & \ldTo_{\rm per}& \\
& & W & &
\end{diagram}
\een
in which the horizontal arrow is induced from the natural inclusion
$X_{\rm small}'\subset X_{\rm small}$ and the two diagonal arrows are given by the monodromy
representations respectively of the covering
and the period maps. On the other hand the lift of the loop
$r''\circ\pi'(L)$ is $r'(L) $,
which is a loop, so the corresponding monodromy transformation  
$$
w:={\rm cov}(r''\circ\pi'(L))\in W
$$ 
fixes the reference point $(\tau^\circ,x^\circ) \in \U'$. Therefore
$w=1$, because the stabilizer of $(\tau^\circ,x^\circ)  $ is trivial. We get that the homotopy class of 
$r''\circ\pi'(L)$ in $\pi_1(X_{\rm small},y^\circ)$ is in the
kernel of the monodromy representation of the period map. Using that
$r''$ is a deformation retract, we get that
$\pi'(L)$ is homotopic to 
$r''\circ\pi'(L)$ in $X'$. Finally, since
$X_{\rm small}$ is a deformation retract of $X$ (see Lemma \ref{X-retract-Y}) we
get that the homotopy class of $\pi'(L) $ in 
$
\pi_1(X,y^\circ)
$
must be in the kernel of the monodromy representation of the period
map.

b)
It remains only to prove that $Z'$ extends analytically to the entire
domain $\U$. The complement of $\U'$ in $\U$ is an analytic
hypersurface. Recalling the Riemann extension theorem, we get that it
is sufficient to prove that the values of $Z'(u)$ are bounded in a
neighborhood of an arbitrary point $u_0\in \U-\U'$. Note that 
$\pi^{\rm aux}(u_0)=:(Q_0,t_0,\lambda_0)$ is a point on the 
discriminant. Then by definition the periods $I^{(-1)}(Q,t,\lambda)\sim
(\lambda-u)^{1/2}$ where $\lambda=u$ is the local equation of the
discriminant near the point $(Q_0,t_0,\lambda_0)$. Therefore, the map
$Z'$ is bounded. 
\qed

\subsection{The Taylor's coefficients $Z^{(n)}$}\label{sec:taylor-coef}

Recall the notation in the proof of Lemma \ref{le:quadr-rel}. Let us
denote by $I^{(n)}$ the matrix whose $(i,j)$ entry is
$(I^{(n)}_{E_i},p^{j-1})$. We claim that the matrix $I^{(-1)}$ can be expressed in terms of the
Wronskian matrix
\ben
\operatorname{Wr} = 
\begin{bmatrix}
Z_1 & Q\partial_Q Z_1 & (Q\partial_Q)^2 Z_1 \\
Z_2 & Q\partial_Q Z_2 & (Q\partial_Q)^2 Z_2 \\
Z_3 & Q\partial_Q Z_3 & (Q\partial_Q)^2 Z_3 
\end{bmatrix}.
\een
Indeed, put
$$
A=A(Q,t,\lambda) = -(-\theta+1/2)\, (\lambda-E\bullet)^{-1}.
$$
The differential equation of the second structure connection can be
written us
\ben
\partial_\lambda I^{(-1)} = -I^{(-1)}A,\quad
Q\partial_Q  I^{(-1)} = I^{(-1)}A\Omega_2,\quad 
\partial_t  I^{(-1)} = I^{(-1)}A \Omega_3,
\een
where $\Omega_i= P^{i-1}\bullet$ is the matrix of quantum
multiplication by $P^{i-1}$. We get  
\ben
Q\partial_Q
\begin{bmatrix}
Z_1\\ Z_2 \\ Z_3
\end{bmatrix} = I^{(-1)} A\, \Omega_2\, 
\begin{bmatrix}
1\\ 0 \\ 0
\end{bmatrix},\quad 
(Q\partial_Q)^2
\begin{bmatrix}
Z_1\\ Z_2 \\ Z_3
\end{bmatrix} = I^{(-1)} (A\, \Omega_2\, A +Q\partial_Q A )\,\Omega_2 
\begin{bmatrix}
1\\ 0 \\ 0
\end{bmatrix}.
\een
Therefore, $\operatorname{Wr}=I^{(-1)} T$, where $T$ is the matrix
with columns 
\ben
\begin{bmatrix}
1\\ 0 \\ 0
\end{bmatrix},\quad 
A\, \Omega_2\, 
\begin{bmatrix}
1\\ 0 \\ 0
\end{bmatrix},\quad 
 (A\, \Omega_2\, A +Q\partial_Q A )\,\Omega_2 
\begin{bmatrix}
1\\ 0 \\ 0
\end{bmatrix}.
\een
The proves of the next two Lemmas involve some long
computations. Although they could be done by hand within acceptable
amount of time, we recommend the use of a computer software such as Mathematica or
Maple. 
\begin{lemma}
The matrix $T$ can be expressed in terms of the genus 0 potential $F$ as
follows
\ben
T=\begin{bmatrix}
1 & \frac{1}{\Delta} (F_{23}\lambda +9F_{33}/2) & \frac{1}{\Delta^2}
t_{13} \\
 0 &  \frac{1}{\Delta} (-\lambda^2/2-3F_{33}t/2) & \frac{1}{\Delta^2}
t_{23} \\
0 &  \frac{1}{\Delta} (-9\lambda/2+3F_{23} t) & \frac{1}{\Delta^2}
t_{33}
\end{bmatrix},
\een
where $\Delta=\operatorname{det}(\lambda-E\bullet)$, 
\ben
t_{13} & = & 
3F_{223}\lambda^4/4+(-3 F_{22} F_{223} + 2 F_{222} F_{23} + 9 F_{233}
- 3 F_{33})\lambda^3/4 \\
&&
+ 
(-12 F_{223} F_{23} - 9 F_{22} F_{233} + 3 F_{22} F_{33} + 9 F_{222} F_{33} - 6 F_{23} F_{233} t + 
  9 F_{223} F_{33} t)\lambda^2/4 \\
&&
+
(-54 F_{23} F_{233} + 27 F_{223} F_{33} - 4 F_{223} F_{23}^2 t + 6 F_{22} F_{23} F_{233} t - 
  9 F_{22} F_{223} F_{33} t \\
&&
- 9 F_{33}^2 t + 6 F_{23} F_{33} (6 + F_{222} t)) \lambda/4+
(81 F_{33}^2 + 36 F_{23}^2 F_{233} t - 72 F_{223} F_{23} F_{33} t \\
&&
- 12 F_{23}^2 F_{33} t + 
  9 F_{22} F_{33}^2 t + 27 F_{222} F_{33}^2 t)/4,
\een
\ben
t_{23}& = &  -F_{222}\lambda^4/4-3F_{223}\lambda^3 + 
(-27F_{233}/4+2F_{223}F_{23}t-3F_{222}F_{33} t/2) \lambda^2+ \\
&&
+(9F_{23}F_{233}t-9F_{223}F_{33}t)\lambda 
+(-3F_{23}^2F_{233}+6F_{223}F_{23}F_{33}-9F_{222}F_{33}^2/4)t^2,
\een
and
\ben
t_{33}
&=&
3\lambda^4/4 -3\lambda^3(F_{22} + 3 F_{222} - 3 F_{223} t)/4 -3 (36
F_{223} + 12 F_{23} + 3 F_{22} F_{223} t \\
&&
- 2 F_{222} F_{23} t - 9 F_{233} t -    3 F_{33} t) \lambda^2/4 -
3 (81 F_{233} + 27 F_{33} - 12 F_{223} F_{23} t \\
&& 
- 4 F_{23}^2 t + 9 F_{22} F_{233} t + 
   3 F_{22} F_{33} t + 9 F_{222} F_{33} t + 6 F_{23} F_{233} t^2 - 
   9 F_{223} F_{33} t^2) \lambda/4 \\
&&
-3 (-54 F_{23} F_{233} t + 81 F_{223} F_{33} t + 4 F_{223} F_{23}^2 t^2 - 
   6 F_{22} F_{23} F_{233} t^2 + 9 F_{22} F_{223} F_{33} t^2 \\
&& 
- 6 F_{222} F_{23} F_{33} t^2)/4.
\een
The homogeneous degree with respect to $Q,t,\lambda$ of the $i$th row of
$T$ is $1-i$. 
\end{lemma}
\proof
Using that 
\ben
\Omega_2=
\begin{bmatrix}
0 & F_{223} & F_{233} \\
1 & F_{222} & F_{223} \\
0 & 1 & 0
\end{bmatrix},\quad
\Omega_3=
\begin{bmatrix}
0 & F_{233} & F_{333} \\
0& F_{223} & F_{233} \\
1 & 0 & 0
\end{bmatrix}
\een
and that $F$ is homogeneous of degree 1 we get 
\ben
\lambda-E\bullet = 
\Omega_2=
\begin{bmatrix}
\lambda & -2F_{23} & -3F_{33} \\
-3 &\lambda - F_{22} & -2F_{23} \\
t & -3 & \lambda
\end{bmatrix}.
\een
We can express $A$ in terms of the partial derivatives of $F$ and
after some long but straightforward computation we get the formulas
stated in the Lemma.
\qed

Another long but straightforward computation yields that 
\ben
\operatorname{det}(T) = -\frac{3}{8\Delta^2}\, \Delta^{(1)} ,
\een
where 
\ben
\Delta^{(1)} = 
\lambda^3 +  3 F_{223} t \lambda^2 
+ 3 (3 F_{233} + F_{33}) t \lambda 
-3 (2 F_{23} F_{233} - 3 F_{223} F_{33}) t^2 .
\een
\begin{lemma}\label{le:operator-L}
a)
There exists an operator
\ben
L(Q,t,\lambda;\partial_Q)=
L_0(Q,t,\lambda) +L_1(Q,t,\lambda) \, (Q\partial_Q)+L_2(Q,t,\lambda) \, (Q\partial_Q)^2 
\een
whose coefficients are rational functions in $\lambda$ depending
analytically on $(Q,t)\in B$ such that
\ben
\partial_t Z(Q,t,\lambda) = L(Q,t,\lambda;\partial_Q)\, Z(Q,t,\lambda).
\een 
b) The coefficients $L_i(Q,t,\lambda)$ have the form 
\ben
L_i(Q,t,\lambda) = \frac{1}{\Delta^{(1)}} \, \ell_i(Q,t,\lambda),\quad
0\leq i\leq 2,
\een
where $\ell_i$ is a polynomial in $\lambda$ of degree $2+i$ whose
coefficients are polynomials in the partial derivatives of
$F$. Moreover, the weigh of $\ell_i$ with respect to the variables
$Q,t,\lambda$ is 4.
\end{lemma}
\proof
Let $Z$ be the column with entries $Z_1,Z_2,Z_3$ and
$\{e_i\}_{i=1}^3\subset \CC^3$ be the standard basis. We have 
\ben
\partial_t Z = \partial_t  I^{(-1)}e_1 = I^{(-1)} A \Omega_3 e_1 =
\operatorname{Wr} T^{-1} A e_3.
\een
The 3rd column $Ae_3$ of the matrix $A$ can be expressed in terms of
the partial derivatives of $F$ as explained above. 
Therefore the coefficients of the differential operator are given by
\ben
\begin{bmatrix} L_0 \\ L_1\\ L_2\end{bmatrix} =\frac{1}{\Delta}\, 
T^{-1} 
\begin{bmatrix} 
3F_{33}\lambda/2 +(4F_{23}^2-3F_{22}F_{33})  \\ 
-F_{23}\lambda -9F_{33}/2 \\ 
-3\lambda^2/2+3F_{22}\lambda/2+9F_{23} \end{bmatrix} .
\een
The rest of the proof is a straightforward computation.
\qed

Let us point out that at $t=0$ we have 
\ben
A=
\frac{1}{\lambda^3-27Q}
\begin{bmatrix}
\lambda^2/2 & 9Q/2 & 3\lambda Q/2 \\
-3\lambda/2 & -\lambda^2/2 & -9Q/2\\
-27/2 & -9\lambda/2 & -3\lambda^2/2
\end{bmatrix},\quad
\Omega_2=
\begin{bmatrix}
0 & 0 & Q \\
1 & 0 & 0 \\
0 & 1 & 0
\end{bmatrix},\quad
\Omega_3=\Omega_2^2,
\een
and 
\beq\label{T:t=0}
T=\frac{1}{\lambda^3-27Q}\, 
\begin{bmatrix}
1 & 9Q/2 & \frac{3Q(2\lambda^3+27Q)}{4(\lambda^3-27Q)} \\
0 & -\lambda^2/2 & -\frac{27Q\lambda^2}{4(\lambda^3-27Q)} \\
0 & -9\lambda/2 & \frac{3\lambda(\lambda^3-108Q)}{4(\lambda^3-27Q)}
\end{bmatrix}.
\eeq
The differential operator takes the form
\ben
L(Q,0,\lambda,\partial_Q)= \lambda^{-2}
\Big( \frac{9Q}{2} +36Q (Q\partial_Q)+2(27Q-\lambda^3)(Q\partial_Q)^2\Big).
\een
\begin{lemma}\label{le:Q-hge}
At $t=0$ the period map satisfies the following differential equation
\ben
(Q\partial_Q)^3 Z(Q,\lambda) = \frac{3Q}{8(\lambda^3-27Q)}\Big( 
5 +46 (Q\partial_Q) + 108 (Q\partial_Q)^2\Big) \, Z(Q,\lambda).
\een
\end{lemma}
\proof
We just need to check that the substitution $Z_i(Q,\lambda) = Q^{-1/6} z_i(x),$ with
$x=\lambda^3/(27Q)$ transforms the above differential equation into
the generalized hypergeometric equation \eqref{hge:3F2}. This however
is a straightforward computation.
\qed

We would like to change the coordinates $(Q,\lambda)\in \CC^*\times
\CC$ using the covering map $\pi$ in Theorem \ref{t1}, i.e.,
\beq\label{pi-coord}
\lambda=2(2\pi/x)^2 E_4(\tau),\quad
Q= \frac{8}{27}(2\pi/x)^6 (E_4(\tau)^3-E_6(\tau)^2). 
\eeq
Recall the Ramanujan's differential equations for the Eisenstein series
\ben
q\partial_q E_2 & = &
\frac{1}{12}(E_2^2-E_4),\\
q\partial_q E_4 & = &
\frac{1}{3}(E_2E_4-E_6),\\
q\partial_q E_6 & = &
\frac{1}{2}(E_2E_6-E_4^2).
\een
\begin{lemma}\label{le:log-Q-der}
Under the change of coordinates \eqref{pi-coord} we have 
\ben
Q\partial_Q =\frac{E_4}{E_6} q\partial_q +\frac{1}{6E_6} (E_2E_4-E_6)
x\partial_x .
\een
\end{lemma}
\proof
After a short computation we get
\ben
\partial_x \lambda & = & -4 (2\pi)^2 x^{-3} E_4(\tau), \\
\partial_xQ & = &  -\frac{16}{9} (2\pi)^6 x^{-7}  (E_4(\tau)^3-E_6(\tau)^2)\\
q\partial_q \lambda  & = &  \frac{2}{3} (2\pi/x)^2
(E_2E_4-E_6),\\
q\partial_q Q & = &  \frac{8}{27} (2\pi/x)^6E_2(E_4^3-E_6^2).
\een
The formula for $Q\partial_Q$ is easy to derive from here.
\qed

\begin{proposition}\label{prop:taylor-coef}
Under the change of coordinates \eqref{pi-coord} we  have 

a) The Taylor coefficient $Z_3^{(0)}(Q,\lambda) = x$ and 
\ben
Z_3^{(n)}(Q,\lambda)\quad \in\quad  x^{1-2n}\,E_6^{-2n}\,  \CC[E_2,E_4,E_6],\quad n>0.
\een

b) We have $Z_2^{(0)}(Q,\lambda) = -2\tau x$ and 
\ben
Z_2^{(n)}+2\tau Z_3^{(n)}(Q,\lambda)\quad \in\quad  
x^{1-2n}\,E_6^{-2n}\,  \CC[E_2,E_4,E_6],\quad n>0.
\een
\end{proposition}
\proof
Using Lemma \ref{le:operator-L} and \ref{le:Q-hge} we get that 
\ben
Z^{(n)}(Q,\lambda) = M^{(n)}(Q,\lambda;\partial_Q)\, Z(Q,\lambda),
\een
where $M^{(n)}$ is a second order differential operator of the form
\ben
M^{(n)}(Q,\lambda;\partial_Q) = \sum_{a=0}^2 M_a^{(n)}(Q,\lambda) \, (Q\partial_Q)^a
\een
whose coefficients are rational functions in $Q$ and $\lambda$ with
poles only at $\lambda=0$ and $\lambda^3-27Q=0$. According to Theorem
\ref{t1}, under the change of coordinates \eqref{pi-coord} we have
\ben
Z^{(0)}(Q,\lambda)  =Z(Q,\lambda) = (\tau^2 x, -2\tau x, x).
\een
Using Lemma \ref{le:log-Q-der} we also have 
\ben
M^{(n)} & = & M^{(n)}_0+M^{(n)}_1 b +M^{(n)}_2(b^2+a(q\partial_q b)) +
\\
&& 
\Big( M^{(n)}_1 a + M^{(n)}_2(a(q\partial_qa) + 2ab) \Big) q\partial_q
+
M^{(n)}_2 a^2 (q\partial_q)^2,
\een
where 
\ben
a:=E_4/E_6,\quad b:=\frac{1}{6}(E_2E_4/E_6-1).
\een
Note that in the above formula for $M^{(n)}$ we have replaced
$x\partial_x$ with 1, because $x\partial_x$ commutes 
with $L$ and it acts on $Z(Q,\lambda)$ by multiplication by 1.
Hence
\ben
Z^{(n)}_3 = (M^{(n)}_0+M^{(n)}_1 b +M^{(n)}_2(b^2+a(q\partial_q b))) x
\een
and 
\ben
Z^{(n)}_2 = -2Z^{(n)}_3 \tau + 
\frac{1}{2\pi\sqrt{-1}} 
\Big( M^{(n)}_1 a + M^{(n)}_2(a(q\partial_qa)+ 2ab) \Big) x
\een
The statements of part a) and b), modulo the order of the poles at $E_4=0$ and
$E_6=0$,  follows from the fact that $Z^{(n)}$ is
homogeneous of degree $n-\frac{1}{2}$, $\tau$ has degree 0 and $x$ has degree
$-\frac{1}{2}$. The statement that $Z^{(n)}_3$ and $Z^{(n)}_2+2\tau Z^{(n)}_3$ do not have a pole at
$E_4=0$ and have a pole of order at most $2n$ at $E_6=0$ follows
from Proposition \ref{prop:hol-lift}. Indeed, according to the Proposition  
the series  
\ben
\sum_{n=0}^\infty Z^{(n)}(\tau,x)\frac{(E_6(\tau)^2 x^{-6} s)^n}{n!}
\een
is convergent for all $(\tau,x,s)\in \U$, so in particular the
coefficient in front of $s^n$ must be holomorphic for all $(\tau,x)\in
\HH\times \CC^*$ and for all $n\geq 0$.  
\qed

The first component $Z_1(Q,t,\lambda)$ of the period map is determined
from the remaining two via the following relation.
\begin{lemma}\label{le:big-q-rel}
We have 
\ben
Z_2^2-4Z_1Z_3 = -32 t.
\een
\end{lemma} 
\proof
The argument is the same as in the proof of Lemma
\ref{le:quadr-rel}. Namely, we have 
\ben
\sum_{i,j=1}^3 \eta^{ij} Z_i Z_j =- \frac{1}{8}\left( Z_2^2-4Z_1Z_3\right)
\een
and the same argument as in Lemma \ref{le:quadr-rel} proves that the
LHS is the $(1,1)$-entry of the matrix
\beq\label{11-entry}
(-\theta+1/2)^{-1} g (\lambda-E\bullet) (-\theta+1/2)^{-1}.
\eeq
The entries of
$\lambda-E\bullet$ can be expressed in terms of the partial derivatives
of the genus zero potential $F$
\ben
\lambda-E\bullet = 
\begin{bmatrix}
\lambda & -2F_{23} & -3F_{33} \\
-3 & \lambda -F_{22}& -2F_{23} \\
t & -3 & \lambda
\end{bmatrix}.
\een 
Note that the $(1,1)$-entry of $g(\lambda-E\bullet)$ is $t$, so the
$(1,1)$-entry of \eqref{11-entry} is $4t$. 
\qed

\subsection{Extension of the period domain}

Recall that we have identified $\CC^3$ with the
space of quadratic forms in two variables (see \eqref{def:phi}). Let
us define an open neighborhood of $\Omega_{\rm small}$ in $\CC^3$ as
the image of the following map:
\beq\label{def:big-Phi}
\Phi: \HH^2\times \CC^* \to \CC^3, (\tau_1,\tau_2,y)\mapsto (z_1,z_2,z_3):=\phi^{-1}(y(v-u\tau_1)(v-u\tau_2)).
\eeq
Recalling the definition of $\phi$ we get that 
\ben
z_1=\tau_1\tau_2\, y,\quad
z_2=-(\tau_1+\tau_2)\, y,\quad z_3=y.
\een

Let us equip $\HH^2\times \CC^*$ with a left $W$-action. If
$w=(g,\sigma)\in W=\operatorname{PSL}_2(\ZZ)\times \{\pm 1\}$, then
we define 
\ben
w\cdot(\tau_1,\tau_2,y) := 
\left( 
\frac{a\tau_1+b}{c\tau_1+d},
\frac{a\tau_2+b}{c\tau_1+d}, \sigma \chi_2(g)
(c\tau_1+d)(c\tau_2+d) y\right).
\een 
\begin{lemma}\label{le:big-mod-action}
Let $w\in \operatorname{GL}(\CC^3)$ be the monodromy transformation
corresponding to an element $(g,\sigma)\in
\operatorname{PSL}_2(\ZZ)\times\{\pm 1\}$  via the map
\eqref{mod-mon}. Then 
\ben
\Phi(\tau_1,\tau_2,y)\cdot w= \Phi((s g^{-1} s^{-1},\sigma)\cdot (\tau_1,\tau_2, y)).
\een
where 
$
s=\begin{bmatrix}
0& 1\\
1& 0
\end{bmatrix}
$ and $(\tau_1,\tau_2,y)\in \HH^2\times \CC^*.$
\end{lemma}
\proof
We prove that the quadratic forms corresponding to the LHS and the RHS
(of the identity that we have to prove) coincide. 
The quadratic form corresponding to the LHS is 
\beq\label{qf-lhs}
\sigma \chi_2(g) y (cu+dv-(au+bv)\tau_1) (cu+dv-(au+bv)\tau_2),
\eeq
where $g=\begin{bmatrix} a & b\\c & d\end{bmatrix}$. Note that 
\ben
\begin{bmatrix}
a& -c\\
-b& d
\end{bmatrix} = 
s\, \begin{bmatrix}
a& b\\
c& d
\end{bmatrix}^{-1}\, s^{-1}
\een
and $\chi_2(sg^{-1}s^{-1}) = \chi_2(g)$. Therefore, the quadratic form
corresponding to the RHS is
\ben
\sigma\chi_2(g) (-b\tau_1+d)(-b\tau_2+d)y
\Big( v-u\frac{a\tau_1-c}{-b\tau_1+d}\Big)
\Big( v-u\frac{a\tau_2-c}{-b\tau_2+d}\Big).
\een
It remains only to verify that the above formula coincides with
\eqref{qf-lhs}. 
\qed

\subsection{Proof of Theorem \ref{t2}}
We would like to invert the period map
\ben
z_i=Z_i(Q,t,\lambda),\quad 1\leq i\leq 3,
\een
i.e., express $(Q,t,\lambda)$ in terms of $(z_1,z_2,z_3)$. 
Recall that 
\ben
(Q,t,\lambda)=\pi^{\rm aux}(\tau,x,t),\quad (z_1,z_2,z_3) = \Phi(\tau_1,\tau_2,y).
\een  
Since the inverse of $\Phi$ is straightforward to find, it is
sufficient to find the relation between the coordinate systems $(\tau,x,t)$ and
$(\tau_1,\tau_2,y)$.

To begin with note that according to Lemma \ref{le:big-q-rel} we have
\ben
t=-\frac{1}{32} (\tau_1-\tau_2)^2 y^2.
\een 
According to Proposition \ref{prop:taylor-coef} we have 
\ben
y=z_3 = x\Big(1+\sum_{n=1}^\infty y_{n}(E_2,E_4,E_6) (tx^{-2})^n\Big)  
\een
and 
\ben
\tau_{12}:=\frac{1}{2}(\tau_1+\tau_2) = -\frac{z_2}{2z_3} = 
\tau + \sum_{n=1}^\infty \widetilde{\tau}_{12,n}(E_2,E_4,E_6) (tx^{-2})^n,
\een
where $y_{n},\widetilde{\tau}_{12,n}\in \CC[E_2,E_4,E_6]E_6^{-2n}$. 
Using the formula for $y$ we can express $x$ in terms of $y,t,$ and
$E_i$ ($i=2,4,6$):
\ben
x=y\Big(1+\sum_{n=1}^\infty x_{n}(E_2,E_4,E_6) (ty^{-2})^n\Big) ,
\een
where $x_n\in \CC[E_2,E_4,E_6]E_6^{-2n}$.
Substituting this into the formula for $\tau_{12}$ we get 
\ben
\tau_{12}= 
\tau + \sum_{n=1}^\infty \tau_{12,n}(E_2,E_4,E_6) (ty^{-2})^n,
\een
where $\tau_{12,n}\in \CC[E_2,E_4,E_6]E_6^{-2n}.$ Using Taylor series expansion
at $\tau_{12}=\tau$ and the Ramanujan's differential equations. We get 
\ben
E_i(\tau_{12}) = E_i(\tau)+ \sum_{n=1}^\infty 
\tau^{(i)}_{12,n}(E_2,E_4,E_6) (ty^{-2})^n,\quad i=2,4,6.
\een
Therefore, we can express $E_i(\tau)$ ($i=2,4,6$) in terms of
$E_i(\tau_{12})$ ($i=2,4,6$)
\ben
E_i(\tau) = E_i(\tau_{12})+ 
\sum_{n=1}^\infty 
\tau^{(i)}_{n}(E_2(\tau_{12}),E_4(\tau_{12}),E_6(\tau_{12}))
(ty^{-2})^n .
\een
Since $ty^{-2}=-(\tau_1-\tau_2)^2/32$ we get inversion formulas of the
following type
\ben
Q & = & \frac{8}{27}
(2\pi/y)^6\, \sum_{n=0}^\infty Q_n(\tau_{12}) (\tau_1-\tau_2)^{2n} ,\\
\lambda & = & 2
(2\pi/y)^2\, \sum_{n=0}^\infty \lambda_n(\tau_{12}) (\tau_1-\tau_2)^{2n}
,\\
t & = & -\frac{1}{32}(\tau_1-\tau_2)^2 y^2,
\een
where $Q_n$ and $\lambda_n$ are polynomial expression in
$E_2(\tau_{12}), E_4(\tau_{12})$, and $E_6(\tau_{12})^{\pm 1}$. 
We have to prove that $Q_n$ and $\lambda_n$ depend
polynomially on $E_i$, i.e., there is no negative powers of
$E_6$. This follows from the fact that the period map for the
second structure connection is locally invertible. 

\begin{lemma}\label{le:Jac-det}
a) The value of the Jacobian determinant 
\ben
\frac{D(Z_1,Z_2,Z_3)}{D(\lambda,Q,t)}=
\begin{bmatrix}
\partial_\lambda Z_1 & \partial_\lambda Z_2 & \partial_\lambda Z_3\\
Q\partial_Q Z_1 & Q\partial_Q Z_2 & Q\partial_Q Z_3\\ 
\partial_t Z_1 & \partial_t Z_2 & \partial_t Z_3
\end{bmatrix} .
\een
at $t=0$ up to a non-zero constant coincides with
$(\lambda^3-27Q)^{-1/2}$. 

b) The value of the Jacobian determinant
\ben
\frac{D(\lambda,Q,t)}{D(\tau,x,s)} =
\begin{bmatrix}
\partial_\tau \lambda & \partial_\tau Q & \partial_\tau t \\
\partial_x \lambda & \partial_x Q & \partial_x t \\
\partial_s \lambda & \partial_s Q & \partial_s t
\end{bmatrix}
\een
at $s=0$ up to a non-zero constant coincides with $E_6(\tau)^2x^{-6}
Q(\lambda^3-27Q)^{1/2}$.  
\end{lemma}
\proof
We will use the notation from section \ref{sec:taylor-coef}.
Using the differential equations of the second structure connection we get 
\ben
\partial_\lambda Z_i = (I^{(0)}_{E_i},1),\quad 
Q\partial_Q Z_i = -(I^{(0)}_{E_i},P),\quad 
\partial_t Z_i = -(I^{(0)}_{E_i} ,P^2).
\een 
Therefore
\ben
\operatorname{det}\, \frac{D(Z_1,Z_2,Z_3)}{D(\lambda,Q,t)} = 
\operatorname{det}\, I^{(0)}= 
\operatorname{det}\, (I^{(-1)} (-\theta+1/2)\, (\lambda-E\bullet)^{-1}).
\een
The above expression should be evaluated at $t=0$. We get
\ben
-\frac{3}{8(\lambda^3-27Q)}\, \operatorname{det}\, I^{(-1)} .
\een
By definition $\operatorname{Wr}=I^{(-1)}T$ and the $T(Q,0,\lambda)$
is given by \eqref{T:t=0}. Therefore,
\ben
\operatorname{det}\, I^{(-1)} = -\frac{8}{3} \lambda^{-3}
(\lambda^3-27 Q)^2\operatorname{det}\,\operatorname{Wr}.
\een
The Jacobian determinant takes the form
\ben
\operatorname{det}\, \frac{D(Z_1,Z_2,Z_3)}{D(\lambda,Q,t)} =
\lambda^{-3}
(\lambda^3-27 Q)\operatorname{det}\,\operatorname{Wr}.
\een
Recall that $Z_i(Q,\lambda) = Q^{-1/6} z_i(x),$ with
$x=\lambda^3/(27Q)$. Therefore, the Wronskian determinant takes the
form
\ben
-Q^{-1/2}\operatorname{det}\,
\begin{bmatrix}
z_1 & (x\partial_x +1/6) z_1 & (x\partial_x +1/6)^2 z_1 \\
z_2 & (x\partial_x +1/6) z_2 & (x\partial_x +1/6)^2 z_2 \\
z_3 & (x\partial_x +1/6) z_3 & (x\partial_x +1/6)^2 z_3 
\end{bmatrix} ,
\een
i.e.,
\ben
\operatorname{det}\,\operatorname{Wr} = 
-Q^{-1/2}\operatorname{det}\,
\begin{bmatrix}
z_1 & x\partial_x z_1 & (x\partial_x )^2 z_1 \\
z_2 & x\partial_x z_2 & (x\partial_x )^2 z_2 \\
z_3 & x\partial_x z_3 & (x\partial_x )^2 z_3 
\end{bmatrix}.
\een
On the other hand, $\{z_i\}_{i=1}^3$ form a basis of solutions of the hypergeometric equation
\eqref{hge:3F2}. Therefore, the above determinant can be expressed
easily in terms of the Wronskian of the differential equation. After a
short computation we get
\ben
\operatorname{det}\,
\begin{bmatrix}
z_1 & x\partial_x z_1 & (x\partial_x )^2 z_1 \\
z_2 & x\partial_x z_2 & (x\partial_x )^2 z_2 \\
z_3 & x\partial_x z_3 & (x\partial_x )^2 z_3 
\end{bmatrix} = 
C\, x (1-x)^{-3/2} = 3\sqrt{3} C\sqrt{-1} 
\frac{Q^{1/2}\lambda^3}{(\lambda^3-27Q)^{3/2}},
\een
where $C$ is a non-zero constant. The precise value of $C$ is
irrelevant, but for the sake of completeness, let us compute it. 
Using the connection formulas in Section \ref{sec:conn-formula}, we
can compute the leading order term of the Wronskian near $x=\infty$
\ben
\operatorname{det}\,
\begin{bmatrix}
z_1 & x\partial_x z_1 & (x\partial_x )^2 z_1 \\
z_2 & x\partial_x z_2 & (x\partial_x )^2 z_2 \\
z_3 & x\partial_x z_3 & (x\partial_x )^2 z_3 
\end{bmatrix} = -\frac{16}{3\sqrt{6}}\, \sqrt{-1} \, x^{-1/2}
+ O(x^{-3/2}).
\een
We get $C=-\frac{16}{3\sqrt{6}}$. Finally, for the Jacobian
determinant we get
\ben
\operatorname{det}\, \frac{D(Z_1,Z_2,Z_3)}{D(\lambda,Q,t)} =
8\sqrt{-2}\, (\lambda^3-27Q)^{-1/2}. 
\een
b) This is an elementary consequence of Ramanujan's differential
equations. 
\qed

\begin{proposition}\label{prop:coeff-anal}
The coefficients
\ben
Q_n,\lambda_n\ \in \ \CC[E_2(\tau_{12}),E_4(\tau_{12}),E_6(\tau_{12})],
\een
i.e., they are quasi-modular forms with respect to
$\tau_{12}\in \HH$. 
\end{proposition}
\proof
Using Proposition \ref{prop:hol-lift} and Lemma \ref{le:Jac-det}, b) we
get that the map $Z^{\rm aux}:\U\to \CC^3$ induces an isomorphism
between an open neighborhood in $\HH\times \CC^*\times \CC$ of 
$$
(\HH\times \CC^*)':=\{(\tau,x)\in \HH\times \CC^* \ |\ E_6(\tau)\neq
0\} 
$$
and an open neighborhood in $\CC^3$ of
\ben
\Omega_{\rm small}' = \{z^*\in \Omega_{\rm small}\ |\  E_6(-z^*_2/(2z_3^*))\neq 0\}.
\een
In particular, the coordinates $(\tau,x,s)$ and $(z_1,z_2,z_3)$ of the
two neighborhoods are biholomorphic. Therefore
\ben
\lambda:=2(2\pi/x)^2 E_4(\tau),\quad 
Q:=\frac{8}{27}(2\pi/x)^6 (E_4(\tau)^3-E_6(\tau)^2),\quad
t:=-\frac{1}{32} (z_2^2-4z_1z_3)
\een
define functions that are 
holomorphic in an open neighborhood in $\CC^3$ of $\Omega_{\rm small}'$. Note that 
by definition the functions $(Q,t,\lambda)$ give an inversion of the
period map. More precisely, if $(Q^*,t^*,\lambda^*)\in X$ with $t^*=0$ and
$z^*=Z(Q^*,t^*,\lambda^*)\in \Omega_{\rm small}$ is a value of the period map
(depending on the choice of a reference path), then in a neighborhood
of $z^*\in \CC^3$ the functions
$(Q,t,\lambda)$ coincide with the unique solution to the equations
\ben
Z_i(Q,t,\lambda) = z_i,\quad 1\leq i\leq 3,
\een
where the branch of $Z_i$ is fixed by $Z_i(Q^*,t^*,\lambda^*) =
z_i^*$.

Clearly $t$ is a holomorphic function on $\CC^3$. We claim that
$\lambda$ and $Q$ extend to holomorphic  functions defined in a
neighborhood of $\Omega_{\rm small}$ in $\CC^3$. The  statement is local, so let
$(Q^*,t^*,\lambda^*)$ with $t^*=0$ be a point on the discriminant and
let $z^*=Z(Q^*,t^*,\lambda^*)\in \Omega_{\rm small}$ be a value of the period
map. Let $\lambda=u(Q,t)$ be the local equation of the discriminant at
the point $(Q^*,t^*,\lambda^*)$. Locally, the components of the period
map can be written as  
\ben
Z_i(Q,t,\lambda) = \frac{1}{2}(E_i|\alpha) Z_\alpha(Q,t,\lambda) +
Z_i^{\rm inv}(Q,t,\lambda)
\een
where $\alpha\in H$ is a vector whose local monodromy around the
discriminant is given by $\alpha\mapsto -\alpha$, $Z_\alpha:=\langle
Z,\alpha\rangle$, and $Z_i^{\rm inv}$ corresponds to the invariant
part of $E_i$, i.e., 
$$
Z_i^{\rm inv}=\langle Z, E_i-(E_i|\alpha)\alpha/2\rangle.
$$ 
On the other hand
\ben
Z_\alpha(Q,t,\lambda) = (\lambda-u(Q,t))^{1/2} \widetilde{Z}_\alpha(Q,t,\lambda)
\een
where $ \widetilde{Z}_\alpha(Q,t,\lambda)$ is holomorphic in a
neighborhood of $(Q^*,t^*,\lambda^*)$ and
$\widetilde{Z}_\alpha(Q^*,t^*,\lambda^*)\neq 0$. Let us choose $i$ such
that $(E_i|\alpha)\neq 0$. Then we get 
\ben
\lambda = u(Q,t) + \mu^2,\quad \mu:= \frac{2}{(E_i|\alpha)}
(Z_i-Z_i^{\rm inv})/\widetilde{Z}_\alpha.
\een
The above equation defines a branched double covering of a
neighborhood of $(Q^*,t^*,\lambda^*)$ and $(Q,t,\mu)$ is a holomorphic
coordinate system  on the double cover. We claim that the local lift
of the period map is an isomorphism. Indeed, the local lift is
a single valued analytic map, because 
\ben
Z_j = \frac{1}{2} (E_j|\alpha) \mu\, \widetilde{Z}_\alpha + Z_j^{\rm
  inv},\quad 1\leq j\leq 3.
\een
We have to check that the corresponding Jacobian determinant does not
vanish at the point $(Q^*,t^*,\mu^*)=(Q^*,0,0).$  Using
\ben
\frac{D(Q,t,\lambda)}{D(Q,t,\mu)}   = 2\mu= 2 (\lambda-u(Q,t))^{1/2},
\een
the chain rule, and Lemma \ref{le:Jac-det}, a) we get 
\ben
\frac{D(Z_1,Z_2,Z_3)}{D(Q,t,\mu)}
(Q^*,0,0)=\frac{2C}{\sqrt{(1-\zeta)(1-\zeta^2)}}\,
\frac{1}{u(Q^*,0)}\neq 0,
\een 
where $C$ is a non-zero constant, $\zeta=e^{2\pi\sqrt{-1}/3}$, and we used that 
\ben
\lambda^3-27Q^* = \prod_{a=0}^2 (\lambda - \zeta^a u(Q^*,0)).
\een
Therefore $Q,t,$ and $\mu$ are
holomorphic in a neighborhood of $z^*\in \CC^3$, which implies that
$Q,t,$ and $\lambda$ are also holomorphic. 

To complete the proof of the proposition we note that 
\ben
\tau_{12} = -\frac{z_2}{2z_3},\quad (\tau_1-\tau_2)^2 = (z_2/z_3)^2
-4(z_1/z_3),\quad y=z_3.
\een
Therefore, $Q,t,$ and $\lambda$ must be holomorphic functions in
$y,\tau_{12},(\tau_1-\tau_2)^2$. In particular $Q_n(\tau_{12})$ and
$\lambda_n(\tau_{12})$ must be holomorphic in $\tau_{12}$, so the
corresponding polynomial expressions in $E_2(\tau_{12})$,
$E_4(\tau_{12})$, and $E_6(\tau_{12})^{\pm 1}$ could not 
have negative powers of $E_6(\tau_{12})$.  
\qed

\subsection{The ring of modular functions}
Note that the ring
$\Gamma(\Omega_{\rm small},\O_{\CC^3})$ is equipped with the action of 
the monodromy group $W$. Let $\Gamma(\Omega_{\rm small},\O_{\CC^3})^W$ be the
ring of $W$-invariant functions. We introduce a subring  of
$W$-invariant functions as follows. Let $\CC[Q,\lambda]\{t\}$ be the ring of power series in $t$
whose coefficients depend polynomially on $Q$ and $\lambda$ and such
that for every $(Q,\lambda)\in \CC^*\times \CC$ the radius of
convergence is non-zero. Then we define
\ben
\mathcal{M}(\Omega,W)= 
\{f\in \Gamma(\Omega_{\rm small},\O_{\CC^3})^W\ |\ f\circ Z \in \CC[Q,\lambda]\{t\}\}.
\een
Proposition \ref{prop:coeff-anal} implies that the tautological
map 
\ben
\mathcal{M}(\Omega,W)\to \CC[Q,\lambda]\{t\},\quad f\mapsto f\circ Z
\een
is an isomorphism. Although the invariant functions corresponding to
$Q$ and $\lambda$ can be find recursively, it will be nice to have a
more intrinsic characterization. Unfortunately we could not achieve
this goal. On the other hand we have managed to find explicitly  invariant functions
$E_4^{(2)},\Delta^{(2)}$ that generate $\mathcal{M}(\Omega,W)$ in the
following sense. If $f\in \mathcal{M}(\Omega,W)$ then 
\ben
f = \sum_{n=0}^\infty c_n(E_4^{(2)},\Delta^{(2)}) t^n,
\een
where $c_n\in \CC[E_4^{(2)},\Delta^{(2)}]$ are some polynomials. The
above equality should be interpreted as equality between formal power
series in $t$. 

In order to find such functions $E_4^{(2)}$ and $\Delta^{(2)}$ it is
enough to construct two $W\times \mu_2$-invariant holomorphic  functions in $\HH^2\times \CC^*$ whose
restrictions to $\HH\times \CC^*$ coincide with
\beq\label{restr-Q-la}
\frac{8}{27}\, (2\pi/z)^6\, (E_4(\tau)^3-E_6(\tau)^2)\quad
\mbox{and}
\quad
2\,(2\pi/z)^2\, E_4(\tau).
\eeq
This could be done easily using the Jacobi theta constants
\ben
\theta_{ab}(0,\tau) =\sum_{n\in \ZZ} \exp\Big(
\pi\sqrt{-1}\left( (n+a/2)^2\tau + (n+a/2)b )\right)\Big),
\quad ab=00,\ 01,\ 10.
\een
For the reader's convenience we have recorded in Table
\ref{tab:tr-rules}  the transformation rules for the theta constants
under the two modular 
transformations $\tau\mapsto \tau+1$ and $\tau\mapsto -1/\tau$. For
more details we refer to \cite{Mu}. 
It is easy to check that 
\ben
E_4(\tau) & = & \frac{1}{2}
(\theta_{00}(\tau)^8+\theta_{10}(\tau)^8+\theta_{01}(\tau)^8)\\
E_4(\tau)^3-E_6(\tau)^2 &  = &  \frac{27}{4} (\theta_{00}(\tau)\theta_{10}(\tau)\theta_{01}(\tau))^8.
\een
\begin{table}
\centering
\caption{Transformation rules}
\label{tab:tr-rules}
\begin{tabular}{l|c|c}
& $\tau\mapsto \tau+1$ & $\tau\mapsto -1/\tau$ \\
\hline
& & \\
$\theta_{00}(\tau)$ & $\theta_{01}(\tau)$ & $(-\mathbf{i}\tau)^{1/2}
\theta_{00}(\tau)$ \\
\hline
& & \\
$\theta_{01}(\tau)$ & $\theta_{00}(\tau)$ & $(-\mathbf{i}\tau)^{1/2}
\theta_{10}(\tau)$ \\
\hline
& & \\
$\theta_{10}(\tau)$ & $e^{2\pi\mathbf{i}/8}\theta_{10}(\tau)$ & $(-\mathbf{i}\tau)^{1/2}
\theta_{01}(\tau)$
\end{tabular}
\end{table}
Let us define 
\ben
E^{(2)}_4(\tau_1,\tau_2,x) :=  (2\pi/x)^2 \, \sum_{ab\in \{00,01,10\}}
\theta_{ab}(\tau_1)^4 \theta_{ab}(\tau_2)^4
\een
and
\ben
\Delta^{(2)}(\tau_1,\tau_2,x) :=  2(2\pi/x)^6 \, \prod_{ab\in \{00,01,10\}}
\theta_{ab}(\tau_1)^4 \theta_{ab}(\tau_2)^4.
\een
It is straightforward to check that $E^{(2)}_4$ and $\Delta^{(2)}$ are
$W\times \mu_2$-invariant holomorphic functions on $\HH^2\times
\CC^*$, so they define $W$-invariant analytic functions on the domain
$\HH^2\times \CC^*/\mu_2$. In particular  $E^{(2)}_4,\Delta^{(2)}\in
\mathcal{M}(\Omega,W)$. 

\appendix\label{sec:appendix}

\section{Genus-0 constraints}
Suppose that we are in the settings of Sections \ref{sec:2nd-str-conn}
and \ref{sec:pv}. Following \cite{G3} we will assume in addition that the Frobenius
structure arises from a set of gravitational descendants. The latter are organized
into a generating function $\F^{(0)}(\tt)$, where
$\tt=(t_{k,i})_{k=0,1,2,\dots}^{i=1,2,\dots,N}$ is a set of formal
variables,  satisfying the following 3 axioms.
\begin{enumerate}
\item[(DE)]
Dilaton Equation:
\ben
\frac{\partial \F^{(0)}}{\partial t_{1,1}} = \sum_{k=0}^\infty \sum_{i=1}^N
t_{k,i}\frac{\partial \F^{(0)}}{\partial t_{k,i}} - 2 \F^{(0)}.
\een
\item[(SE)]
String Equation:
\ben
\frac{\partial \F^{(0)}}{\partial t_{0,1}} = \frac{1}{2}(t_0,t_0) + \sum_{k=0}^\infty \sum_{i=1}^N
t_{k+1,i}\frac{\partial \F^{(0)}}{\partial t_{k,i}},
\een
where $t_0=\sum_{i=1}^N t_{0,i}\phi_i$.
\item[(TRR)]
Topological Recursion Relations:
\ben
\frac{\partial^3 \F^{(0)}}{
  \partial t_{k_1+1,i_1} \partial t_{k_2,i_2}\partial t_{k_3,i_3}} =
\sum_{j_1,j_2=1}^N
\frac{\partial^2 \F^{(0)}} {\partial t_{k_1,i_1} \partial t_{0,j_1}}
g^{j_1 j_2}
\frac{\partial^3 \F^{(0)}}{
  \partial t_{0,j_2} \partial t_{k_2,i_2}\partial t_{k_3,i_3}},
\een
where $g_{ij}=(\phi_i,\phi_j)$ is the matrix of the Frobenius pairing
and $g^{ij}$ are the entries of the inverse matrix.
\end{enumerate}

\noindent
Finally, we are going to consider only the case when the intersection pairing is non-degenerate.

\subsection{The Heisenberg vertex operator algebra}
In this section we recall the main construction from
\cite{BM}. Although the work in \cite{BM} is in the settings of simple
singularities the generalisation to an arbitrary semi-simple Frobenius
manifold is straightforward (see \cite{Mi2}, Section 5).
Let us equip the vector space $H[s,s^{-1}]\oplus \CC$ with the
structure of a Heisenberg Lie algebra such that
\ben
[\alpha s^m,\beta s^n]=m\delta_{m,-n}(\alpha|\beta),\quad
[\alpha s^m, c]=0,\quad
\alpha,\beta\in H,\quad m,n\in \ZZ,\quad c\in \CC.
\een
Let $\mathfrak{F}:=\operatorname{Sym}(H[s^{-1}]s^{-1})$ be the corresponding
Fock space, i.e., the unique irreducible highest weight representation
with highest weight vector $1$ defined by
\ben
(\alpha s^n ) 1 = 0,\quad n\geq 0. 
\een
Following \cite{BM} we introduce the W-algebra $\mathcal{W}\subset \mathfrak{F}$
as the kernel of all screening operators
\ben
\mathcal{W}=\{w\in \mathfrak{F}\ |\ e^{\alpha}_{(0)} w=0\ \forall
\mbox{ reflection vectors }\alpha\}
\een
where the linear operators $e^{\alpha}_{(n)}$ are defined as the
Fourier coefficients of the following vertex operator:
\ben
\Gamma^\alpha(\zeta):=\zeta^{-|\alpha|^2/2}
e^{\sum_{n>0} (\alpha  s^{-n}) \frac{\zeta^n}{n} }
e^{\sum_{n<0} (\alpha  s^{-n}) \frac{\zeta^n}{n} }=:
\sum_{n\in \ZZ} e^{\alpha}_{(n)} \zeta^{-n-1} ,
\een
where $|\alpha|^2=(\alpha|\alpha)$. The vector space $\mathfrak{F}$ is equipped
with the structure of a Vertex Operator Algebra (VOA) such that the
{\em state field correspondence} 
\ben
Y(a,\zeta)=\sum_{n\in \ZZ} a_{(n)} \zeta^{-n-1},
\quad a\in \mathfrak{F},
\quad a_{(n)}\in \operatorname{End}(\mathfrak{F})
\een
is defined by
\ben
Y(\alpha s^{-1},\zeta) := \sum_{n\in \ZZ} (\alpha s^n) \zeta^{-n-1}
\een
and the following {\em operator product expansion} formula holds:
\ben
Y(a_{(n)}b,\zeta) = \operatorname{Res}_{\zeta_1=\zeta} d\zeta_1
(\zeta_1-\zeta)^n Y(a,\zeta_1)Y(b,\zeta).
\een
Let us fix a basis $\{\phi_i\}_{i=1}^N$ of $H$ and define the
following space of formal power series
\ben
\mathfrak{M}=\CC_\hbar[\![q_0,q_1+\mathbf{1},q_2,\dots]\!]
\een
where $\CC_\hbar:=\CC(\!(\hbar^{1/2})\!)$ and
$q_k=(q_{k,1},\dots,q_{k,N})$ is a sequence of formal vector 
variables. The vector space $\mathfrak{M}$ is equipped with the structure of a
twisted VOA module structure over the VOA $\mathfrak{F}$ in the following
way. Following Givental we equip the vector space $H[z,z^{-1}]$ with a
symplectic form
\beq\label{sympl-form}
\Omega(f,g) = \operatorname{Res}_{z=0} dz (f(-z),g(z)).
\eeq
Let $H[z,z^{-1}]\oplus \CC$ be the Heisenberg Lie algebra with bracket
\ben
[f(z),g(z)]=\Omega(f,g),\quad f,g \in H[z,z^{-1}].
\een
We define a representation of this Heisenberg Lie algebra on $\mathfrak{M}$ as
follows
\ben
\phi_i z^k\mapsto (\phi_i z^k)\sphat:= - \sqrt{\hbar}
\frac{\partial}{\partial q_{k,i}},\quad
\phi^i(-z)^{-k-1}\mapsto (\phi^i(-z)^{-k-1})\sphat:=
q_{k,i}/\sqrt{\hbar}, 
\een
where $\{\phi^i\}_{i=1}^N$ is a basis of $H$ dual to
$\{\phi_i\}_{i=1}^N$ with respect to the Frobenius pairing. Put
\ben
Y^\mathfrak{M}(\alpha s^{-1},\lambda) :=
(\partial_\lambda \widetilde{\mathbf{f}}_\alpha(\lambda,z))\sphat,
\een
where
\ben
\widetilde{\mathbf{f}}_\alpha(\lambda,z):=\sum_{n\in \ZZ}
\widetilde{I}^{(n)}_\alpha(\lambda) (-z)^n.
\een
We define $Y^\mathfrak{M}(a,\lambda)$ for all $a\in \mathfrak{F}$ in such a way that the
following operator product expansion formula holds
\ben
Y^\mathfrak{M}(a_{(n)}b,\lambda) = \operatorname{Res}_{\lambda_1=\lambda}
d\lambda_1 (\lambda_1-\lambda)^n Y^\mathfrak{M}(a,\lambda_1)Y^\mathfrak{M}(b,\lambda). 
\een
\subsection{Genus-0 reduction}
The total descendant potential of a semi-simple Frobenius manifold $M$
is a formal power series of the type
\ben
\D(\hbar,\qq):=\exp\Big(\sum_{g=0}^\infty \F^{(g)}(\qq) \hbar^{g-1}\Big)
\een
where $\qq=(q_0,q_1,\dots)$ is the sequence of formal variables from
above and $\F^{(g)}(\qq)$ ($g\geq 0$) is the so called {\em genus-g
  descendant} potential. By definition
\ben
\F^{(g)}(\qq)\ \in\ \O_M(M-\mathcal{K})[\![q_1+\mathbf{1},q_2,\dots]\!],
\een
where $\mathcal{K}\subset M$ is the analytic hypersurface of all non-semisimple
points and $q_0=(q_{0,1},\dots,q_{0,N}) $ should be identified with
a flat coordinate system on $M$. For precise definitions and more
details we refer to \cite{G1,G2}. Let us fix a flat coordinate system
defined in a neighborhood of a semi-simple point $t^\circ\in M$, s.t.,
all coordinates $t^\circ_i=0$. Then by taking the Taylor series
expansion at $q_0=0$ we identify $\D(\hbar,\qq)$ with an element in
the Fock space $\mathfrak{M}$

First of all note that if $w\in \mathfrak{F}$ is monodromy invariant
then $Y^\mathfrak{M}(w,\lambda)\D$ is a formal power series in $\qq$ whose coefficients are formal
Laurent series in $\hbar$ whose coefficients  are Laurent series in
$\lambda^{-1}$.  The proof of Theorem 1.1 in \cite{BM}, Section 8
is straightforward to generalize to the current settings. We get that
if $w\in \W$, then
$Y^\mathfrak{M}(w,\lambda)\D$ is a formal power series in $\qq$, formal Laurent
series in $\hbar$ whose coefficients are polynomials in $\lambda$. For
brevity in the latter case we say that $Y^\mathfrak{M}(w,\lambda)\D$ is {\em regular}
in $\lambda$. Note that the regularity is equivalent to a sequence of
differential operator constraints for $\D$. Indeed expanding
$Y^\mathfrak{M}(w,\lambda)$ into a Laurent series in $\lambda^{-1}$ yields a
sequence of differential operators acting on $\mathfrak{M}$. The regularity means
that the differential operators in front of negative powers of
$\lambda$ annihilate $\D$.

Let us compute the leading order term in the Laurent series expansion
in $\hbar$ of $Y^\mathfrak{M}(w,\lambda)\D/\D$. It is convenient to embed $H\to
\mathfrak{F}$ via $\alpha\mapsto \alpha s^{-1}$. Let us fix a basis
$\alpha_1,\dots,\alpha_N$ of $H$ and denote by
\ben
\partial^k\alpha_i := (-\partial_s)^k (\alpha_i s^{-1}) = k! \alpha_i
s^{-k-1}. 
\een
Then the Fock space
\ben
\mathfrak{F}= \CC[\partial^k\alpha_i\ |\ 1\leq i\leq N,\ k\geq 0].
\een
We will refer to $\partial^k\alpha_i$ for $k>0$ as {\em jet
  variables}. 
For a given monomial
\ben
(\partial^{k_1}\alpha_1 )^{l_1}\cdots (\partial^{k_N}\alpha_N )^{l_N},
\quad
k_i,l_i\in \ZZ_{\geq 0},
\een
we define its weight to be $l_1(k_1+1)+\cdots +l_N(k_N+1)$ and its degree to
be $l_1+\cdots +l_N$. 
Let us introduce a grading of $\mathfrak{F}$ such that the homogeneous
piece $\mathfrak{F}_n$ is spanned by monomials of weight $n$. The
W-algebra $\W$ is a graded subspace of $\mathfrak{F}$, i.e., if $w\in
\W$ and $w=\sum_{n=0}^\infty w_n$ is a decomposition into homogeneous
components $w_n\in \mathfrak{F}_n$ then $w_n\in \W$.

Suppose now that
$w\in \W$ is homogeneous, i.e., $w\in \mathfrak{F}_n$ for some $n$.
Note that $w=w^{(n)}+w^{(n-1)}+\cdots + w^{(1)}$, where $w^{(i)}$ is a
linear combination of monomials of degree $i$. Moreover the term
$w^{(n)}=f(\alpha_1,\dots,\alpha_N)\in \CC[\alpha_1,\dots,\alpha_N]$. 
Recalling the operator product expansion formula we get 
\ben
Y^\mathfrak{M}(\partial^k \alpha,\lambda) = \partial_\lambda^k
Y^\mathfrak{M}(\alpha,\lambda). 
\een
Therefore
\ben
Y^\mathfrak{M}(\partial^k \alpha,\lambda) \D(\hbar,\qq) =
\left(\left.\hbar^{-1/2} \partial_\lambda^k
    \phi_\alpha(\lambda,\qq,\pp)
  \right|_{\pp=\partial_\qq\F^{(0)}} +
  \cdots \right)\D(\hbar,\qq),
\een
where the dots stand for terms that involve higher powers of
$\hbar$, $\pp=(p_0,p_1,\dots)$ is a sequence of vector variables 
$p_k=(p_{k,1},\dots,p_{k,N})$, the substitution $\pp=\partial_\qq\F^{(0)}$
means $p_{k,i}=\frac{\partial \F^{(0)}}{\partial q_{k,i}}$, and  
\beq\label{def-fia}
\phi_\alpha(\lambda,\qq,\pp) := \sum_{k=0}^\infty
\Big(
(\widetilde{I}_\alpha^{(-k)}(\lambda),\phi_i) q_{k,i} + (-1)^{k+1}
(\widetilde{I}_\alpha^{(k+1)}(\lambda),\phi^i) p_{k,i} \Big). 
\eeq
Note that $Y^\mathfrak{M}(w,\lambda)\D/\D$ is a Laurent series in
$\hbar^{1/2}$ whose leading order term is
\beq\label{lot}
\hbar^{-n/2}
\left.
f(\phi_{\alpha_1}(\lambda,\qq,\pp),\dots,
\phi_{\alpha_N}(\lambda,\qq,\pp))
\right|_{\pp=\partial_\qq\F^{(0)}}.
\eeq
We get that if $w\in \W$ is homogeneous of weight $n$ then the
expression \eqref{lot} is regular in $\lambda$, where $f=w^{(n)}$ is
the degree $n$ part of $w$. By definition the regularity condition
means that in the Laurent series expansion in
$\lambda^{-1}$ all coefficients in front of negative powers of
$\lambda$ must vanish. On the other hand, every $\alpha_i\in H$ is
a linear function on $H^*$. 
Therefore the polynomial $f$ defines a holomorphic function on $H^*$.
The main result in this Appendix can be stated as follows. 
\begin{proposition}\label{pr:mod-W}
If $f\in \CC[\alpha_1,\dots,\alpha_N]$ is a $W$-invariant polynomial,
then $f\in \M(\Omega,W)$ and the expression \eqref{lot} is regular. 
\end{proposition}
\begin{remark}
We expect that Proposition \ref{pr:mod-W} can be generalized in the
following way:  if $f\in \O_{H^*}(\Omega)$ is a
$W$-invariant holomorphic function then the expression \eqref{lot} makes sense and
it is regular in $\lambda$ if and only if $f\in \M(\Omega,W)$. 
\end{remark}
\begin{remark}
  If the monodromy group $W$ is not finite (which is almost always the
  case), then the ring of modular functions $\M(\Omega,W)$ contains
  very few polynomials. Therefore the W-algebra defined in \cite{BM}
  is not big enough to characterize the total descendant
  potential. Nevertheless there is still some hope that the
  construction from \cite{BM} can be generalized, because we can replace
  the Fock space $\mathfrak{F}$  with a larger one, such as 
  \ben
  \O_{H^*}(\Omega)[\partial^k \alpha_i\ :\ k>0,\ 1\leq i\leq N]. 
  \een
  It will be interesting to find out if the modular functions for
  $\PP^2$ found in this paper can be extended to W-constraints for the
  total descendant potential of $\PP^2$.
\end{remark}

\subsection{Givental's symplectic space formalism}\label{sec:lagr_cone}
The proof of Proposition \ref{pr:mod-W} relies on the properties of a
certain Lagrangian cone (see \cite{G3}). Let us recall the necessary
background.

Put $\H=H(\!(z^{-1})\!)$ and let us define a symplectic structure on
$\H$ via formula \eqref{sympl-form}. We have $\H=\H_+\oplus \H_-$,
where $\H_+:=H[z]$ and $\H_-:=H[\![z^{-1}]\!]z^{-1}$ are Lagrangian
subspaces. The formal variables $\qq$ are identified with coordinates
on $\H_+$ via
\ben
\qq\mapsto \qq(z):=\sum_{k=0}^\infty \sum_{i=1}^N q_{k,i} \phi_i z^k.
\een
The linear structure on $\H_+$ gives a natural identification of each
tangent space
\ben
T_\qq\H_+\cong \H_+,\quad \frac{\partial}{\partial q_{k,i}}\mapsto
\phi_i z^k.
\een
Using the symplectic form $\Omega$ we identify
\ben
\H_-\cong (\H_+)^*,\quad v\mapsto \Omega(v,\cdot).
\een
Therefore the cotangent bundle $T^*\H_+\cong \H$. Note that under the
identification $T_\qq\H_+\cong \H_+$ the 1-form
$dq_{k,i}=\Omega(\phi^i (-z)^{-k-1},\cdot)$. The axioms of
gravitational descendants can be reformulated in the following
geometric way. Let us identify $\F^{(0)}(\tt)$ with a function on
$\H_+$ via the substitution $\tt(z)=\qq(z)+ z$. This change of
variables is known as the {\em dilaton shift}. Let
\ben
J(\qq,z):= \qq(z) +\sum_{k=0}^\infty
\sum_{i=1}^N\frac{\partial\F^{(0)}}{\partial q_{k,i}} (\qq) \phi^i(-z)^{-k-1}.
\een
Intuitively as $\qq$ varies in $\H_+$ the values of $J(\qq,z)$ defines
a subset $\L\subset \H$, which 
under the isomorphism $\H\cong T^*\H_+$ coincides
with the graph of the differential $d_\qq\F^{(0)}$. The 3 axioms
DE, SE, and TRR are equivalent to the following properties of $\L$:
\begin{enumerate}
\item[(1)] $\L$ is a quadratic cone
\item[(2)] If $\ff\in \L$ is a smooth point then the tangent
space $L=T_\ff\L$ is Lagrangian, $zL=L\cap \L$, and if $\mathbf{g}\in
zL$, then $T_\mathbf{g}\L=L$.
\end{enumerate}
Formally the properties of $J(\qq,z)$ can be stated as follows. Let us
recall that if the Frobenius structure comes from gravitational
descendants then there is a natural choice of a calibration
\ben
(S(t,z)\phi_i,\phi_j):=(\phi_i,\phi_j) +\sum_{k=0}^\infty z^{-k-1}\, 
\left.
  \frac{\partial^2 \F^{(0)}}{\partial t_{k,i}\partial t_{0,j}}
\right|_{\tt(z)=t}.
\een
There exists
\ben
\tau(\tt)\in H[\![\tt]\!],\quad v(\tt,z)\in H[\![\tt,z]\!] z 
\een
such that
\ben
J(\qq,z) = S(\tau(\tt),z)^{-1}  v(\tt,z). 
\een
Note that $J(\qq,z)=-z+J(\tt,z)$ and that $S_1(t) 1= t$. Multiplying
both sides of the above identify by  $S(\tau(\tt),z)$ and comparing
the coefficients in front of $z^0$ yields
\ben
t_0 + S_1(\tau) t_1+S_2(\tau) t_2 +\cdots = \tau.
\een
This equation allows us to solve for $\tau$ in terms of
$\tt$. Finally $v(\tt,z)= [S(\tau(\tt),z) (\tt(z)-z)]_+$, where
$[\cdot ]_+: \H\to \H_+$ is the projection. In other words,
the formal series $J(\qq,z)$ is uniquely determined from the
calibration $S(t,z)$.
\subsection{Proof of Proposition \ref{pr:mod-W}}
Recalling the definition \eqref{def-fia} we get that 
\ben
\left.\phi_\alpha(\lambda,\qq,\pp)
\right|_{\pp=\partial_\qq\F^{(0)}}=
\partial_\lambda \Omega(\widetilde{\ff}_\alpha(\lambda,z), J(\qq,z)).
\een
There exists $\tau(\qq)$ and $v(\qq,z)$ such that
$J(\qq,z)=S(\tau(\qq),z)^{-1}v(\qq,z)$ (see Section
\ref{sec:lagr_cone}). The calibration $S$ is a symplectic
transformation. Therefore
\ben
\left.\phi_\alpha(\lambda,\qq,\pp)
\right|_{\pp=\partial_\qq\F^{(0)}} =
\partial_\lambda \Omega(\ff_\alpha(\tau(\qq),\lambda,z),v(\qq,z))=
\Omega(\ff_\alpha(\tau(\qq),\lambda,z),v(\qq,z) z^{-1}),
\een
where we used that
\ben
\ff_\alpha(t,\lambda,z):=\sum_{n\in \ZZ}
I^{(n)}_\alpha(t,\lambda)(-z)^n=
S(t,z) \widetilde{\ff}_\alpha(\lambda,z).
\een
Let $\mathfrak{m}_k$ ($k\geq 1$) be the ideal in
$\O_M(M)[\![q_1+\mathbf{1},q_2,\dots,]\!]$ generated by
$(q_k+\mathbf{1}\delta_{k,1},q_{k+1},\dots)$. Note that $\tau(\qq)\in
q_0+\mathfrak{m}_1$ and
\ben
v(\qq,z)z^{-1} = -\mathbf{1} + \sum_{k=0}^\infty v_k (\qq) z^k,
\quad v_k\in \mathfrak{m}_k.
\een
We have
\ben
\left.\phi_\alpha(\lambda,\qq,\pp)
\right|_{\pp=\partial_\qq\F^{(0)}} =
\left.\Big( - (I^{(-1)}_\alpha(t,\lambda) ,\mathbf{1})+
  \sum_{k=0}^\infty (I^{(-k-1)}_\alpha(t,\lambda),v_k)
  \Big)\right|_{t=\tau(\qq)}
\een
Suppose now that $f\in \CC[\alpha_1,\dots,\alpha_N]$ is a homogeneous
polynomial of degree $n$. Note that the expression \eqref{lot} is
obtained from the polynomial $f(\alpha_1,\dots,\alpha_N)$ via two
substitutions: first
\beq\label{subst-1}
\alpha_i\mapsto - (I^{(-1)}_\alpha(t,\lambda) ,\mathbf{1})+
\sum_{k=0}^\infty (I^{(-k-1)}_\alpha(t,\lambda),v_k),
\quad 1\leq i\leq N,
\eeq
and second $t\mapsto \tau(q)$. We claim that the first substitution
yields an expression regular in $\lambda$. This would complete the proof of
the proposition, because the second substitution preserves the
regularity property. Suppose that we make the substitution
\eqref{subst-1} in $f(\alpha_1,\dots,\alpha_N)$. Let us fix $t\in
M$. Then the resulting expression is a formal 
power series in $q_1+\mathbf{1}, q_2,\dots$ whose coefficients are
polynomial expressions in $(I^{(-k)}_{\alpha_i}(t,\lambda),\phi_j)$
with $k>0$. Since $f$ is a $W$-invariant polynomial, each coefficient
is a single-valued holomorphic function on
$\CC\setminus{\{u_1,\dots,u_N\}}$, where $u_i$ are the eigenvalues of
the operator $E\bullet_t$. Since the periods have finite order pole at
$\lambda=\infty$ the regularity condition will be established if we
manage to prove that the coefficients are holomorphic at
$\lambda=u_i$. This however follows from the Riemann's extension
theorem, because the period vectors
$(I^{(-k)}_{\alpha_i}(t,\lambda),\phi_j)$ for $k>0$ are bounded in a
neighborhood of $\lambda=u_i$. Note that if $v_k=0$ for $k\geq 0$ then
this argument proves that $f\in \M(\Omega,W)$.
\qed

\bibliographystyle{amsalpha}

\end{document}